\documentclass{article}

\usepackage{arxiv}
\usepackage{subfigure}
\usepackage[utf8]{inputenc} 
\usepackage[T1]{fontenc}    
\usepackage{url}            
\usepackage{booktabs}       
\usepackage{amsfonts}       
\usepackage{nicefrac}       
\usepackage{microtype}      
\usepackage{lipsum}

\usepackage{fixltx2e}
\usepackage{amsmath}
\usepackage[normalem]{ulem}
\usepackage{amssymb}
\usepackage{amsthm}
\usepackage{color}
\usepackage{pifont}
\newcommand{\RN}[1]{%
	\textup{\uppercase\expandafter{\romannumeral#1}}%
}
\usepackage{listings}
\usepackage{color} 
\definecolor{mygreen}{RGB}{28,172,0} 
\definecolor{mylilas}{RGB}{170,55,241}

\newtheorem{assumption}{Assumption}[section]
\newtheorem{example}{Example}[section]

\newcommand{\quoteIt}[1]{{``#1''}}
\renewcommand{\Re}{{\mathbb{R}}}
\newcommand{\1}{{\mathbf{1}}}
\newcommand{\mymbox}[1]{\mbox{\scriptsize #1}}
\newcommand{\U}{\mathcal{U}}

\newcommand{\myzeta}{\zeta}
\newcommand{\price}{p}
\DeclareMathOperator*{\esssup}{ess\,sup}

\newcommand{\removed}[1]{{}}
\newcommand{\ErickCommentsForFuture}[1]{{}}

\newcommand{\SaeedComments}[1]{}

\newcommand{\Tau}{\mathcal{T}}
\newcommand{\modified}[1]{{ #1}}

\newcommand{\saeedmodified}[1]{{\color{magenta} #1}}

\newcommand{\erickdelete}[1]{{}}
\newcommand{\saeedZeta}[1]{\zeta}

\theoremstyle{plain}
\newtheorem{theorem}{Theorem}[section]
\newtheorem{corollary}[theorem]{Corollary}
\newtheorem{lemma}[theorem]{Lemma}
\newtheorem{proposition}[theorem]{Proposition}
\theoremstyle{definition}
\newtheorem{definition}{Definition}
\theoremstyle{remark}
\newtheorem{remark}{Remark}

\usepackage{graphicx}

\title{Equal Risk Pricing and Hedging of Financial Derivatives with Convex Risk Measures}

\author{
  Saeed Marzban \\
  HEC Montr\'eal, Montr\'eal,\\
  H3T 2A7, Canada\\
  \texttt{saeed.marzban@hec.ca} \\
   \And
 Erick Delage \\
  HEC Montr\'eal, Montr\'eal,\\
  H3T 2A7, Canada\\
  \texttt{erick.delage@hec.ca} \\
  \And
  Jonathan Yumeng Li \\
  Telfer School of Management, University of Ottawa,\\
  Ottawa, Ontario K1N 6N5, Canada\\
  \texttt{jonathan.li@telfer.uottawa.ca} \\
}

\begin{document}
\maketitle

\begin{abstract}
In this paper, we consider the problem of equal risk pricing and hedging in which the fair price of an option is the price that exposes both sides of the contract to the same level of risk. Focusing for the first time on the context where risk is measured according to convex risk measures, we establish that the problem reduces to solving independently the writer and the buyer's hedging problem with zero initial capital. By further imposing that the risk measures decompose in a way that satisfies a Markovian property, we provide dynamic programming equations that can be used to solve the hedging problems for both the case of European and American options. All of our results are general enough to accommodate situations where the risk is measured according to a worst-case risk measure as is typically done in robust optimization. 
Our numerical study illustrates the advantages of equal risk pricing over schemes that only account for a single party, pricing based on quadratic hedging  (i.e. $\epsilon$-arbitrage pricing), or pricing based on a fixed equivalent martingale measure (i.e. Black-Scholes pricing). 
In particular, the numerical results confirm that when employing an equal risk price both the writer and the buyer end up being exposed to risks that are more similar and on average smaller than what they would experience with the other approaches.
\end{abstract}


\section{Introduction}
One of the main challenges in pricing and hedging financial derivatives is that the market is often incomplete and thus there exists unhedgeable risk that needs to be further accounted for in pricing. In such a market, the price of a financial derivative cannot be set according to non-arbitrage theory as traditionally exploited in \cite{black1973pricing,merton1973theory,cox1979option,king2002duality}. Modern approaches to incomplete-market pricing can be broadly divided into two main categories.
The first one involves pricing a derivative based on a fixed \quoteIt{risk-neutral} martingale measure, either obtained from calibrating against market data \cite{hull1987pricing,heston1993closed,amin1993jump},  by minimizing the distance to a physical measure \cite{Delbaen95thevariance-optimal}, or by marginal indifference pricing \cite{brennan79}. The second category involves methods that rely on identifying the indifference price of a risk averse hedging problem, including for example good deals bound \cite{Jaschke2001}, expected utility indifference pricing \cite{CARR2001131}, or the quadratic hedging models \cite{FoellmerHedging,schweizer1996approximation,gourieroux1998mean,bertsimas2001hedging}. We refer readers to \cite{schweizer1999guided} and \cite{Staum}  for comprehensive surveys of these methods. 



In this paper, our focus is on studying a pricing method known as equal-risk pricing (ERP), which was first introduced in the recent work of \cite{guo2017equal}. The method can be considered falling into the second category mentioned above in that it involves the formulation of risk-averse hedging problems. In particular, it takes into account the risk preferences of both sides of a contract and seeks a fair unique transaction price that would ensure the minimal risk exposures (according to the formulated risk averse hedging problems) of both sides of a contract are equal. In \cite{guo2017equal}, special attention was paid to the case where risk is measured based on an expected disutility framework and where the market is incomplete due to no-short selling constraints on the hedging positions. They proved the existence and the uniqueness of the equal-risk price and provided pricing formulas for European and American options with payoffs that are monotonic in the underlying asset price. In the case where the constraints are lifted, they showed that the equal-risk price coincides with the price resulting from a complete market model.

\begin{figure}
\begin{center}
\begin{minipage}{150mm}
\subfigure[The option prices]{
\resizebox*{70mm}{!}{\includegraphics{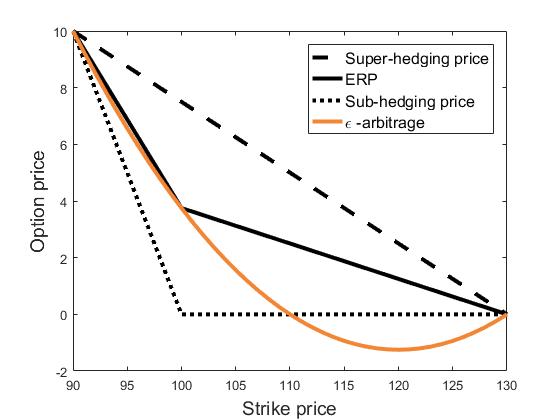}}\label{fig:LbUbA}}
\subfigure[The hedging loss results]{
\resizebox*{70mm}{!}{\includegraphics{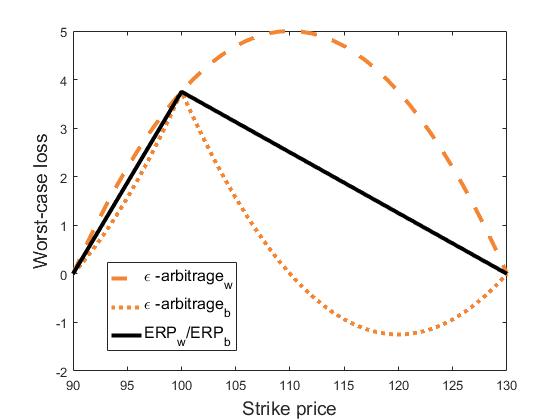}}\label{fig:LbUbB}}
\caption{Comparison of prices and hedging loss in the simple one period European call option pricing example. (a) shows the upper and lower bound of the no-arbitrage interval, together with the equal risk and $\epsilon$-arbitrage prices. (b) shows the worst-case loss incurred by each each party of the contract under their respective optimal hedging strategies. \label{fig:LbUb}}
\end{minipage}
\end{center}
\end{figure}

To put into perspective the strength of ERP, we should emphasize that most pricing methods focus only on a single side of the contract when formulating risk-averse hedging problems. The minimum price that a writer is willing to take according to a writer's hedging problem is however generally higher than the highest price that a buyer is willing to pay according to a buyer's hedging problem. Hence, there is a lack of mechanism to suggest a \quoteIt{transaction} price, i.e. acceptable to both the writer and the buyer. ERP provides such a mechanism by suggesting that a transaction should occur at a price which leaves both the writer and buyer with equivalent risk exposure.
To better illustrate this point, one can consider the example of pricing an European call option in the context where hedging can only occur at time zero. We further assume that the risk free rate is zero, and that the underlying stock price starts at a value of $100\$$ while its value at exercise time is known to be uniformly distributed over $[90,\,130]$. In this context, a risk averse writer might require that the price of an at-the-money option be set as high as $7.5\$$ to fully cover her risk while the buyer can use the same argument to require a price of 0\$. When a worst-case risk measure is used for both parties, one can show that the ERP allows the two parties to settle for the price of $3.75\$$ which exposes both of them to the same risk, i.e. $3.75\$$. 
%
%
Alternatively, one could suggest a transaction price based on a quadratic hedging scheme such as $\epsilon$-arbitrage pricing (see its application with worst-case risk measure in \cite{bandi2014robust}), yet  as shown in Figure \ref{fig:LbUb}, such paradigms can propose prices that leaves both parties with surprisingly uneven risk, giving in some case even rise to arbitrage opportunities (c.f. the negative price for strike prices between 110 and 130).
We refer the reader to Appendix \ref{app:analyticalSolution} for details of the analysis presented in this figure.

The contribution of the paper can be summarized as follows:
\begin{itemize}
\item We extend the definition of ERP to the set of all monotone risk \modified{measures} that can be interpreted as certainty equivalent measures (i.e. $\rho(t)=t$ for all $t\in\Re$). This class of risk measure includes the set of convex risk measures for which we establish for the first time that ERP \modified{is arbitrage-free under weak conditions and} actually reduces to computing the center of a so-called fair price interval (FPI). In comparison to the work of \cite{guo2017equal} which focused on an expected disutility framework that employs a fixed equivalent martingale measure, our generalized framework allows an arbitrary, and possibly different, probability measure to be used by each party, and corrects for the fact that the expected disutilities experienced by the two different parties are intrinsically non comparable. 
\item In the case of discrete time hedging, we show how the boundaries of such a fair price interval can be obtained by using dynamic programming for both European and American options as long as the convex risk measures employed by the two parties are one-step decomposable and satisfy a Markovian property (see Section \ref{sec:erpDef} for proper definitions). These dynamic programs are amenable to numerical computation given that they employ a finite dimensional state space. In the case of American options, they will also provide a different price depending on whether the buyer is willing to commit up front to an exercise strategy.
\item In the context where the underlying asset follows a geometric Brownian motion, we show for the first time how robust optimization can motivate the use of worst-case risk measures that only consider a subset of the outcome space. Similarly to risk neutral risk measures (i.e. that measure risk using expected value), these worst-case risk measures are easily interpretable and will satisfy the properties needed for dynamic programming to be used. On the other hand, unlike risk neutral measures, they also provide risk aware hedging policies. Our numerical experiments also indicate that the fair price interval might converge, as the number of rebalancing periods increases, to the Black-\modified{Scholes} price when an uncertainty set inspired by the work of \cite{bernhard2003robust} is properly calibrated in a market driven by a geometric Brownian motion. If supported theoretically, such a property would close the gap between risk neutral pricing and risk averse discrete-time hedging using worst-case risk measures. 
\item We present the first numerical study that provides evidence that equal risk prices allow both the writer and the buyer to be exposed to risks that are more similar and on average smaller than what they might
experience with risk neutral or quadratic hedging prices. In particular, \modified{when a worst-case risk measure is used, the risk inequity for the higher quantiles of each party's final loss will be reduced by a factor between 2 and 10 (depending on the type of option) compared to $\epsilon$-arbitrage  and Black-\modified{Scholes} prices}. This is done while keeping the average risk among the two parties to a similar or better level.
\end{itemize}

The paper is organized as follows. In Section \ref{sec:Motivation} we formally define the equal risk pricing framework and demonstrate that this price coincides with the mid point of the fair price interval when risk is captured using convex risk measures.
In Section \ref{sec:EuropeanAmerican}, we focus on the context of discrete-time option pricing and derive the dynamic programming equations that can be used to compute the equal risk price for European and American contingent claims. Next, in Section \ref{sec:numInvestigate}, an application of equal risk pricing is presented where the risk attitude of both writer and buyer is captured by so-called worst-case risk measures. A numerical study is also presented to validate the quality of prices obtained using the equal risk pricing paradigm both from the point of view of risk exposition for the parties and fairness. Finally, we conclude the paper in Section \ref{sec:conlusion}. We further refer the reader to an extensive set of Appendices describing detailed arguments supporting all propositions and lemmas presented in this article.

\section{The Equal Risk Pricing Framework}\label{sec:Motivation}
This section presents equal risk pricing framework and provides an interpretation of the price resulting from this model. In particular, we introduce the use of risk measures in pricing and hedging options based on this framework.

\subsection{The Equal Risk Pricing Model}\label{sec:erpDef}
To present the equal risk framework, we consider a model of the market proposed by \cite{xu2006risk}. Namely, we assume that the market is frictionless, i.e. there is no transaction cost, tax, etc. The filtered probability space is defined as $ (\Omega,\mathcal{F},\mathbb{F}=(\mathcal{F}_t)_{0\le t\le T},\mathbb{P}) $ and there is a money market account with zero interest rate, for simplicity and a risky asset $S_t$, with $0\le t \le T$, which is $ \mathcal{F}_t$-measurable. As in \cite{xu2006risk}, we assume that the risky asset $S_t$ is a locally bounded real-valued semi-martingale process.  Furthermore, the set of equivalent local martingale measures for $S_t$ is assumed non-empty to exclude arbitrage opportunity.

The set of all admissible self financing hedging strategies with the initial capital $ \price_0 $ is shown by $ \mathcal{X}(\price_0)$:
\begin{equation}\nonumber
\mathcal{X}(\price_0)=\left\{X_t \left| \exists \xi_s, \exists c\in\Re, \quad X_t=\price_0+\int_{0}^{t}\xi_s dS_s \geq c, \quad \forall t \in [0,T] \right. \right\}\,,
\end{equation}
in which, for each $t$, the decision  $ \xi_t $ is $ \mathcal{F}_t $-measurable and represents the number of shares of the risky asset in the portfolio, $ X_t $ is the accumulated wealth, and for simplicity, we assume that the risk free rate is zero. \modified{Although we impose very few restrictions on the hedging strategies in the set $\mathcal{X}(\price_0)$, as mentioned in \cite{xu2006risk}, the assumption of locally bounded real-valued semi-martingale $S_t$ allows many jump-diffusion and pure-jump models to be considered for the price process thus already giving rise to the possibility of an incomplete market. 
Alternatively, other definitions of $\mathcal{X}(\price_0)$ could be used here to model different characteristics of the market, e.g. discrete trading times (see Section \ref{sec:EuropeanAmerican}), or transaction costs, etc., without affecting the nature of our discussion.}

\modified{
We consider in this paper a class of payoff functions $F(\{S_t\}_{0 \leq t \leq T})$ that admit the formulation of 
$F(S_T,Y_T)$ where $Y_t$ is an auxiliary fixed-dimensional stochastic process that is $\mathcal{F}_t$-measurable. This class of payoff functions is common in the literature, for example in \cite{bertsimas2001hedging} and is more easily amenable to numerical methods (see Section 3 for more detailed discussions in a Markovian setting).
Here are a few examples, where we denote by $\{t_k\}_{k=1}^N$ a set of discrete time points.}

\begin{enumerate}
	\item \textbf{Options on the maximum value reached by a stock.}  \modified{The option pays off the maximum stock price reached during $\{t_k\}_{k=1}^N$ and is defined by}: 
	\begin{equation}\nonumber
		F\modified{\left(\max_{k=1,...,N}S_{t_k}\right)}\,,
	\end{equation}
	which is a function of the whole history of the stock price. In order to define the payoff as a function of some variables at the current time, $ Y_{t} $ can be defined as:
	\begin{equation}\nonumber
		Y_{t}=\max_{k:t_k \le t}S_{t_k}\,.
	\end{equation}
	Using this definition, the payoff of the option at the maturity can be written as a function of $ Y_T$.
	\item \textbf{Asian options.}  The payoff of an asian option is a function of the average stock price \modified{during $\{t_k\}_{k=1}^N$}:
	\begin{equation}\nonumber
		F\left(\frac{1}{N}\sum_{k=1}^{N}S_{t_k}\right)\,.
	\end{equation} 
	Letting $ Y_{t} $ have the following form:
	\begin{equation}\nonumber
		Y_{t}=\frac{1}{N}\sum_{k:t_k \le t}S_{t_k}\,.
	\end{equation}
	Again the payoff of the option at the maturity can be written as a function of $ Y_T $.
\end{enumerate}
For a broader range of options that can be as well recast into the above general form of payoff function, we refer readers to \cite{bertsimas2001hedging}.

Knowing the option payoff at maturity and assuming that the risk aversion of both participants are respectively characterized by two risk measures $\rho^w$ and $\rho^b$, i.e.  that map any random liability in $\mathcal{L}_p(\Omega,\mathcal{F}_{T},\mathbb{P})$\modified{, which one wishes to minimize,} to the set of real numbers (or infinity) and capture the amount of risk that is perceived by the participants, it is possible to define the minimal risk achievable by the option writer and the buyer as follows:
\begin{subequations}\label{equation:riskMeasures}
	\begin{align}
	&\varrho^w(\price_0)=\inf_{X\in\mathcal{X}(\price_0)} \rho^{w}(F(S_T,Y_T)-X_T) \label{equation:riskMeasures1}\\
	&\varrho^b(\price_0)=\inf_{X\in\mathcal{X}(-\price_0)} \rho^{b}(-F(S_T,Y_T)-X_T)\,,\label{equation:riskMeasures2}
	\end{align}
\end{subequations}
where $\price_0\in\Re$ represents the price that is charged by the writer to the buyer for committing to pay an amount $F(S_T,Y_T)$ to the buyer at time $T$. The quantities $ \varrho^w(\price_0), \varrho^b(\price_0) $ are the minimal risks associated to the optimal hedging of the writer and the buyer, respectively. In equation \eqref{equation:riskMeasures1}, the writer is receiving $ \price_0 $ as the initial payment and implements an optimal hedging strategy for the liability captured by $F(S_T,Y_T)$. On the other hand, in \eqref{equation:riskMeasures2} the buyer is assumed to borrow $\price_0 $ in order to pay for the option and then to manage a portfolio that will minimize the risks associated to his final wealth $F(S_T,Y_T)+X_T$. Note that while in practice the buyer usually might not buy an option by short-selling the risk free asset and might not optimize a portfolio with the intent of hedging the option, equation \eqref{equation:riskMeasures2} identifies the minimal risk that he could achieve by doing so, which can certainly serve as an argument in negotiating the price of the option given that this is always a possibility for him.

Following the work of \cite{guo2017equal}, the notion of minimal risk achievable for both participants can in turn be used to define an equal risk price as follows.
\begin{definition}\label{def:erp}
(Equal risk price) Given that both the writer's risk measure, $ \rho^w$, and buyer's risk measure $ \rho^b $ are interpretable as certainty equivalents, i.e.:
\begin{equation}\label{eq:constantRisk}
\forall c \in \mathbb{R}, \quad \rho(c)=c\,,
\end{equation}
and are monotone\footnote{Technically speaking, we also require both risk measures to satisfy Fatou's property and to satisfy $\rho(X)=\lim_{m\rightarrow \infty} \rho(\min(X,m))$ when $X$ is uniformly bounded from below while it should satisfy $\rho(X)=\lim_{m\rightarrow \infty} \rho(\max(X,-m))$ when $X$ is uniformly bounded above if the limit is finite otherwise be considered undefined (see \cite{xu2006risk} for details).}, i.e. having $ X,Y $ representing costs,
\begin{equation}\label{eq:Monotone}
\forall\,X,\,Y, \quad \;X\geq Y \mbox{a.s.} \Rightarrow\rho(X)\geq\rho(Y)\,,
\end{equation}
then the equal risk price is defined as the unique $\price_0^* $ that satisfies
\begin{equation}\label{eq:ERP}
 \varrho^w(\price_0^*)=\varrho^b(\price_0^*) \in \Re\,,
\end{equation}
when such a unique price exists.
\end{definition}
The reason for imposing equation $ (\ref{eq:constantRisk}) $ is to make sure that the units of $ \varrho^w $ and $ \varrho^b $ are comparable, i.e. that risk is expressed in the units of equivalent certain payoffs. Note that this assumption is not imposed in \cite{guo2017equal} where the notion of equal risk price can become arbitrary, e.g. when $\rho$ measures expected utility since utility functions are defined only up to positive affine transformations. 
Note also that this definition of equal risk price holds for general European options yet we will later present a similar definition for American options as well.

Besides the equal risk price, in an incomplete market with two risk averse market participants, another relevant and closely related concept takes the form of the following \quoteIt{fair price interval} (c.f. \cite {IntervalMarket}).

\begin{definition}
(Fair price interval) Given a writer's risk measure $\rho^w$ and buyer's risk measure $\rho^b$, the fair price interval is defined as the interval of prices for which both the writer and the buyer  are unable to exploit the market to completely hedge the risk of the contract they have agreed upon. Mathematically, the fair price interval takes the form $[\price_0^b,\price_0^w] $, where $\price_0^b=\sup\{\price_0 | \varrho^b(\price_0) \le 0 \} $ and $\price_0^w=\inf\{\price_0 | \varrho^w(\price_0) \le 0 \} $.
\end{definition}

\modified{It is worth differentiating the FPI from the no-arbitrage interval. In particular, the latter is defined as the interval $[\bar{\price}_0^b,\bar{\price}_0^w]$, such that:
\[\bar{\price}_0^b:=\sup\{\price_0|\exists X\in\mathcal{X}(-\price_0), F(S_T,Y_T)+X_T\geq 0 \mbox{ a.s.}\}\]
and
\[\bar{\price}_0^w:=\inf\{\price_0|\exists X\in\mathcal{X}(\price_0), F(S_T,Y_T)-X_T\leq 0 \mbox{ a.s.}\}.\]
One can easily exploit the fact that the risk measures are interpretable as certainty equivalents and monotone to show that the FPI, which accounts for the fact that the two parties are not arbitrarily risk averse, is necessarily a subset of the no-arbitrage interval. 
While the no-arbitrage price interval is always guaranteed to be non-empty, this is not necessarily the case for the FPI. An empty FPI captures the existence of a price for which both the writer and buyer end up being exposed to a negative risk thus making the ERP paradigm less relevant.
}

Note also that both the equal risk price and fair price interval can only be measured if the risk measures $\rho^w$ and $\rho^b$ are known. 
In practice, this might require both parties involved to provide supporting evidence for their respective choice of risk measure in the form of historical decisions that were taken using such measures. In the rest of the paper, we make the assumption that the true risk measures are known by each party.

\subsection{Equal Risk Pricing with Convex Risk Measures}

Since the work of \cite{artzner1999coherent}, it is now common to define coherent risk measures as risk measures that satisfy the following properties, where $ X $ and $ Z $  represent two random liabilities:
\begin{itemize}
	\item Monotonicity: if $ X \le Z \quad a.s.$ then $ \rho(X) \le \rho(Z) $
	\item Subadditivity: $ \rho(X+Z) \le \rho(X) + \rho(Z) $
	\item Positive homogeneity: If $ \lambda \ge 0 $, then $ \rho(\lambda X)=\lambda \rho(X) $
	\item Translation invariance: If $ m \in \mathbb{R} $, then $ \rho(X+m)=\rho(X)+m $
	\item Normalized risk: $\rho(0)=0$.
\end{itemize}
The first property naturally applies because if at any possible state that may happen the amount of  liability $ X $ is less than the liability $ Z $, then the risk of $ X $ is less than the risk of $ Z $. The second property is specifying that diversification does not increase the risk. The third property is implying that the risk of a position is linearly proportional to its size. Finally the last property implies that the addition of a sure amount to a random liability will decrease the risk by the same amount. By relaxing \modified{positive} homogeneity and \modified{subadditivity} with the following \modified{convexity} property, the family of risk measures becomes known as the larger family of ``convex risk measures" (\cite{follmer2011stochastic}):
\begin{itemize}
	\item Convexity: $ \rho(\lambda X+(1-\lambda)Z) \le \lambda \rho(X)+(1-\lambda)\rho(Z) $, for $ 0 \le \lambda \le 1 $.
\end{itemize}

Without loss of generality, In order to ensure that an equal risk price exists, we impose that participants are unable to design self financing hedging strategies that reach arbitrarily low risks.
\begin{assumption}\label{ass:AssRiskMeasure}
	We assume that the risk measures, $ \rho^w$ and $ \rho^b $ satisfy a \quoteIt{bounded market risk} assumption, i.e.
	\begin{equation}\nonumber
	0 \geq \inf_{X \in \mathcal{X}(0)} \rho^{w}(-X_T) > - \infty, \quad 	0 \geq \inf_{X \in \mathcal{X}(0)} \rho^{b}(-X_T) > - \infty\,.
	\end{equation}
In particular, if the risk measures are coherent, then this assumption implies that\footnote{Indeed, for a coherent risk measure, we have that $\inf_{X \in \mathcal{X}(0)} \rho^{w}(-X_T) <0$ implies that $\inf_{X \in \mathcal{X}(0)} \rho^{w}(-X_T) = -\infty$ because of positive homogeneity.
}
 \begin{equation}\nonumber
	\inf_{X \in \mathcal{X}(0)} \rho^{w}(-X_T) =0, \quad \inf_{X \in \mathcal{X}(0)} \rho^{b}(-X_T) =0\,.
	\end{equation}
\end{assumption}

Note that this assumption was also made in \cite{xu2006risk} (see Assumption 2.3) and reflects the fact that a participant believes that he cannot make an arbitrarily large risk-adjusted profit from trading in this market. We argue that in the context of equal risk price, it is made without loss of generality since if either risk measure violates the condition, then one should simply conclude that there exists no equal risk price as defined in Definition \ref{def:erp}. This is due to the fact that for all $\price_0$, we would have for example that
\begin{align*}
\varrho^w(\price_0)&=\inf_{X \in \mathcal{X}(\price_0)} \rho^{w}(\modified{F(S_T,Y_T)}-X_T)\leq  \inf_{X \in \mathcal{X}(\price_0)} (1/2)\rho^{w}(2\modified{F(S_T,Y_T)})+(1/2)\rho(-2X_T) \\
&=  (1/2)\rho^{w}(2\modified{F(S_T,Y_T)})-\price_0+\inf_{X \in \mathcal{X}(0)}(1/2)\rho(-2X_T)=-\infty \notin \Re\,. 
\end{align*}

An interesting conclusion can be drawn regarding the relation between the equal risk price and the fair price interval when both risk measures are convex risk measures. 

\begin{proposition}\label{pr:TrInv}
	Given that both $\rho^w$ and $\rho^b$ are convex risk measures, an equal risk price exists if and only if the fair price interval is bounded. Moreover, when it exists it is equal to:
\[\price_0^*=(\varrho^w(0)-\varrho^b(0))/2\,,\]
which is the center of the fair price interval if the latter is non-empty.
\end{proposition}

Based on the \modified{Proposition} \ref{pr:TrInv}, when using convex risk measures, the equal risk price can simply be found by evaluating the two boundaries of the fair price interval.


Following up on an important concern raised about the $\epsilon$-arbitrage pricing approach, based on a result from \cite{xu2006risk}, we can actually  confirm the fact that for convex risk measures  that satisfy the bounded market risk property, the equal risk price is arbitrage-free under weak conditions.
\begin{lemma}\label{thm:ERPisNoArbitrage}
If the fair price interval exists and is non-empty and both $\rho^w$ and $\rho^b$ are convex risk measures, then the equal risk price lies in the no-arbitrage price interval.
\end{lemma}



In what follows, we will show how the result of Proposition \ref{pr:TrInv} can be further exploited to identify the equal risk price of both European and American style options using dynamic programming in a context where hedging is implemented at discrete time points.

\begin{remark}\label{rem:nonTI}
We should note here that in the case where the risk measures do not satisfy the translation invariance property, one can still exploit the above observation that the equal risk price falls within the fair price interval \modified{and is therefore arbitrage-free assuming non-emptiness of this interval}. Namely, if such a price exists, 
one can identify it by employing a bisection algorithm that can establish $\Delta(\price_0):=\varrho^w(\price_0)-\varrho^b(\price_0)=0$. 
The convergence of a bisection method can rely on the fact that $\Delta(\price_0)$ is non-increasing and that it is greater or equal to zero at $\price_0^b$ and lower or equal to zero at $\price_0^w$. \modified{Finally, some guidance regarding the derivation of dynamic programming equations for this more general context can be found in Appendix \ref{app:DP:nonTransInv}.}
\end{remark}

\section{Discrete Dynamic Formulations for Equal Risk Pricing Framework}\label{sec:EuropeanAmerican}

In contexts where trading can only occur at specific periods of time $\{t_k\}_{k=0}^{K-1}\subset [0, T[$, one typically redefines the set of all admissible self financing hedging strategies in terms of the wealth accumulated at each period:
\begin{equation}\nonumber
\bar{\mathcal{X}}(\price_0)=\left\{X:\Omega\rightarrow \Re^K \left| \exists \{\xi_k\}_{k=0}^{K-1}, \quad X_k=\price_0+\sum_{k'=0}^{k-1} \xi_{k'} \Delta S_{k'+1} , \quad \forall k=1,\dots,K \right. \right\}\,,
\end{equation}
where $ \Delta S_{k+1}=S_{t_{k+1}}-S_{t_k} $ and where, for each $k=0,\dots,K-1$, the hedging strategy $\xi_k$ is a random variable adapted to the filtration $ \bar{\mathbb{F}}=(\mathcal{F}_{t_0},...,\mathcal{F}_{t_{K-1}}) $ and captures the number of \modified{shares of the} risky assets held in the portfolio during the period $[t_k,\;t_{k+1}]$. Finally, we assume that all random variables of interest in the discrete hedging problem are well-behaved.
\begin{assumption}\label{ass:niceRV}
There exists some $p\in[1,\,\infty]$ such that all $X\in\bar{\mathcal{X}}(\price_0)$ is such that for all $k=1,\dots,K$ we have that $X_k\in \mathcal{L}_p(\Omega,\mathcal{F}_{t_k},\mathbb{P})$ and that the payoff function $F(S_T,Y_T)\in\mathcal{L}_p(\Omega,\mathcal{F}_{T},\mathbb{P})$. 
 \end{assumption}

In particular, this assumption allows us to make use of a decomposability property of risk measures, which as shown in \cite{Ruszc06:CRM,Ruszczynski2010,pichler2018risk} 
is a key concept for producing a dynamic formulation for problems \eqref{equation:riskMeasures1} and \eqref{equation:riskMeasures2}. 

\begin{definition}\label{def:onestepDecomp}{(One-step decomposable risk measures)}
	The measure $ \rho:\mathcal{L}_p(\Omega,\mathcal{F}_{T},\mathbb{P})\rightarrow\Re$ is ``one-step decomposable" if there exists a set of risk measures  $\{\rho_{k}\}_{k=0}^{K-1}$ such that $\rho(X)=\rho_0(\rho_1(\cdots\, \rho_{K-2}(\rho_{K-1}(X))\cdots)$ and where each measure $\rho_{k}:\mathcal{L}_p(\Omega,\mathcal{F}_{k+1},\mathbb{P})\rightarrow \mathcal{L}_p(\Omega,\mathcal{F}_{k},\mathbb{P})$ is a conditional risk mapping (as defined in \cite{Ruszc06:CRM}, i.e. it satisfies the following properties:
	\begin{itemize}
	\item Conditional convexity : $\forall\,\theta\in[0,\,1],\,\forall\,X,Y\in\mathcal{L}_p(\Omega,\mathcal{F}_{k+1},\mathbb{P}),\,\rho_k(\theta Y+(1-\theta)X)\leq  \theta\rho_k(Y)+(1-\theta)\rho_k(X) \mbox{ a. s. }$
	\item Conditional monotonicity : $\forall\,X,Y\in\mathcal{L}_p(\Omega,\mathcal{F}_{k+1},\mathbb{P}),\,\,Y\geq X \mbox{ a. s. }\Rightarrow \rho_k(Y)\geq \rho_k(X) \mbox{ a. s. }$
\item	Conditional translation invariance : $\forall\,X\in\mathcal{L}_p(\Omega,\mathcal{F}_{k},\mathbb{P}),\,Y\in\mathcal{L}_p(\Omega,\mathcal{F}_{k+1},\mathbb{P}),\,\rho_k(X+Y)=X+\rho_k(Y) \mbox{ a. s. }$
	\end{itemize}
Additionally, a coherent risk measure is said to be \quoteIt{one-step coherently decomposable} if each measure $\rho_k$ also satisfies
	\begin{itemize}
	\item Conditional scale invariance : $\forall\,\alpha\geq 0,\,\forall\,X\in\mathcal{L}_p(\Omega,\mathcal{F}_{k+1},\mathbb{P}),\,\rho_k(\alpha X)= \alpha\rho_k(X) \mbox{ a. s. }$
	\end{itemize}\end{definition}

Among all risk measures that are one-step decomposable, a special class of risk measures can be shown to be especially attractive from a computational point of view. We will refer to these measures as Markovian risk measures which can be used when the filtered measurable space $(\Omega,\mathcal{F})$ is a progressively revealed product space.

\begin{definition}\label{ass:naturalFiltProbSpace}
The filtered probability space $(\Omega,\mathcal{F},\mathbb{F},\mathbb{P})$ is said to be supported on a progressively revealed product space if there exists a sequence $\{(\Omega_k,\Sigma_k)\}_{k=1}^K$ such that $(\Omega,\mathcal{F})$ is the product space, i.e. $\Omega:=\times_{k=1}^K \Omega_k$ and $\mathcal{F}:=\otimes_{k=1}^K \Sigma_k$, and $\mathbb{F}$ is the natural filtration in this space, i.e. $\mathbb{F}:=\{\sigma(\pi_{k'}:k'\leq k)\}_{k=1}^K$, where $\pi_k(\omega):=\omega_k$.
\end{definition}

With this definition in hand, we are now ready to define the class of Markovian risk measures. We note that a similar class of risk measure was proposed in \cite{Ruszczynski2010} in the context of a Markov decision processes. We however simplify the definition by exploiting the fact that conditional risk mappings are unaffected by decisions.

\begin{definition}\label{def:MarkovianRiskMeasure}(Markovian risk measure)
Given that $(\Omega,\mathcal{F},\mathbb{F},\mathbb{P})$ is supported on a progressively revealed product space defined through some $\{(\Omega_k,\Sigma_k)\}_{k=1}^K$, a one-step decomposable risk measure is said to be Markovian if there exists 
an $\mathbb{F}$ measurable stochastic process $\theta_k:\Omega \rightarrow \Re^m$, with $k=0,\dots,K$, that follows some dynamics $\theta_{k+1}(\omega)=f(\theta_k(\omega),\omega_k)$ for some $f:\Re^m\times \Omega_k\rightarrow\Re^m$, and 
some $\bar{\rho}_k:\mathcal{L}_p(\Omega_{k+1},\Sigma_{k+1},\mathbb{P}_{k+1})\times\Re^m\rightarrow \Re$, with $\mathbb{P}_{k+1}$ the marginalization of $\mathbb{P}$ on $\Omega_{k+1}$  such that:
\[\rho_k(X,\omega)=\bar{\rho}_k(\bar{\Pi}_k(X,\omega),\theta_k(\omega))\,,\]
where $\bar{\Pi}_k(X,\omega)$ is a random variable in $\mathcal{L}_p(\Omega_{k+1},\Sigma_{k+1},\mathbb{P}_{k+1})$ defined as $\bar{\Pi}_k(X,\omega,\bar{\omega}_{k+1})=X(\omega_{1:k},\bar{\omega}_{k+1},\omega_{k+1:K})$.
\end{definition}

Finally, given that the decomposable risk measure in this paper will be used in an arbitrage free financial market, one can formulate an assumption that imposes the bounded market assumption \ref{ass:AssRiskMeasure} on each conditional risk mapping.

\begin{assumption}\label{ass:BCMR}{(Bounded conditional market risk)}
	A one-step decomposable risk measure $\rho$ is said to express \quoteIt{bounded conditional market risk} if each conditional risk mapping $\rho_{k}$ satisfies the following properties:
\[0 \geq \inf_{\xi_k,\dots,\xi_{K-1}} \rho_{k,K}(-\sum_{\ell=k}^K\xi_\ell\Delta S_{\ell+1}) > - \infty \mbox{ a. s.}\,,\]
where $\rho_{k,K}(X):=\rho_k(\rho_{k+1}(\cdots\rho_{K-1}(X)\cdots)$. Furthermore, if it is conditionally scale invariant then $\inf_{\xi_k,\dots,\xi_{K-1}} \rho_{k,K}(-\sum_{\ell=k}^K\xi_\ell\Delta S_{\ell+1}) =0$.
\end{assumption}

In what follows, we derive dynamic equations that can be used to compute the equal risk price of European and American style options in discrete time trading. We further exploit the translation invariance and Markovian properties to reduce the dimension of the state space required to formulate the Bellman equations. \modified{We also conclude this section with an example of such equations when employing a recursive conditional value-at-risk measure.}

\ErickCommentsForFuture{It would be nice to include here an example (or two) of Markovian risk measure with bounded conditional market risk. The simplest would be the expected value when the reference measure satisfies the markovian property, or simply the robust measure with box uncertainty set.}

\subsection{European Style Options}\label{subsec:EUoption}

In order to evaluate the equal risk price of options of the form of  $F(S_T,Y_T)$ in discrete time with convex risk measures, as described in Proposition \ref{pr:TrInv}  one should solve problems \eqref{equation:riskMeasures1} and \eqref{equation:riskMeasures2} under the feasible set of strategies $\bar{\mathcal{X}}(0)$. Interestingly, the work of
\cite{pichler2018risk} provide simple arguments for deriving useful dynamic equations in contexts where the risk measures are one-step decomposable risk measures.

\begin{proposition}\label{prop:bellman}
Given that $ \rho^w $ and $ \rho^b $ are one-step decomposable risk measures as defined in Definition \ref{def:onestepDecomp}, then $\varrho^w(0)=V_0^w$ and $\varrho^b(0)=V_0^b$, where each $V_k^w$ and $V_k^b$ for $k=0,\dots,K$ are defined as follow: \\
\textbf{Writer's model:}
\begin{subequations}\label{eq:writerEuropDP}
\begin{align}
&V^w_k(\omega):=\inf_{\xi_k}\rho^w_k(-\xi_k \Delta S_{k+1}+ V^w_{k+1},\omega)\,,\,k=0,\dots,K-1\\
&V^w_K(\omega):=F(S_{\modified{K}}(\omega),Y_{\modified{K}}(\omega))\,.
\end{align}
\end{subequations}
\textbf{Buyer's model:}
\begin{subequations}\label{eq:buyerEuropDP}
\begin{align}
&V^b_k(\omega):=\inf_{\xi_k}\rho^b_k(- \xi_k \Delta S_{k+1} + V^b_{k+1},\omega)\,,\,k=0,\dots,K-1\\
&V^b_K(\omega):=-F(S_{\modified{K}}(\omega),Y_{\modified{K}}(\omega))\,\modified{,}
\end{align}
\end{subequations}
and assuming that  each $V_k^w\in\mathcal{L}_p(\Omega,\mathcal{F}_{k},\mathbb{P})$ and $V_k^b\in\mathcal{L}_p(\Omega,\mathcal{F}_{k},\mathbb{P})$.
Furthermore, the minimal risk hedging policy for both the writer and the buyer can be described respectively as:
\begin{align*}
\xi_k^{w*}(\omega)\in\arg\min_{\xi_k}\rho^w_k(-\xi_k \Delta S_{k+1}+ V^w_{k+1},\omega),\,\forall\,k=1,\dots,K-1\\
\xi_k^{b*}(\omega)\in\arg\min_{\xi_k}\rho^b_k(-\xi_k \Delta S_{k+1}+ V^b_{k+1},\omega),\,\forall\,k=1,\dots,K-1\,.
\end{align*}
\end{proposition}

We then get that if the filtered probability space is supported on a progressively revealed product space, if both $\Delta S_k$ and $\Delta Y_k:=Y_k-Y_{k-1}$ are measurable on $\Sigma_k$, such that they co-exist in both $\mathcal{L}_p(\Omega,\mathcal{F},\mathbb{P})$ and $\mathcal{L}_p(\Omega_{k},\Sigma_{k},\mathbb{P}_{k})$, and if  both $\rho^w$ and $\rho^b$ satisfy the Markovian assumption with respect to $\theta^w$ and $\theta^b$ respectively, then we can derive  finite dimensional Bellman equations that allow us to compute the equal risk price. These can be defined as follows:
\[\tilde{V}^w_K(S_K,Y_K,\theta_K^w):= F(S_K,Y_K)\,,\]
and recursively
\[\tilde{V}^w_k(S_k,Y_k,\theta_k^w):= \inf_{\xi_k} \bar{\rho}_k(-\xi_k\Delta S_{k+1} + \tilde{V}(S_k+\Delta S_{k+1},Y_k+\Delta Y_{k+1},f_k(\theta_k^w)),\theta_k)\,,\]
where each $\tilde{f}_k(\theta_k^w)$ can be considered a random variable in $\mathcal{L}_p(\Omega_{k+1},\Sigma_{k+1},\mathbb{P}_{k+1})$. These equations have the property that:
\[V^w_K(\omega)=\tilde{V}^w_K(S_K(\omega),Y_K(\omega),\theta_K^w(\omega))\,,\]
and recursively that if $V^w_{k+1}(\omega)=\tilde{V}^w_{k+1}(S_{k+1}(\omega),Y_{k+1}(\omega),\theta_{k+1}^w(\omega))$, then we have that:
\begin{align*}
&V^w_k(\omega)=\inf_{\xi_k}\rho^w_k(-\xi_k \Delta S_{k+1}+ V^w_{k+1},\omega)\\
&=\inf_{\xi_k}\rho^w_k(-\xi_k \Delta S_{k+1}+ \tilde{V}^w_{k+1}(S_{k+1},Y_{k+1},\theta_{k+1}^w),\omega)\\
&=\inf_{\xi_k}\bar{\rho}^w_k(\bar{\Pi}_k(-\xi_k \Delta S_{k+1}+ \tilde{V}^w_{k+1}(S_{k+1},Y_{k+1},\theta_{k+1}^w)),\theta_k^w(\omega)),\\
&=\inf_{\xi_k}\bar{\rho}^w_k(-\xi_k \Delta S_{k+1}+ \tilde{V}^w_{k+1}(S_k(\omega)+\Delta S_{k+1},Y_k(\omega)+\Delta Y_{k+1},f(\theta_k^w)),\theta_k^w(\omega))\\
&=\tilde{V}_k(S_k(\omega),Y_k(\omega),\theta_k(\omega))\,.
\end{align*}
From these derivations we see that $\varrho^w(0)=V^w_0=\tilde{V}^w_0(S_0,Y_0,\theta_0^w)$. In the case of the buyer, similar derivations lead to the Bellman equations:
\[\tilde{V}^b_K(S_K,Y_K,\theta_K^b):= -F(S_K,Y_K)\]
and 
\[\tilde{V}^b_k(S_k,Y_k,\theta_k^b):= \inf_{\xi_k} \bar{\rho}_k^b(-\xi_k\Delta S_{k+1} + \tilde{V}\modified{^b_{k+1}}(S_k+\Delta S_{k+1},Y_k+\Delta Y_{k+1},f_k(\theta_k^b)),\theta_k^b)\,,\]
which can be used to compute $\varrho^b(0)=V^b_0=\tilde{V}^b_0(S_0,Y_0,\theta_0^b)$. We can therefore conclude that $\price_0^*=(\tilde{V}^w_0(S_0,Y_0,\theta_0^w)-\tilde{V}^b_0(S_0,Y_0,\theta_0^b))/2$.

\subsection{American Style Option}\label{sec:amOption}

Contrary to European options, the exercise time of American options is flexible and up to the buyer's decision. Therefore, in the equal risk model we need to consider the interaction between an optimal exercise time and one-step decomposable risk measures. Similarly as in \cite{pichler2018risk}, we will define the exercise time as a \quoteIt{stopping time} adapted to the filtration $\mathbb{F}$, i.e. that it is a random variable  $\tau:\Omega\rightarrow \{0,\dots,K\}$, such that $ \{\omega: \tau(\omega)=t\} \in \mathcal{F}_t$, for all $\forall t=\{0,...,T\}$. Considering that the option payoff is now $ F_t(S_t,Y_t) \in L_p(\Omega,\mathcal{F}_t,\mathbb{P})$ at time $t$ if and only if $\tau=t$, we let $F_\tau(S_\tau,Y_\tau)$ capture the new payoff function which is defined as follows:
\begin{equation}\nonumber
F_\tau(S_\tau,Y_\tau):=\sum_{t=0}^{T}\mathbb{I}_{\{\tau=t \}}F_t(S_t,Y_t)\,,
\end{equation}
where $ \mathbb{I}_{\{\tau=t \}} $ is the indicator function, which is one for $ \tau=t $, and zero otherwise. We also redefine the set of self-financing hedging strategy to make its relation to $\tau$ explicit:
\modified{\[	\bar{\mathcal{X}}_\tau(\price_0):= \left\{X^\tau:\Tau\times\Omega\rightarrow \Re^K \middle| \begin{array}{l}\exists X_0^\tau=\price_0,  \forall k=1,\dots,K-1,\exists \xi_k, \{\hat{\xi}_k^i\}_{i=0}^{k}\quad\quad\quad\quad\\ X_{k+1}^\tau(\tau)=X_{k}^\tau(\tau)+(\xi_{k}\1\{\tau>k\}+\sum_{i=0}^{k}\hat{\xi}_{k}^i\1\{\tau=i\}) \Delta S_{k+1},\,\forall \tau\in\Tau\end{array} \right\}\,,\]}
\modified{where $\Tau$ is the set of all exercise time process, and} where each $\xi_k$ and $\hat{\xi}_k^i$ are $\mathbb{F}_k$-adapted.
Specifically, $\xi_k$ models the hedging strategy that is implemented at time $k$ when exercise has not occurred yet while $\hat{\xi}_k^i$ models the hedging strategy that is implemented at time $k$ when exercise occurred in period $k'=i$. 

\begin{remark}
We need to emphasize the fact that in most of the recent literature, hedging is considered to stop once the option is exercised. We intentionally omit making this assumption up front  and choose to model the possibility of hedging for both the writer and the buyer throughout the horizon. We will later show that when $\rho^w$ and $\rho^b$ are one-step coherently decomposable, the two approaches become equivalent, i.e. one can consider that $\hat{\xi}_k^i=0$ for all $k=0,\dots,T-1$ and all $i=0,\dots,k$ (see Section \ref{sec:NoHedging} for further discussion).
\modified{In cases where the assumption does not hold, we consider important to model hedging beyond exercise time in the buyer problem in order to avoid having incentives to delay exercise time simply to be able to benefit from later market opportunities. Similarly, in the writer's problem, if hedging stops at exercise time, the worst-case exercise time policy could be biased towards zero in order to prevent the writer from benefiting from later market conditions.}\modified{We also note that in any case, the analysis that follows can straightforwardly be adapted to a definition of the set of self-financing hedging strategies that explicitly enforces no hedging beyond exercise time.} 
\end{remark}

In this context, the definition of the equal risk framework needs to be adapted to account for the presence of $\tau$. In what follows, we will consider two formulations.

\begin{definition}\label{def:ERP_wCommit}
\textbf{ERP with Commitment.} Given that both the writer's risk measure, $ \rho^w$, and buyer's risk measure $ \rho^b $ are interpretable as certainty equivalents
and are strictly monotone with respect to certain amounts, then the equal risk price with commitment, when it exists, is defined as the unique $ \price_0^* $ for which there exists a stopping time policy $\tau^*$ that satisfies:
\begin{align*}
& \varrho^w(\price_0^*,\tau^*)=\varrho^b(\price_0^*,\tau^*) \in \Re&&\&&&\tau^* \in \arg\min_{\tau}\varrho^b(\price_0^*,\tau)\,,
\end{align*}
where
	\begin{align}
	&\varrho^w(\price_0,\tau)=\inf_{\modified{X^\tau}\in\bar{\mathcal{X}}_\tau(\price_0)} \rho^{w}(F(S_{\tau},Y_{\tau})-\modified{X_K^\tau(\tau)}) \label{equation:riskMeasuresCommit1}\\
	&\varrho^b(\price_0,\tau)=\inf_{\modified{X^\tau}\in\bar{\mathcal{X}}_\tau(-\price_0)} \rho^{b}(-F(S_{\tau},Y_{\tau})-\modified{X_K^\tau(\tau)})\,.\label{equation:riskMeasuresCommit2}
	\end{align}
\end{definition}

In simple terms, Definition \ref{def:ERP_wCommit} reflects the assumption that the buyer of the option commits to following a risk minimizing exercise strategy at the moment of purchasing the option. With this information, the writer can be more effective in hedging the option which, as will be shown, has the effect of giving rise to a lower equal risk price then when no commitment is made by the buyer. \modified{While we note that in practice, it might not be interesting for a buyer to commit upfront to an exercise strategy, the notion of ERP with commitment can serve the purpose of assessing the \quoteIt{cost of non-commitment}, which is a concept that is unique to the pricing of American options in an incomplete markets  (because of the multiplicity of arbitrage-free prices) and which can help interpreting the ERP without commitment.}

\begin{definition}\label{def:ERP_NoCommit}
\textbf{ERP without Commitment.} Given that both the writer's risk measure, $ \rho^w$, and buyer's risk measure $ \rho^b $ are interpretable as certainty equivalents
and are strictly monotone with respect to certain amounts, then the equal risk price without commitment, when it exists, is defined as the unique $ \price_0^* $ that satisfies:
\begin{align*}
 \varrho_\tau^w(\price_0^*)=\varrho_\tau^b(\price_0^*) \in \Re\,,
\end{align*}
where
	\begin{align}
	&\varrho_\tau^w(\price_0)=\inf_{\modified{X^\tau}\in\bar{\mathcal{X}}_\tau(\price_0)} \sup_\tau\rho^{w}(F(S_{\tau},Y_{\tau})-\modified{X_K^\tau(\tau)}) \label{equation:riskMeasuresNoCommit1}\\
	&\varrho_\tau^b(\price_0)=\inf_{\modified{X^\tau}\in\bar{\mathcal{X}}_\tau(-\price_0)} \inf_\tau\rho^{b}(-F(S_{\tau},Y_{\tau})-\modified{X_K^\tau(\tau)})\,.\label{equation:riskMeasuresNoCommit2}
	\end{align}
\end{definition}

Note that in this definition, the writer is unaware of the exercise strategy that will be employed by the buyer. He therefore considers the minimal risk of entering in this contract agreement as being the risk achieved by following an optimal hedging strategy that is adapted to both the filtration and the information about $\tau$ that is progressively revealed. 

We now demonstrate how the equal risk price necessarily increases when passing from the \quoteIt{with commitment} to \quoteIt{without commitment} framework.

\newcommand{\wncstar}{p_{\mymbox{nc}}^*}
\newcommand{\wcstar}{p_{\mymbox{c}}^*}
\newcommand{\tcstar}{\tau_{\mymbox{c}}^*}

\begin{lemma}\label{thm:noComitIsMoreExpansive} 
Given any American type option, the ERP with commitment $\wcstar$ is always smaller or equal to the ERP without commitment $\wncstar$.
\end{lemma}

In what follows, we derive dynamic programming equations that can be used to compute the ERP in both type of settings. Section \ref{sec:NoHedging} then exploits these equations to establish that when the risk measures are coherently decomposable, risk cannot be reduced by hedging beyond the exercise time.

\subsubsection{Bellman Equations for Equal Risk Price with Commitment}

We start with a simple lemma that extends the result of Proposition \ref{pr:TrInv} to the context of an American option with commitment.

\begin{lemma}\label{thm:midPriceWithCommitment} 
	Given that both $\rho^w$ and $\rho^b$ are convex risk measures, an equal risk price exists if and only if the fair price interval defined as:
\[[-\varrho^b(0,\tau_0),\,\varrho^w(0,\tau_0)]\,,\]	
where $\tau_0\in\arg\min_{\tau}\varrho^b(0,\tau)$, is bounded. Moreover, when it exists it is equal to the center of this interval which can be calculated as:
\[\price_0^*:=(\varrho^w(0,\tau_0)-\varrho^b(0,\tau_0))/2\,.\]
\end{lemma}

Lemma \ref{thm:midPriceWithCommitment} indicates that to evaluate the ERP, one needs to be able to compute $\varrho^w(0,\tau)$ and $\varrho^b(0,\tau)$ for any fixed exercise policy $\tau$, and to identify a procedure that can solve the optimal exercise time problem: $\min_{\tau}\varrho^b(0,\tau)$. As for the case of European options, all these elements can be characterized using dynamic programming equations.

\begin{proposition}\label{thm:AmWithCommitDP} 
Given that $ \rho^w $ and $ \rho^b $ are one-step decomposable risk measures as defined in \modified{Definition} \ref{def:onestepDecomp}, then $\varrho^w(0,\tau_0)=V_0^w(\tau_0)$ and $\varrho^b(0,\tau_0)=V_0^b(0)$, where for any exercise strategy $\tau$, each $V_k^w(\tau)$, $V_k^b(0)$, and $V_k^b(1)$ for $k=0,\dots,K$ are defined as follow: \\
\textbf{Writer's model:}
\begin{equation}
\begin{aligned}
&V^w_k(\tau,\omega):=\inf_{\xi_k}\rho^w_k( V^w_{k+1}(\tau)-\xi_k \Delta S_{k+1},\omega)+\1\{\tau(\omega)=k\}F(S_{k}(\omega),Y_k(\omega))\\
&V^w_K(\tau,\omega):=\1\{\tau(\omega)=K\}F(S_{K}(\omega),Y_K(\omega))\,.
\end{aligned}
\end{equation}
\textbf{Buyer's model:}
\begin{align}
&V^b_k(1,\omega):=\inf_{\xi_k}\rho^b_k(- \xi_k \Delta S_{k+1} + V^b_{k+1}(1),\omega)\label{eq:AmComitBuyerV1}\\
&V^b_k(0,\omega):=\min(V^b_k(1,\omega)-F(S_k(\omega),Y_k(\omega)),\;\inf_{\xi_k}\rho^b_k(- \xi_k \Delta S_{k+1} + V^b_{k+1}(0),\omega))\label{eq:AmComitBuyerV2}\\
&V^b_K(\bar{Z}_K,\omega):=-(1-\bar{Z}_K\modified{(\omega)})F(S_{K}(\omega),Y_K(\omega))\,.\label{eq:AmComitBuyerV3}
\end{align}
and assuming that each $V_k^w(\tau)$, $V_k^b(0)$, and $V_k^b(1)$ are in $\mathcal{L}_p(\Omega,\mathcal{F}_{k},\mathbb{P})$.
Furthermore, a feasible candidate for $\tau_0$ can be found using
\begin{align}
\tau_0(\omega)=\min\{k=0,\dots,K| V^b_k(0,\omega) = V^b_k(1,\omega)-F(S_k(\omega),Y_k(\omega)) \} \,.\label{eq:AmComitTau}
\end{align}
Finally, given that the option is sold at the equal risk price $\price_0^*=(V_0^w(\tau_0)-V_0^b(0))/2$ based on an exercise strategy $\tau_0$, the minimal risk hedging policy for both the writer and the buyer can be described respectively as:
\begin{align*}
\hat{\xi}_k^{i*}(\tau,\omega)\equiv\xi_k^{*}(\tau,\omega)\in\arg\min_{\xi_k}\rho^w_k( V^w_{k+1}(\tau)-\xi_k \Delta S_{k+1},\omega),\,\forall\,i=0,\dots,k,\,k=1,\dots,K-1
\end{align*}
for the writer, and
\begin{align}
\hat{\xi}_k^{i*}(\omega)\equiv\xi_k^{*}(\omega)\in\arg\min_{\xi_k}\rho^b_k(-\xi_k \Delta S_{k+1}+ V^b_{k+1}(\1\{\tau_0(\omega)\leq k\}),\omega),\,\forall\,i=0,\dots,k,\,k=1,\dots,K-1\label{eq:AmComitHedgeBuyer}
\end{align}
for the buyer.
\end{proposition}

In order for the evaluation of $\varrho^w(0,\tau_0)$ and $\varrho^w(0,\tau_0)$ to be computable numerically, it becomes essential to identifying Bellman equations on a finite dimensional state space. When the Markovian assumption holds for both $\rho^w$ and $\rho^b$ with respect to some process $\theta_k$, these can be derived as follows. For the buyer's problem, we have that:
\begin{equation}
\begin{aligned}
&\tilde{V}^b_k(1,S_k,Y_k,\theta_k):=\inf_{\xi_k}\bar{\rho}^b_k(- \xi_k \Delta S_{k+1} + \tilde{V}^b_{k+1}(1,S_k+\Delta S_{k+1},Y_k+\Delta Y_{k+1},f(\theta_k),\theta_k)\\
&\tilde{V}^b_k(0,S_k,Y_k,\theta_k):=\min(\tilde{V}^b_k(1,S_k,Y_k,\theta_k)-F(S_K,Y_K),\\
&\;\inf_{\xi_k}\bar{\rho}^b_k(- \xi_k \Delta S_{k+1} + \tilde{V}^b_{k+1}(0,S_k+\Delta S_{k+1},Y_k+\Delta Y_{k+1},f(\theta_k,\omega_{k+1}),\theta_k))\\
&\tilde{V}^b_K(\bar{Z}_K,S_k,Y_k,\theta_k):=-(1-\bar{Z}_K)F(S_K,Y_K)\,.
\end{aligned}\label{eq:buyerMarkovian}
\end{equation}
In order to obtain an optimal exercise policy, one can first observe that 
with
\begin{align*}
\tau_0(\omega)&=\min\{k=0,\dots,K| \tilde{V}^b_k(0,S_k(\omega),Y_k(\omega),\theta_k(\omega))\\ &=\tilde{V}^b_k(1,S_k(\omega),Y_k(\omega),\theta_k(\omega))-F(S_k(\omega),Y_k(\omega)) \} \,.
\end{align*}
Yet, when letting $Z_k:=\1\{\tau_0=k\}$ and $\bar{Z}_k:=\1\{\tau_0<k\}$, we can define: 
\[g_k(\bar{Z}_k,S_k,Y_k,\theta_k):=\1\{(\bar{Z}_k=0)\;\&\; (\tilde{V}^b_k(0,S_k,Y_k,\theta_k) =\tilde{V}^b_k(1,S_k,Y_k,\theta_k)-F(S_K,Y_K)) \}\]
so that 
\[Z_k(\omega)=g_k(\bar{Z}_k(\omega),S_k(\omega),Y_k(\omega),\theta_k(\omega))\,,\]
and 
\[\bar{Z}_{k+1}(\omega)=\bar{Z}_k(\omega)+g_k(\bar{Z}_k(\omega),S_k(\omega),Y_k(\omega),\theta_k(\omega))\,.\]
This implies that $\tau_0(\omega)=\sum_{k=0}^K g_k(\bar{Z}_k(\omega),S_k(\omega),Y_k(\omega),\theta_k(\omega))$ which can be implemented by exploiting the Bellman equations.
We can then proceed with describing the reduced equations for the writer's problem:
\begin{align*}
&\tilde{V}^w_k(\bar{Z}_k,S_k,Y_k,\theta_k):=\\
&\quad\inf_{\xi_k}\bar{\rho}^w_k(-\xi_k \Delta S_{k+1}+\tilde{V}^w_{k+1}(\bar{Z}_k+g_k(\bar{Z}_k,S_k,Y_k,\theta_k),S_k+\Delta S_{k+1},Y_k+\Delta Y_{k+1}),f_k(\theta_k)),\theta_k)+\\
&\quad\quad\quad\quad\quad\quad\quad g_k(\bar{Z}_k,S_k,Y_k,\theta_k)F(S_k,Y_k)\\
&\tilde{V}^w_K(\bar{Z}_K,S_K,Y_K,\theta_K):=g_K(\bar{Z}_K,S_K,Y_K,\theta_K)F(S_K,Y_K)\,,
\end{align*}
so that 
\[V_k^w(\tau_0,\omega)= \tilde{V}_k^w(\bar{Z}(\omega),S_k(\omega),Y_k(\omega),\theta_k(\omega))\,.\]

\subsubsection{Bellman Equations for Equal Risk Price without Commitment}

In the context of a contract where the buyer does not commit to a specific exercise policy, Proposition \ref{pr:TrInv} extends straightforwardly yet we provide the details in the following lemma for completeness.

\begin{lemma}\label{thm:midPriceWithoutCommitment} 
	Given that both $\rho^w$ and $\rho^b$ are convex risk measures, an equal risk price exists if and only if the fair price interval defined as $[-\varrho_\tau^b(0),\,\varrho_\tau^w(0)]$, is bounded. Moreover, when it exists it is equal to the center of this interval which can be calculated as $\price_0^*:=(\varrho_\tau^w(0)-\varrho_\tau^b(0))/2$.
\end{lemma}

The main difference between this case and the case with commitment is that in order to compute $\varrho_\tau^w(0)$, now there is a need to further determine the worst-possible exercise policy that 
the writer would hedge against. Since the whole hedging problem for the writer now takes the form of a minimax optimization problem, additional care has to be taken to ensure the decisions of hedging and exercising (the options) are executed in the right order when formulating the Bellmanequations. In particular, we proceed by fixing first the hedging decisions and identifying recursive equations that solve the worst-case exercise time problem (\cite{pichler2018risk}). We then use the arguments based on the interchangeability principle in dynamic programming  (see \cite{pichler2018risk}) to establish that the hedging decisions that minimize the recursive equations globally can be obtained from decisions that minimize the recursive equations stage-wise. The details can be found in the appendix and this leads to the following dynamic programming equations. On the other hand, it is not hard to confirm that the computation of $\varrho_\tau^b(0)$ for the buyer coincides with the computation required in the case of commitment. 

\begin{proposition}\label{thm:DPwithoutCommitment} 
Given that $ \rho^w $ and $ \rho^b $ are one-step decomposable risk measures as defined in Definition \ref{def:onestepDecomp}, then $\varrho_\tau^w(0)=V_0^w(0)$ and $\varrho^b(0)=V_0^b(0)$, each $V_k^w(0)$ and $V_k^w(1)$ for $k=0,\dots,K$ are defined as follow: \\
\textbf{Writer's model:}
\begin{align}
&V^w_k(1,\omega):=\inf_{\xi_k}\rho^w_k(- \xi_k \Delta S_{k+1} + V^w_{k+1}(1),\omega)\label{eq:erpnoComV1}\\
&V^w_k(0,\omega):=\max(V^w_k(1,\omega)+F(S_k(\omega),Y_k(\omega)),\;\inf_{\xi_k}\rho^w_k(- \xi_k \Delta S_{k+1} + V^w_{k+1}(0),\omega))\label{eq:erpnoComV2}\\
&V^w_K(\bar{Z}_K,\omega):=(1- \bar{Z}_K\modified{(\omega)})F(S_K(\omega),Y_K(\omega))\,,\label{eq:erpnoComV3}
\end{align}
while $V_k^b(0)$ and $V_k^b(1)$ are defined as in equations \eqref{eq:AmComitBuyerV1}-\eqref{eq:AmComitBuyerV3} and assuming that each $V_k^w(0)$, $V_k^w(1)$, $V_k^b(0)$, and $V_k^b(1)$ are in $\mathcal{L}_p(\Omega,\mathcal{F}_{k},\mathbb{P})$.
Furthermore, given that the option is sold at the equal risk price $\price_0^*=(V_0^w(0)-V_0^b(0))/2$, a minimal risk hedging policy for the writer can be described as:
\begin{align*}
\hat{\xi}_k^{i*}(\bar\omega)\equiv\xi_k^{*}(\omega)\in\arg\min_{\xi_k}\rho^w_k(-\xi_k \Delta S_{k+1}+ V^w_{k+1}(\1\{\tau\leq k\}),\omega),\,\forall\,i=0,\dots,k,\,k=1,\dots,K-1\,,
\end{align*}
where $\tau$ is the observed exercise strategy.
In the case of the buyer, a risk minimizing hedging strategy is as in equation \eqref{eq:AmComitTau} while a risk minimizing exercise strategy can be found using equation \eqref{eq:AmComitHedgeBuyer}.
\end{proposition}

When the Markovian assumption holds with respect to some process $ \theta_k $, we can again derive finite dimensional equations. In particular, for the buyer's problem, these are exactly as presented in equations (\ref{eq:buyerMarkovian}). On the other hand, for the writer's problem, we have that:
\begin{equation}
\begin{aligned}
&\tilde{V}^w_k(1,S_k,Y_k,\theta_k):=\inf_{\xi_k}\bar{\rho}^w_k(- \xi_k \Delta S_{k+1} + \tilde{V}^w_{k+1}(1,S_k+\Delta S_{k+1},Y_k+\Delta Y_{k+1},f(\theta_k),\theta_k)\\
&\tilde{V}^w_k(0,S_k,Y_k,\theta_k):=\max(\tilde{V}^w_k(1,S_k,Y_k,\theta_k)+F(S_K,Y_K),  \\
&\inf_{\xi_k}\bar{\rho}^w_k(- \xi_k \Delta S_{k+1} + \tilde{V}^w_{k+1}(0,S_k+\Delta S_{k+1},Y_k+\Delta Y_{k+1},f(\theta_k),\theta_k))\\
&\tilde{V}^w_K(\bar{Z}_K,S_k,Y_k,\theta_k):=(1-\bar{Z}_K)F(S_K,Y_K)\,.
\end{aligned}\label{eq:writerMarkovian}
\end{equation}

\subsubsection{On the Value of Hedging Beyond the Exercise Time}\label{sec:NoHedging}

As pointed out in the beginning of Section 3, our dynamic programming (DP) formulations of the hedging problem are more general in that they allow for the possibility of hedging after the exercise of the options. This in principle provides the opportunities for both the writer and buyer to seek further risk reduction. But at the same time it adds additional complexity to the DP formulation, which becomes computationally more costly to solve than the DP that assumes no hedging after exercise of the options. In this section, we identify the condition under which hedging beyond the exercise time actually does not reduce risk. In particular, based on our general DP formulation, we find that it is actually optimal to stop hedging after the exercise time if the employed risk measure is coherent.   

\begin{corollary}\label{thm:AmStopHedge}
If $\rho^w$ is one-step coherently decomposable, then it becomes optimal for the writer to terminate the hedging strategy at the exact moment that the American option is exercised. The same applies to the buyer.
\end{corollary}

As detailed in Appendix \ref{app:proof:AmStopHedge}, this observation is closely related to the assumption of bounded market risk, in which case there exists no risk reduction opportunity when measured according to a coherent risk measure. Since in this case hedging beyond exercise time adds no value, one can simply employ a DP formulation that assumes that hedging stops at the exercise time. 

\ErickCommentsForFuture{It might be interesting to provide a clear and simple example where the non-coherently decomposable risk measure produces a hedging beyond exercise. Something like a single-stage example where the exercise is one step 1, yet the convex risk measure proposes an investment in the risky asset. This could wait until the revision...}

In the next section, we elaborate on a specific class of coherently decomposable risk measure, referred to as \quoteIt{worst-case risk measures}. We further provide numerical evidence on the quality of prices obtained using such risk measure both from the point of view of risk exposure and fairness.

\subsection{\modified{Recursive Conditional Value-at-Risk Example}}

In this section, we provide a specific example of the Markovian counterpart of a popular one-step decomposable risk measure. We demonstrate how
our results can be applied to this risk measure so
as to write down the corresponding Bellman equations. 

We start by assuming the stochastic processes $S_{k}$ and $Y_{k}$
admit the following recursive representation, which is common in many
applications:
\[
S_{k+1}=f(S_{k},\epsilon_{k+1}),\;Y_{k+1}=g(Y_{k},\epsilon_{k+1}),
\]
for some $f:\Re\times \Re^{n_\epsilon}\rightarrow\Re$ and $g:\Re^{n_y}\times \Re^{n_\epsilon}\rightarrow\Re^{n_y}$, and $(\epsilon_{1},...,\epsilon_{K})$ is a realisation of the progressively revealed outcome space $\Omega:=\times_{k=1}^{K}\mathbb{R}^{n_\epsilon}$ equipped with probability measure $\mathbb{P}$ and natural filtration $\mathbb{F}$.

\begin{definition}
(Recursive conditional value-at-risk) Given a random variable
$X$ and a process  $\{ \beta_k\}_{k=0}^{K-1}$ which is $\mathcal{F}_k$ measurable, i.e. $\beta_k:=\Re^{n_\epsilon^k}\rightarrow \Re$, the recursive conditional value-at-risk measure $\rho$
is a one-step decomposable risk measure obtained using a conditional value-at-risk measure, defined as
\[
\rho_k(X,\omega)=\inf_{t}\;t+\frac{1}{1-\beta_k(\epsilon_1(\omega),\dots,\epsilon_k(\omega))}\mathbb{E}\left[(X-t)^{+}\mid|\epsilon_1(\omega),\dots,\epsilon_k(\omega)\right],
\]
as the conditional risk mapping.
\end{definition}

Note that the recursive conditional value-at-risk measure defined above only qualifies, in its general form, as a Markovian risk measure if one considers $\theta_k:=[\epsilon_1^T\;\;\dots\;\;\epsilon_k^T]^T$. This can quickly give rise to the curse of dimensionality when constructing and solving the associated DP formulation. To circumvent this issue, a common practice is to assume that the $\epsilon_k$ process satisfies the Markov property, i.e. $\mathbb{E}[X|\epsilon_1(\omega),\dots,\epsilon_k(\omega)]=\mathbb{E}[X|\epsilon_k(\omega)]$ for all $X\in\mathcal{L}_p(\Omega,\mathcal{F}_{k+1},\mathbb{P})$. One however also needs an additional assumption about the $\beta_k$ process such that $\beta_k=h(\beta_{k-1},\epsilon_k)$ for some $h:\Re\times\Re^{n_\epsilon}\rightarrow \Re$, in order to satisfy the Markovian risk measure assumption under a process $\theta_k:=[\beta_{k-1} \;\;\epsilon_{k-1}^T]^T$.

We can now summarize the Bellman equations that can be derived for
the case of a recursive conditional value-at-risk that is Markovian with respect to $\theta_k$ by following the result and discussions in Section 3.1. Namely, the writer problem's Bellman equations for the case
of European options can be written as follows:
\begin{align*}
\tilde{V}_{0}^{w}(S_{0}&,Y_{0})=\inf_{\xi,t}\;t+\frac{1}{1-\beta_{0}}\mathbb{E}[(-\xi(f(S_{0},\epsilon_1)-S_0)+\tilde{V}_{1}^{w}(f(S_{0},\epsilon_1),g(Y_{0},\epsilon_1),h(\beta_{0},\epsilon_1),\epsilon_1)-t)^{+}],
\\
\tilde{V}_{k}^{w}(S_{k}&,Y_{k},\beta_k,\epsilon_k)=\\
&\inf_{\xi,t}\;t+\frac{1}{1-\beta_{k}}\mathbb{E}[-\xi( f(S_{k},\epsilon_{k+1})-S_k)+\tilde{V}_{k+1}^{w}(f(S_{k},\epsilon_{k+1}),g(Y_{k},\epsilon_{k+1}),h(\beta_{k},\epsilon_{k+1}),\epsilon_{k+1})-t)^{+}|\epsilon_{k}]
\end{align*}
and $\tilde{V}_{K}^{w}(S_{K},Y_{K},\beta_{K},\epsilon_K)=F(S_{K},Y_{K})$. Similar
Bellman equations can be derived for the buyer and we omit them for brevity. 

In the case of American option, we can follow the result and discussions
in Section 3.2.2. to write down the following Bellman equations for the
buyer:
\begin{align*}
\tilde{V}_{k}^{b}(0,&S_{k},Y_{k},\beta_{k},\epsilon_k)  =\min\bigg(-F(S_{k},Y_{k}),\;\\
 & \inf_{\xi,t}\;t+\frac{1}{1-\beta_{k}}\mathbb{E}[(-\xi_{k}(f(S_{k},\epsilon_{k+1})-S_k)+\tilde{V}_{k+1}^{b}(0,f(S_{k},\epsilon_{k+1}),g(Y_{k},\epsilon_{k+1}),h(\beta_{k},\epsilon_{k+1}),\epsilon_{k+1})-t)^{+}|\epsilon_k]\bigg)
\end{align*}
and
\[
\tilde{V}_{K}^{b}(\bar{Z}_{K},S_{K},Y_{K},\beta_{K},\epsilon_K)=-(1-\bar{Z}_{K})F(S_{K},Y_{K}),
\]
where we exploited the fact that $\tilde{V}_{k}^{b}(1,S_{k},Y_{k},\beta_{k}\epsilon_k)=0$ since the conditional value-at-risk conditional risk mapping is coherent and Corollary \ref{thm:AmStopHedge} applies.
We omit the writer's equations for brevity. 

The arguments used above can be employed for many other recursive risk measures as long as $\rho_k$ is the conditional analog of a law-invariant coherent risk measure. 
\begin{example}
One obtains a recursive mean semi-deviation measure when using a conditional risk mapping $\rho_k$
defined as $\rho_{k}(X)=\mathbb{E}[X|{\cal F}_{k}]+\kappa_k\mathbb{E}[((X-\mathbb{E}[X|{\cal F}_{k}])_{+})^{r}|{\cal F}_{k}]^{\frac{1}{r}}$,
where $\kappa_k$ is ${\cal F}_{k}$-measurable.
\end{example}
\begin{example}
One obtains a recursive mean CVaR measure when using a conditional risk mapping $\rho_k$
defined as $\rho_{k}(X)= \mathbb{E}[X|{\cal F}_{k}]+\kappa_k\left(\inf_t\,t+\frac{1}{1-\beta_k}\mathbb{E}[(X-t)^+|{\cal F}_{k}]\right)$,
where $\kappa_k>0$ and $\beta_k\in [0,\,)1$ are ${\cal F}_{k}$-measurable.
\end{example}

On the other hand, it is worth emphasizing that one-step decomposable risk measure that are constructed based on the composition of law invariant coherent risk measures as suggested above are not law invariant unless the conditional mappings are expectation or worst-case risk measures \cite{SHAPIRO2012436}. This motivates us, in the following section, to focus our numerical study on the latter class of risk measures.

\section{Numerical Study with Worst-case Risk Measures}\label{sec:numInvestigate}


In this section, we provide necessary details of implementing the equal risk pricing model in the case where the risk measure takes the form of a worst-case risk measure. In particular, such form of risk measures has been considered in the literature of robust optimization, which requires the specification of an uncertainty set $U$ over which the worst-case loss is calculated. The $\epsilon$-arbitrage pricing model mentioned earlier in the introduction is one example that employs an uncertainty set motivated by central limit theorem. While the $\epsilon$-arbitrage pricing model does not distinguish between the writer's and the buyer's loss, the equal risk pricing model proposed in this paper does, and one of our goals in this section is to demonstrate numerically the strength of the equal risk pricing model over the $\epsilon$-arbitrage pricing model. We will also benchmark the equal risk pricing model against the Black-Scholes pricing model in the case of European option, and against the binomial pricing model in the case of American option.

To facilitate the comparisons between the aforementioned models, we start by considering a market of assets that are driven by a Geometric Brownian Motion (GBM). We assume that the asset returns can only be observed at a set of uniformly distributed time points on the interval $[0,T]$ such that each time point $t_k:=kT/K$, $k=1,...,K$. Without loss of generality, we can write $S_{t_k} = S_0 \Pi_{l=1}^k(1+r_k)$ to denote the dynamic of asset price given a random vector of observed returns taking values in $\Re^K$ and an initial asset price $S_0$.
In order to formalize worst-case risk measures over such a market, we consider an outcome space $\Omega:= \Re^K$ and an associated filtered probability space  $(\mathbb{R}^K,{\cal B}(\mathbb{R}^K),\bar{\mathbb{F}},\bar{\mathbb{P}})$, where ${\cal B}(\mathbb{R}^K)$ is the Borel $\sigma$-algebra on $\mathbb{R}^K$, and  $\bar{\mathbb{F}}:=\{ \sigma(r_{k'}:k'\leq k)\}$ is the natural filtration. We let $\bar{\mathbb{P}}$ be the probability measure that captures
\[(1+r_{k})\sim({\rm i.i.d.\;}{\rm Lognormal}(\mu T/K,\sigma^{2}T/K),\;k=1,...,K\,,\]
where $\mu$ and $\sigma$ are the statistics of the GBM per unit of time $T$. Note that this filtered probability space is supported on a progressively revealed product space as defined in Definition \ref{ass:naturalFiltProbSpace}. For the sake of convenience, we reformulate the hedging decision problem in terms of how much money is invested in the risky asset at each time point, denoted by $\zeta_0,\dots,\zeta_{K-1}$, instead of the number of \modified{shares of the} risky assets, i.e. $\zeta_k=\xi_k S_{t_k}$. This leads to the following equation representing the evolution of wealth:
\begin{equation}\nonumber
X_{k}=\price_0+\sum_{k'=0}^{k-1} \zeta_{k'} r_{k'+1}, \forall k=1,...,K. 
\end{equation}

In this numerical study, we will assume that the writer and buyer are employing a risk measure that is motivated by robust optimization. In particular, we will assume that they are concerned about the worst-case performance for realizations that arise in a predefined uncertainty set $\U$. We therefore define a worst-case risk measure as:
\[
\rho(X)=\mbox{ess}\sup_{\mathbb{U}(\U)}X\,,
\]
where $\U \subset ]-1,\,\infty[^K$ is compact and regular closed, and where $\mathbb{U}(\U)$ refers to the uniform distribution over $\U$. In what follows, we will simplify presentation by employing the notation from robust optimization with $\U$ as the so-called uncertainty set:
\[
\rho(X)=\sup_{r\in \U}X(r)\,.
\]
Clearly, this risk measure is necessarily monotone, translation invariant, and
coherent. Moreover, it is also one-step decomposable using:
\[\rho_k(X,r):=\left\{\begin{array}{cl}\sup_{r'\in \U:r'_{1:k}=r_{1:k}} X(r')&\mbox{if $\exists\,r'\in\U,\,r'_{1:k}=r_{1:k}$}\\ X([r_{1:k}^T\;\; 0_{k+1:K}^T]^T) & \mbox{otherwise}\end{array}\right.\,,\]
where $r_{1:k}\in\Re^k$ refers to the first $k$-th first terms of $r$, and where 
$X([r_{1:k}^T\;\; 0_{k+1:K}^T]^T)$ is short for
\[\inf_{\epsilon>0}\esssup_{r'\in ]-1,\,\infty[^K:r'_{1:k}=r_{1:k},\|r'_{k+1:K}\|_\infty\leq \epsilon} X\,.\]
Note that the conditional measure that is used for the case where $\nexists\,r'\in\U,\,r'_{1:k}=r_{1:k}$ can be arbitrary if one is only interested in calculating $\varrho(0)$ given that the latter is unaffected by the level of loss when $r\notin \U$. In practice however, one might get a \quoteIt{better} hedging policy by employing a more risk-aware measure than $X([r_{1:k}^T\;\; 0_{k+1:K}^T]^T)$. 
Indeed, one can confirm that $\rho$ can equivalently be described as:
\begin{align*}
\rho(X)&= \sup_{r^1\in \U} \sup_{r^2\in\U:r_1^1=r_1^2}  \sup_{r^K\in\U:r_{1:K-1}^{K-1}=r_{1:K-1}^{K}} X(r^K)\\
 &= \rho_0(\rho_1(\cdots\rho_{K-1}(X)\cdots)\,.
\end{align*}
In many cases, the one-step decomposable risk measure $\rho_k$ can be further shown to satisfy the Markov property, e.g. with the uncertainty sets presented in the following sections. One can then follow the   discussion in the Section \ref{sec:EuropeanAmerican} to write down the dynamic programming equations for both cases of European and American options. 

In all of our experiments, we consider an option with maturity $T=1$ (year) that is written over an asset with $\mu=0.0718$ (annualized mean), $\sigma = 0.1283$ (annualized volatility), and with an initial price $S_0 = 1000$. Our choices of values for $\mu$ and $\sigma$ come from \cite{franccois2014optimal} where they were calibrated on historical data of the S\&P 500 index.

\subsection{Comparison with $\epsilon$-arbitrage Pricing}\label{sec:abprice}

We present in this section the results of comparing the equal risk pricing model with the $\epsilon$-arbitrage pricing model proposed in \cite{bandi2014robust}. Recall that the uncertainty set $\U$ employed in \cite{bandi2014robust} admits the following form motivated by the central limit theorem:


\begin{equation}\label{eq:cltUS}
\U_1=\left \{r \in \mathbb{R}^K \middle| \left| \frac{\sum_{\ell=1}^{k}\log (1+r_\ell)-\mu kT/K}{\sigma\sqrt{kT/K}} \right| \le \Gamma, \forall k \in \{1,...,K\} \right \}\,,
\end{equation}
where $K$ is the number of periods up to the maturity of the option, and $\Gamma$ denotes the "budget" of uncertainty at each time point $t_k$. Unfortunately, the above uncertainty set cannot be directly applied in the equal risk pricing model, since its associated worst-case risk measure does not necessarily satisfy the bounded conditional  market risk, i.e. Assumption \ref{ass:BCMR}. We show in the following how the set can be slightly modified so that it satisfies Assumption \ref{ass:BCMR}.

\begin{lemma} \label{lemu1}
	Given that $\mu \sqrt{kT/K}/\sigma  \leq \Gamma$ for all $k\in\{1,\dots,K\}$, then the worst-case risk measure $\rho$ that exploits $ \U'_1=\U_1 \cap \mathcal{W}$ with
	\[	\mathcal{W}=\left \{r \in \mathbb{R}^K \middle| \max_{k'\in\{k+1,\dots,K\}} \left| \frac{\sum_{\ell=1}^{k}\log (1+r_\ell)-\mu k'T/K}{\sigma\sqrt{k'T/K}} \right| -\Gamma\le 0, \forall k \in \{1,...,K\}\right \}
\]
	 satisfies both assumptions \ref{ass:AssRiskMeasure} and \ref{ass:BCMR}.
\end{lemma}
It is worth noting that the above set is smaller than the original set $\U_1$, as it excludes the sample paths that can lead to infinitely small risk. But as shown in the appendix, the above modified set is in some sense 
the ``largest" subset of $\U_1$ that make the worst-case risk measure satisfy assumptions \ref{ass:AssRiskMeasure} and \ref{ass:BCMR}. It is not hard to confirm that when using $\U_1'$, the worst-case risk measure is  Markovian with respect to  $\theta_k:=\sum_{\ell=1}^k \log(1+r_{\ell})$ (\modified{see appendices \ref{proveMarkov} and \ref{app:implementationDetails} respectively for a proof and details about the implementation of the dynamic program}).

The parameter that needs to be further determined in our experiments is the budget parameter $\Gamma$. To do so, we start by first sampling $100000$ price paths from the GBM and then calibrating $\Gamma$ so that the uncertainty set
would cover at least $95\%$ of the paths.
In Table \ref{table:ERPriceAndEpsilonPrice}, we present the option prices generated from the equal risk and the $\epsilon$-arbitrage pricing models for various values of $K$ and different types of options: In-The-Money (ITM), At-The-Mone (ATM), and Out-of-The-Money (OTM). The table also presents the fair price intervals.

From Table \ref{table:ERPriceAndEpsilonPrice}, we can make a few observations about the prices generated from the two models. Firstly, in the case of OTM, the prices generated from the $\epsilon$-arbitrage pricing model are consistently lower than the prices generated from the equal risk pricing model. This is consistent with what was observed for the single period example in the introduction. Recall that in the case of single period (see Figure \ref{fig:LbUb}), the $\epsilon$-arbitrage prices were always smaller or equal to ERP and differed most significantly from ERP when the options was out-of-the-money. Indeed, we see from Table \ref{table:ERPriceAndEpsilonPrice} that in the case of ITM and ATM, the prices of the two models are more similar (without any clear dominance), but in the case of OTM options, the ERP is always significantly bigger than the $\epsilon$-arbitrage price. This confirms that the $\epsilon$-arbitrage pricing model can generate unrealisticly low prices even in a multi-period hedging problem. Secondly, one can notice in Table \ref{table:ERPriceAndEpsilonPrice} that the FPI lower bounds always take the value of zero, i.e. the buyer's perception of minimal hedging risk is invariant to the number of rebalancing periods. While this may seem counter-intuitive, we can actually find an explanation by taking a closer look at the structure of the uncertainty set $\U_1'$. Namely, the set only imposes upper bounds on the variations of the underlying asset process. It turns out however that for the buyer's optimal hedging strategy, the worst paths are paths where the prices stay constant. These paths remain feasible regardless of the value of $\Gamma$. This explains why the lower bounds always reach the lowest possible value, i.e. zero, regardless of the number of hedging periods. Lastly, in Table \ref{table:ERPriceAndEpsilonPrice},  we provide also the prices generated from the Black-Scholes formula, and one can notice that the prices from equal risk pricing are always higher than the Black-Scholes prices. This can also be explained by the conservativeness of the FPI lower bounds, which drives up the ERP. We will discuss in the next section how such an issue might be resolved with a different choice of uncertainty set.

\begin{figure}
\begin{center}
\begin{minipage}{160mm}
\subfigure[Average (ITM)]{
\resizebox*{50mm}{!}{\includegraphics{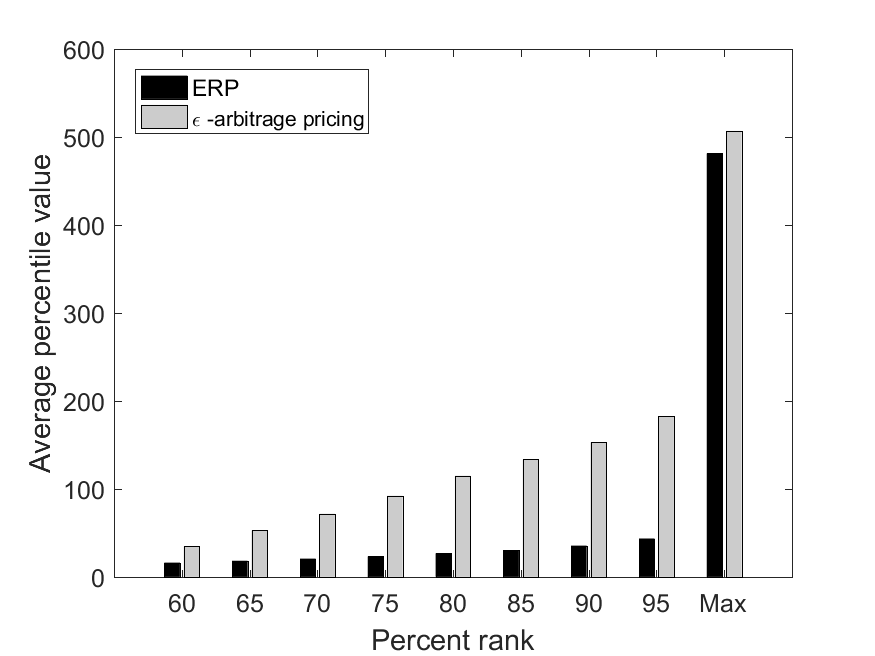}}}
\subfigure[Average (ATM)]{
\resizebox*{50mm}{!}{\includegraphics{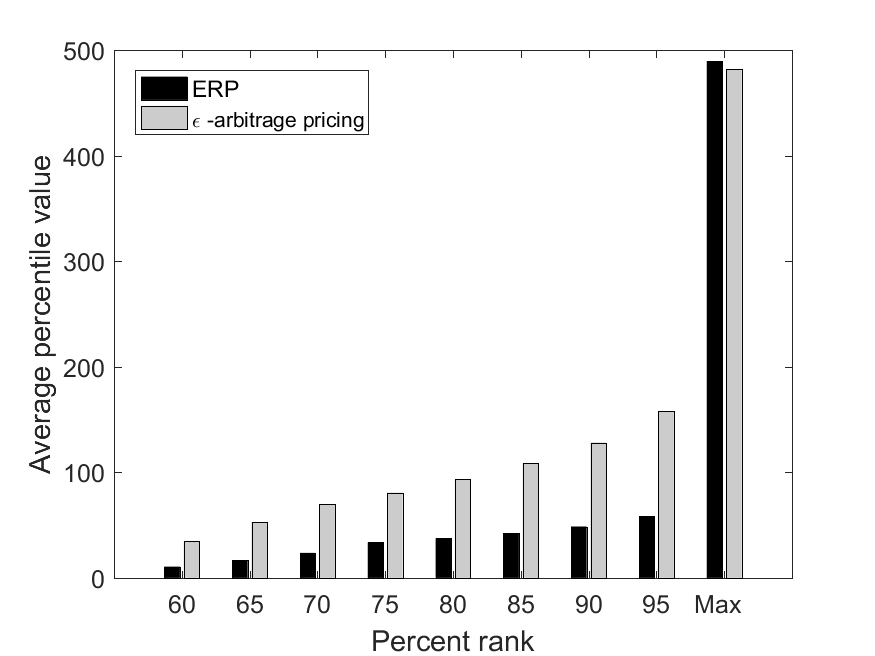}}}
\subfigure[Average (OTM)]{
\resizebox*{50mm}{!}{\includegraphics{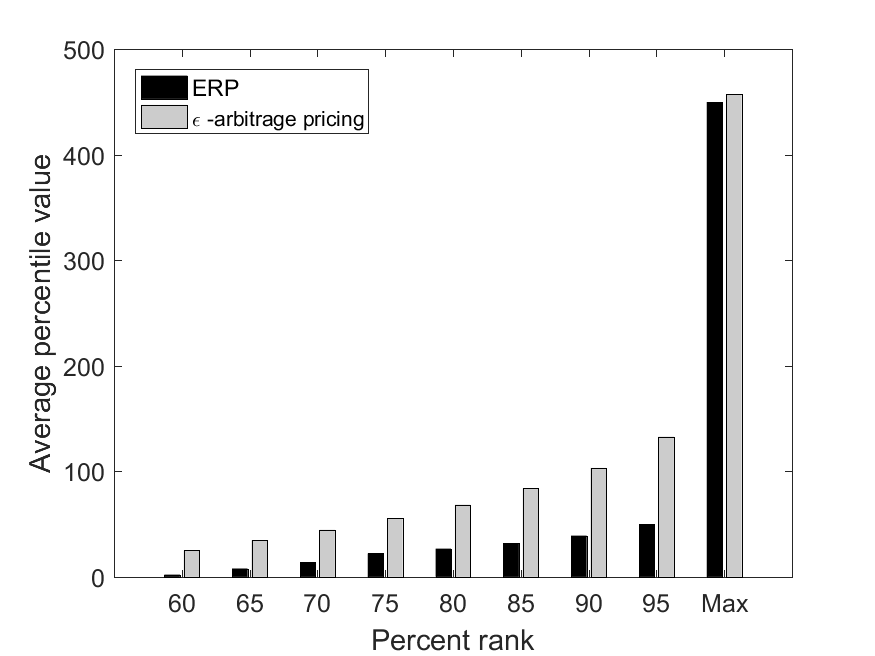}}}

\subfigure[Difference (ITM)]{
\resizebox*{50mm}{!}{\includegraphics{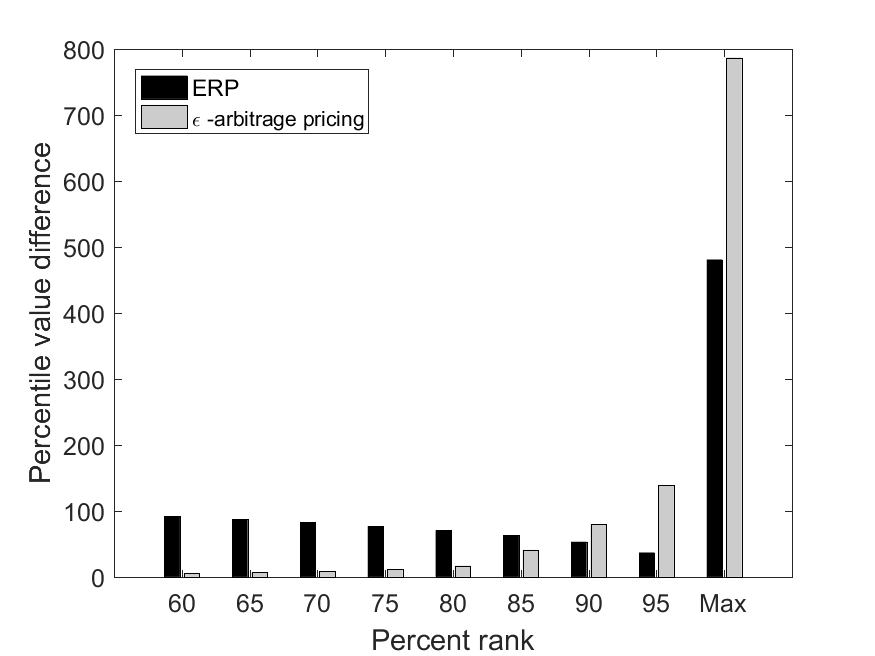}}}
\subfigure[Difference (ATM)]{
\resizebox*{50mm}{!}{\includegraphics{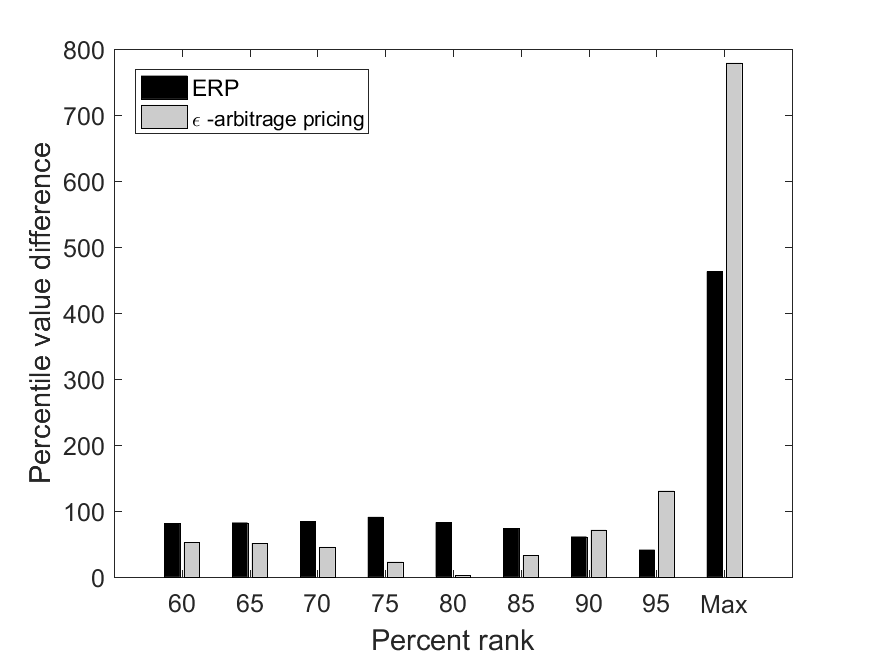}}}
\subfigure[Difference (OTM)]{
\resizebox*{50mm}{!}{\includegraphics{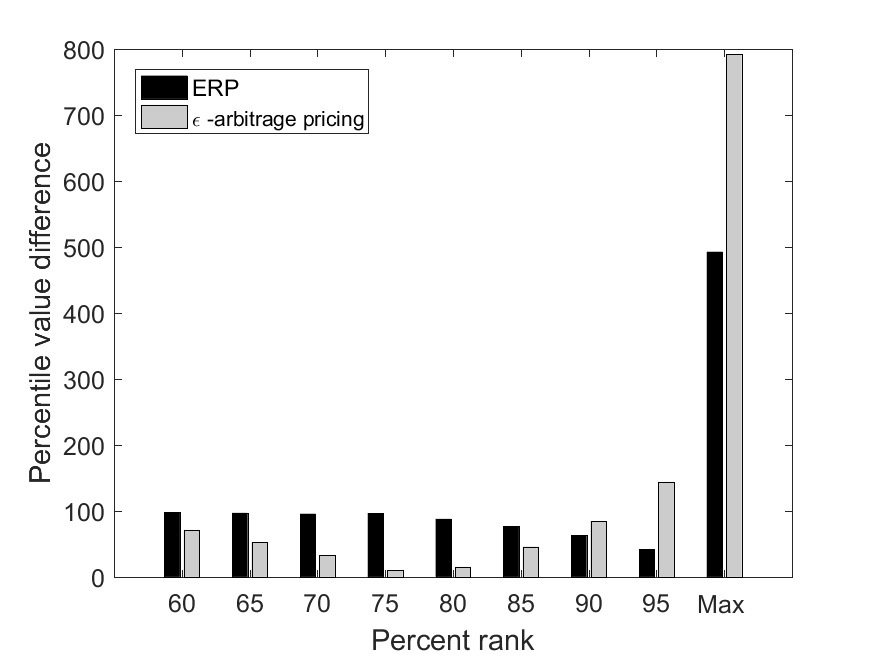}}}

\caption{Comparison of hedging performance achieved under $\epsilon$-arbitrage and equal risk pricing of an European call option with $K=16$ rebalancing periods under a worst-case risk measure that accounts for $\U_1$. (a),(b), and (c) present for different percentile ranks $q$, the average among the $q$-percentile of the loss incurred for the writer and buyer of the ITM, ATM, and OTM options respectively. (d), (e), and (f) present the difference between the same $q$-percentile losses for different percentile rank.
\modified{Note that the \quoteIt{Max} rank refers to the worst-case sample path.}
\label{fig:HedgingResultsEpsilonER}}
\end{minipage}
\end{center}
\end{figure}

We compare also the risk exposure and level of fairness achieved by the transaction prices and hedging strategies produced from the two models. In particular, in our experiments we first simulate a set of $100000$ different sample paths for the risky asset and then for each path we implement the optimal hedging strategy of each model starting with an initial capital that accounts for the transaction price. We record the hedging loss (for both the writer and the buyer) resulting from each sample path and compare different quantiles of the realized losses for both the writer and the buyer. For each quantile level of interest, we compare two different metrics: the average of the quantile value among the writer and buyer’s loss, and the absolute difference between each party's quantile value. Figure \ref{fig:HedgingResultsEpsilonER} presents these metrics for options with different moneyness levels. As seen in Figure \ref{fig:HedgingResultsEpsilonER} (d),(e),(f), the hedging strategy and transaction price suggested by the ERP model leads to lower differences between the two parties' losses when considering quantiles above 90\%. This is clear evidence that ERP is better at sharing the risks among the two parties. It is worth noting that for lower quantiles, $\epsilon$-arbitrage becomes more attractive in this regard which can be explained by the fact that our worst-case risk measures that are used by ERP are insensitive to the performance achieved at lower quantiles. From Figure \ref{fig:HedgingResultsEpsilonER} (a),(b),(c), we see another strength of the ERP model, namely that it does have the ambition of producing optimal risk averse hedging strategies for the two parties together with the ERP. Indeed, this is not the case of the $\epsilon$-arbitrage pricing model, which searches for a single hedging strategy that minimizes the worst-case absolute ``deviation" of the cumulated wealth from the payout.

\begin{table}
\begin{center}
\small
\caption{The prices resulting from ERP with $\U_1'$, $\epsilon$-arbitrage pricing and the Black-Scholes models for options written on an asset with the initial price of 1000, expected annual return of 0.0718, annual standard deviation of 0.1283, strike prices of 950 (ITM), 1000 (ATM), and 1050 (OTM), and one year of maturity. For ERP, the fair price interval is also presented.}
{\begin{tabular}{@{}lcccc|cccc|cccc}
\toprule
&\multicolumn{4}{c}{\textbf{ITM}}&\multicolumn{4}{c}{\textbf{ATM}}&\multicolumn{4}{c}{\textbf{OTM}}\\
\hline
Periods & 16       & 25       & 49       & 100  & 16       & 25       & 49       & 100& 16       & 25       & 49       & 100\\
\hline
$\Gamma$ & 2.63 &	2.70   &	2.79 &	 2.87  & 2.63 &	2.70   &	2.79 &	 2.87 & 2.63 &	2.70   &	2.79 &	 2.87     \\
\hline
FPI-Upper	&	159.20	&	152.95	&	155.88	&	173.90	&	105.72	&	99.36	&	100.81	&	110.90	&	91.73	&	85.21	&	86.69	&	95.45	\\
ERP	&	79.60	&	76.47	&	77.94	&	86.95	&	52.86	&	49.68	&	50.41	&	55.45	&	45.86	&	42.61	&	43.34	&	47.73	\\
FPI-Lower	&	0.00	&	0.00	&	0.00	&	0.00	&	0.00	&	0.00	&	0.00	&	0.00	&	0.00	&	0.00	&	0.00	&	0.00	\\
$\epsilon$-arbitrage price	&	78.40	&	73.60	&	75.20	&	91.20	&	57.60	&	52.80	&	56.00	&	62.40	&	32.00	&	25.60	&	27.20	&	38.40	\\
BS	&	78.80	&	78.80	&	78.80	&	78.80	&	51.15	&	51.15	&	51.15	&	51.15	&	31.17	&	31.17	&	31.17	&	31.17	\\
\hline
\end{tabular}}
\label{table:ERPriceAndEpsilonPrice}
\end{center}
\end{table}

\subsection{Comparison with the Black-Scholes}

In the previous section, we highlighted how the FPI lower bound becomes overly conservative when employing $\U_1$. We believe this explains why the ERP did not show sign of convergence to the Black-Scholes price even when the number of rebalancing periods became large. Given that such a convergence property is appealing when the market is actually based on a GBM process,  in this section, we address  this issue by employing a different uncertainty set that is now motivated by the work of \cite{bernhard2003robust}, namely:
\[\U_2= \left\{r \in \mathbb{R}^K \middle| \sum_{\ell=1}^{sN}r_\ell^2 \in [\sigma^2 sT/S - \Gamma \sqrt{sN}, \sigma^2 sT/S + \Gamma \sqrt{sN}], \forall s \in \{1,...,S\}\right\}\,,\]
with $\Gamma$ small enough so that $\U_2\subset [-1,\infty[^K$. Here we consider the time horizon to be partitioned into $S$ intervals of duration $T/S$, and each interval consists of a set of $N:=K/S$ periods at which the portfolio can be rebalanced. Note that unlike for the set $\U_1$, the set $\U_2$ constrains both the maximum and minimum long term observed deviations.
The main motivation behind the above set is that in the case $\Gamma = 0$, we can expect based on \cite{bernhard2003robust} that the FPI will converge to Black-Scholes price as both $K$ and $N$ converge to infinity. On the other hand, for finite values of $K$ and $N$, the \quoteIt{so-called} budget $\Gamma$ of the set $\U_2$ allows to characterize a meaningful confidence region for the trajectory of the risky asset process. Lastly, as shown in Appendix \ref{proveConditionalBoundedMktRisk}, one can verify that the worst-case risk measure with $\U_2$ satisfies the bounded conditional market risk property and that is Markovian with respect to $\theta_k:=\sum_{\ell=1}^k r_\ell^2$ (see Appendix \ref{proveMarkov} for details).

In our experiments, we set the number of partitions to the square root of the total number of rebalancing periods, i.e. $S=\sqrt{K}$. We calibrate again $ \Gamma$ so that the set $\U_2$ contains $ 95\% $ of simulated price paths. Table \ref{table:ERPriceAndBSPriceCalibrated} presents the equal risk and the fair price intervals against the Black-Scholes prices. From the table, we now see some evidence that the price generated from equal risk pricing is likely to converge to the Black-Scholes price. In particular, one can notice for each type of option that as the total number of rebalancing periods increases, both the upper and lower bounds  evolve monotonically towards the Black-Scholes price, thus driving the equal risk prices  closer and closer to it. 
Unlike with $\U_1$, we see that the FPI lower bounds are now sensitive to the total number of rebalancing periods. Indeed, the lower bound on the total deviation in $\U_2$ allows the buyer of call options to have a less conservative perception of hedging risk. Finally, it is worth noting that the resulting equal risk prices  tend to be slightly higher than the Black-Scholes prices. Indeed, from a practical point of view, this margin can be interpreted as a  \quoteIt{risk premium} on the Black-Scholes price that compensates for the uncertainty that is unaccounted for by the Black-Scholes formula.

\begin{table}
\begin{center}
\small
\caption{The option prices resulting from the equal risk (ERP) and the Black-Scholes (BS) models by using $ \U_2 $. The table also shows the calibrated $\Gamma$ and the upper and lower bounds of the fair price interval.}
{\begin{tabular}{@{}lcccc|cccc|cccc}
\toprule
& \multicolumn{4}{c}{\textbf{ITM}} & \multicolumn{4}{c}{\textbf{ATM}} & \multicolumn{4}{c}{\textbf{OTM}}\\
\hline
Periods & 16       & 49       & 100      & 225  & 16       & 49       & 100      & 225 & 16       & 49       & 100      & 225      \\
\hline
$\Gamma$ & 0.083 &	0.022   &	0.006 &	0.004  & 0.083 &	0.022   &	0.006 &	0.004 & 0.083 &	0.022   &	0.006 &	0.004      \\
\hline
FPI-Upper	&	125.19	&	113.66	&	107.04	&	102.18	&	99.31	&	88.61	&	81.66	&	75.96	&	83.00	&	69.01	&	61.57	&	55.77	\\
ERP	&	87.60	&	84.10	&	83.36	&	83.24	&	49.65	&	54.96	&	55.38	&	54.99	&	41.50	&	36.21	&	35.25	&	34.88	\\
FPI-Lower	&	50.00	&	54.55	&	59.67	&	64.30	&	0.00	&	21.30	&	29.11	&	34.03	&	0.00	&	3.41	&	8.93	&	13.98	\\
BS price	&	78.80	&	78.80	&	78.80	&	78.80	&	51.15	&	51.15	&	51.15	&	51.15	&	31.17	&	31.17	&	31.17	&	31.17	\\
\hline
\end{tabular}}
\label{table:ERPriceAndBSPriceCalibrated}
\end{center}
\end{table}



\begin{figure}
\begin{center}
\begin{minipage}{160mm}
\subfigure[Average (ITM)]{
\resizebox*{50mm}{!}{\includegraphics{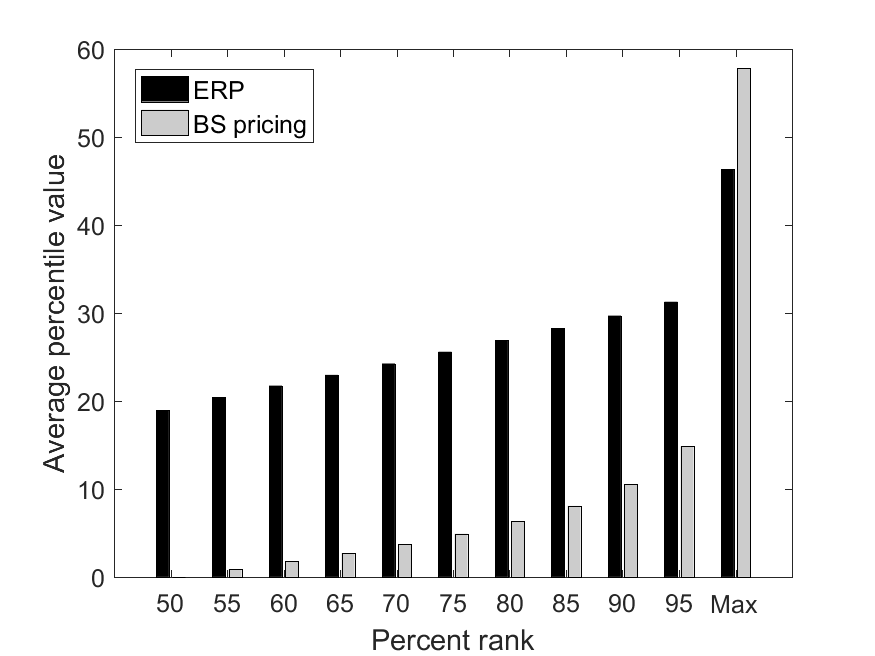}}}
\subfigure[Average (ITM)]{
\resizebox*{50mm}{!}{\includegraphics{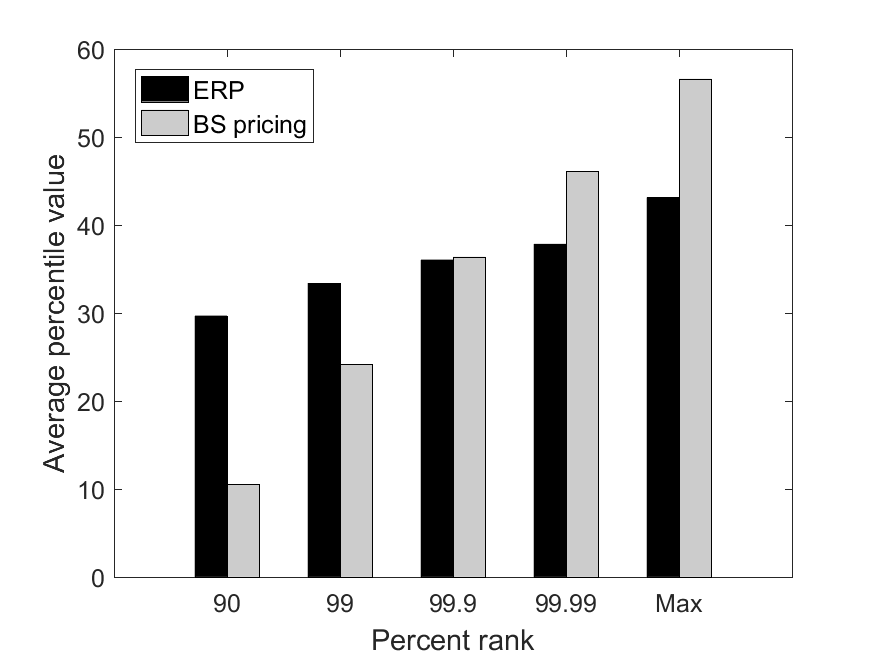}}}
\subfigure[Difference (ITM)]{
\resizebox*{50mm}{!}{\includegraphics{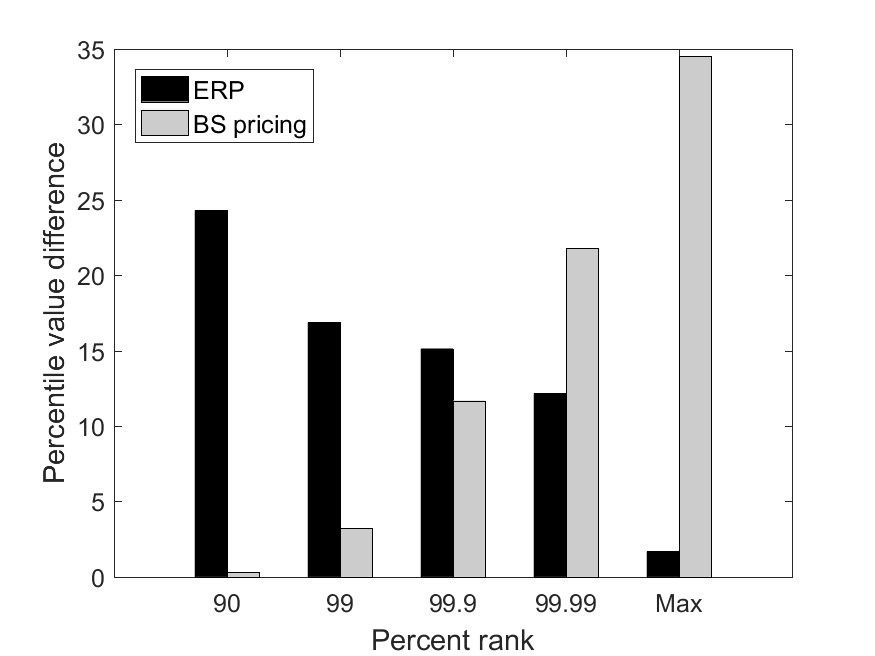}}}

\subfigure[Average (ATM)]{
\resizebox*{50mm}{!}{\includegraphics{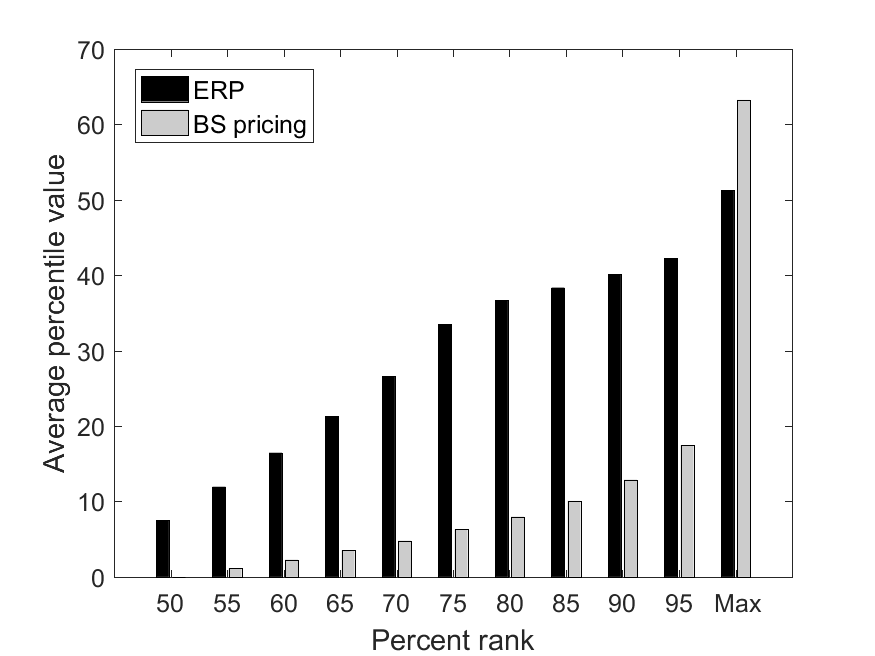}}}
\subfigure[Average (ATM)]{
\resizebox*{50mm}{!}{\includegraphics{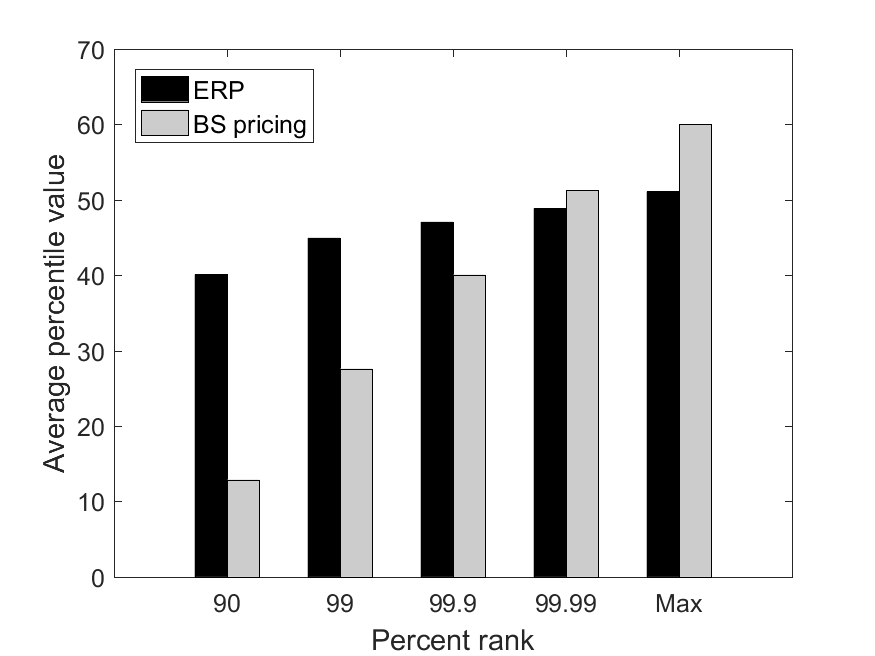}}}
\subfigure[Difference (ATM)]{
\resizebox*{50mm}{!}{\includegraphics{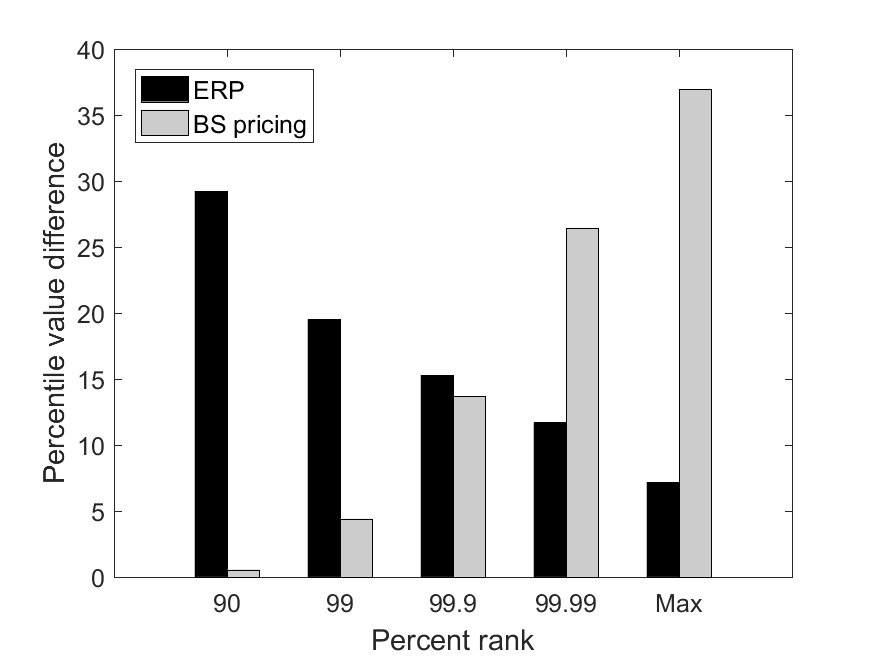}}}

\subfigure[Average (OTM)]{
\resizebox*{50mm}{!}{\includegraphics{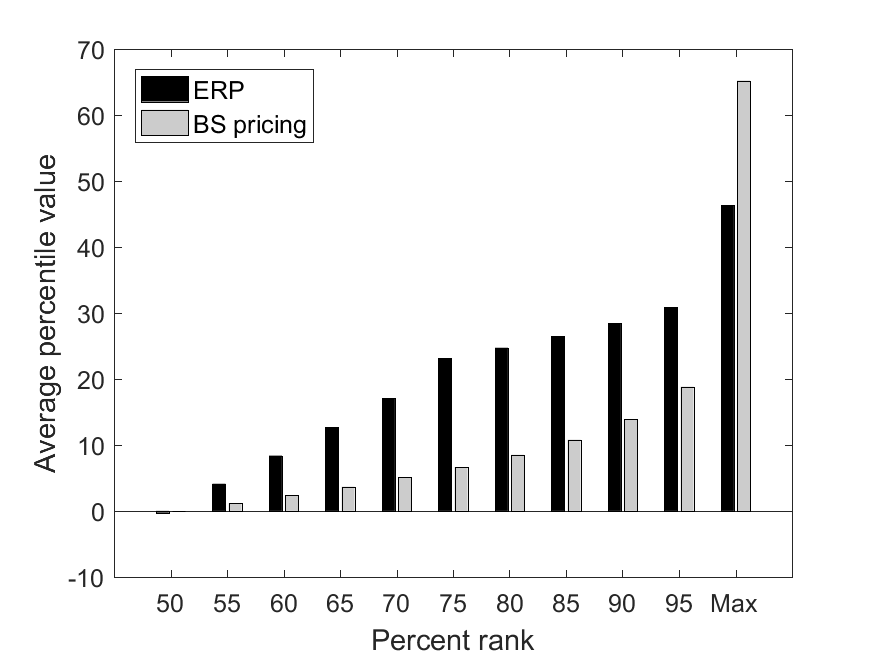}}}
\subfigure[Average (OTM)]{
\resizebox*{50mm}{!}{\includegraphics{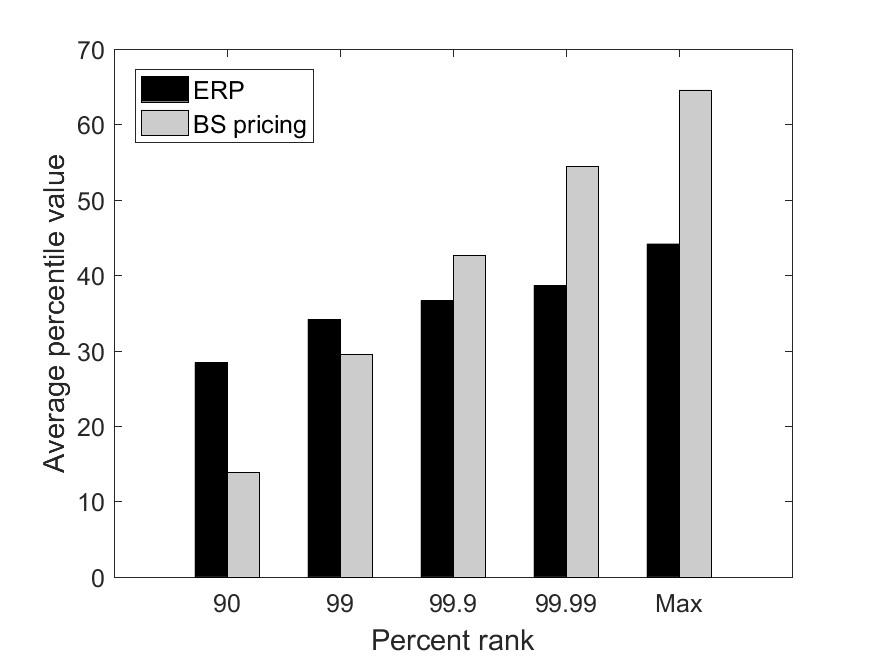}}}
\subfigure[Difference (OTM)]{
\resizebox*{50mm}{!}{\includegraphics{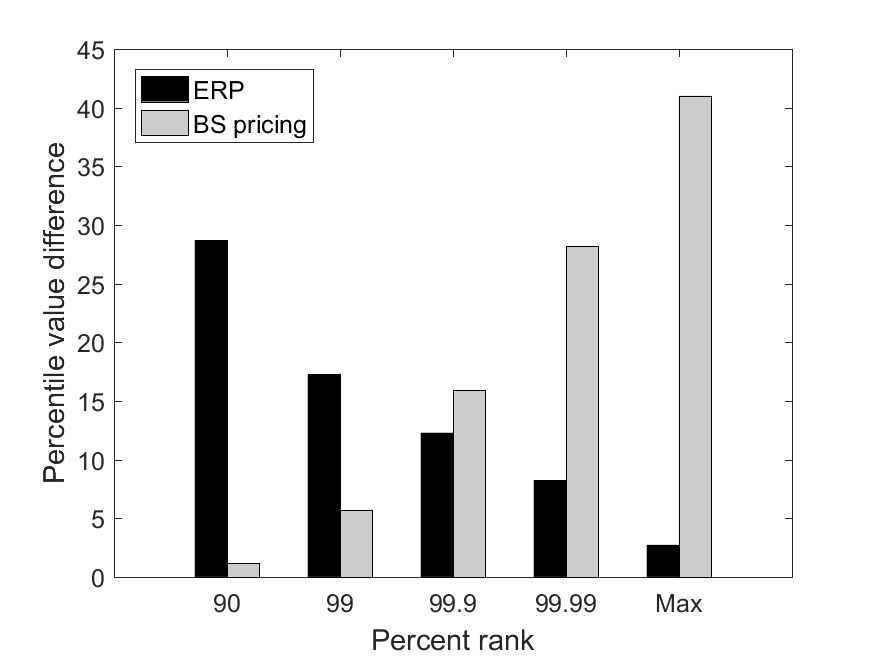}}}

\caption{
Comparison of hedging performance achieved under the Black-Scholes and the equal risk pricing of a European call option with $K=16$ rebalancing periods under a worst-case risk measure that accounts for $\U_2$. (a),(d), and (g) present for different percentile ranks $q$, the average among the $q$-percentile of the loss incurred for the writer and buyer of the ITM, ATM, and OTM options respectively. (b),(e), and (h) presents similar information but focusing on higher percentiles. (c), (f), and (i) present the difference between the same $q$-percentile losses. \modified{Note that the \quoteIt{Max} rank refers to the worst-case sample path.} \label{fig:HedgingResultsHighPartition}}
\end{minipage}
\end{center}
\end{figure}

We provide also in Figure \ref{fig:HedgingResultsHighPartition} the comparison of hedging performances between the equal risk pricing model and the Black-Scholes model, in the case $K=16$. As in the case of comparing with $\epsilon$-arbitrage pricing, we report the performances in terms of both the average and the absolute difference of the writer and buyer’s quantiles of their realized loss distribution under their respective hedging strategy. In particular, here we provide these metrics for quantiles \modified{starting from $99\%$} in order to emphasize what happens at the tail of the loss distributions. \modified{For these figures, the last group of bars is labeled with \quoteIt{Max} to show the worst-case value of the metrics in all samples. This is in line with the type of risk measure that is used in this section.} The results for smaller quantile levels are also provided for average losses (see Figure \ref{fig:HedgingResultsHighPartition} (a),(d),(g)) to present a complete picture. We see that hedging according to the Black-Scholes model actually performs fairly well across a wide range of lower level quantiles, which is not surprising given the market is assumed to follow the GBM assumed by Black-Scholes. Unlike the Black-Scholes model, the ERP model employs a worst-case risk measure that controls the risk in the tail of the loss distributions. As shown in the figures with higher quantile levels, hedging and pricing according to ERP model does indeed become the best scheme when focusing on those regions in terms of both the averages and the differences of risks for the two parties.

\begin{figure}
\begin{center}
\begin{minipage}{150mm}
\subfigure[Writer]{
\resizebox*{70mm}{!}{\includegraphics{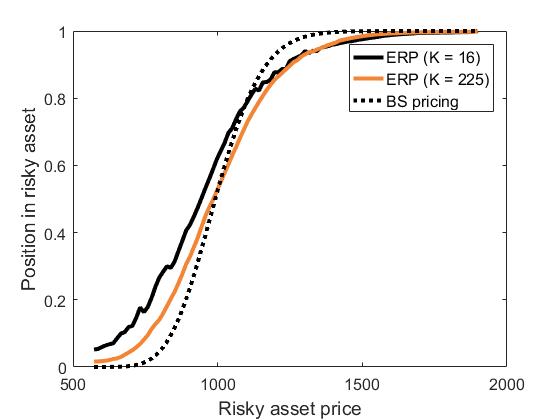}}}
\subfigure[Buyer]{
\resizebox*{70mm}{!}{\includegraphics{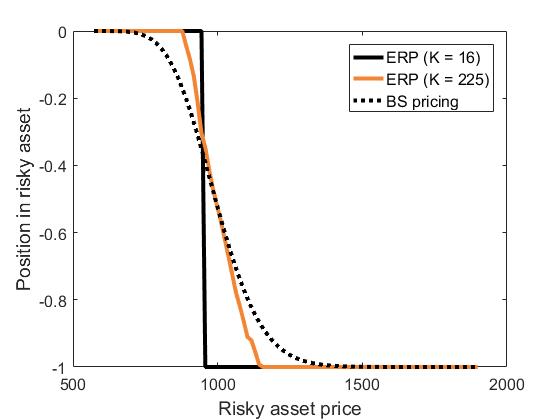}}}



\caption{Comparison of hedging strategies for a European ATM Call option under different number of rebalancing periods. (a) presents the optimal strategies for the writer under the Black-Scholes and the equal risk pricing models for time $k=0$. (b) presents the same for the buyer. \label{fig:HedgingResultsATM2D}}
\end{minipage}
\end{center}
\end{figure}

We continue with Figure \ref{fig:HedgingResultsATM2D}, which presents the hedging strategies at time \modified{$k=0$} proposed by the equal risk pricing and the Black-Scholes pricing for an ATM call option under 16 and 225 rebalancing periods for different asset prices. The first observation is that as $K$ increases, the hedging strategy seems to resemble more the strategy obtained from the Black-Scholes model for both the writer and the buyer. This complements the observation that the ERP appeared to converge to the Black-Scholes price. For lower values of $K$, the hedging strategy for the buyer differs significantly from Black-Scholes hedging because of the larger uncertainty about the risky asset's price process. Specifically, it swings from fully shorting the risky asset to keeping only the risk free asset. The latter strategy becomes optimal because $\Gamma$ is large enough to allow the risky asset to evolve exactly as the risk free one, which also drags the lower bound of FPI to zero as discussed in Section \ref{sec:abprice}. \modified{The second observation is that the hedging strategies of the equal risk model are less sensitive for the writer and more sensitive for the buyer to the variations in the underlying stock price at time $k=0$. In particular, Figure \ref{fig:HedgingResultsATM2D}(a) shows that the ERP model provides a hedging strategy with a lower slope for the writer of the option compared to the Black-Scholes.  On the other side, the opposite is happening for the option buyer. The equal risk hedging strategy is more sensitive to the asset price compared to Black-Scholes strategy.

Finally, we present additional information about the different strategies in Table \ref{tbl:hedgingStrategies}. In particular, the table presents the mean of the average portfolio turnover, computed as $\frac{\sum_{k=1}^{K-1} | \xi_{k+1} - \xi_k |}{K-1}$, over 100000 sample paths for each strategy. This statistic provides evidence that the writer incurs less rebalancing while the opposite is true for the buyer. For completeness, we also present in the table the mean of the average number of shares of the risky asset that are held by each strategy, together with the mean and standard deviation of the respective portfolio values at the maturity of the option.}

\begin{table}
\begin{center}
\small
\caption{\modified{Comparison of the hedging strategies resulting from the equal risk model and the Black-Scholes for $K=225$}}
{\begin{tabular}{@{}lccc|ccc|ccc}
\toprule
\textbf{Measure}&\multicolumn{3}{c}{\textbf{ITM}} & \multicolumn{3}{c}{\textbf{ATM}} & \multicolumn{3}{c}{\textbf{OTM}}\\
\hline
&ERP$_w$&ERP$_b$&BS pricing&ERP$_w$&ERP$_b$&BS pricing&ERP$_w$&ERP$_b$&BS pricing\\
\hline
Average position in risky asset - Mean & 0.7060 & 0.7997 & 0.7501 & 0.5950 & 0.6172 & 0.60954 & 0.4815  & 0.4360 & 0.4606 \\
Average portfolio turnover - Mean & 0.0136 & 0.0200 & 0.0168 & 0.0154 & 0.026 &  0.0203 & 0.0160 & 0.0272 &  0.0210\\
Terminal portfolio value - Mean & 135.15 & 142.43 & 136.68 & 100.05 & 102.45 &  98.88 & 72.39 & 70.00 &  67.91\\
Terminal portfolio value - STD & 113.37 & 131.87 & 122.48 & 101.93 & 118.98 & 110.33 & 89.17 & 102.29 &  95.14\\
\hline
\end{tabular}}
\label{tbl:hedgingStrategies}
\end{center}
\end{table}

\subsection{The Case of American Options}

In this section, we take a further step to benchmark equal risk pricing model against a binomial tree model in the case of American option. For the same reason discussed in the previous section, we assume the worst-case risk measures are defined according to the uncertainty set $\U_2$. Here, we consider put options rather than call options, as the former has attracted more attention in the literature treating American options.
 
The calibration of the uncertainty set $\U_2$, i.e. $\Gamma$ is done in the same fashion as in the previous section. We implement the equal risk model for both the case of with commitment (see Definition \ref{def:ERP_wCommit}) and without commitment (see \modified{Definition} \ref{def:ERP_NoCommit}). We summarize in Table \ref{table:ERPriceAndBSPriceCalibratedAme} all the prices and FPI bounds generated from the model against the option prices generated from the binomial tree model.

As we expected (see Lemma \ref{thm:noComitIsMoreExpansive} ), the equal risk prices with commitment are always smaller or equal to the prices without commitment. We can also confirm that the differences of the prices between the two cases result from the differences in their respective upper bound prices, since their lower bounds are similar. The results also show that as the number of rebalancing periods $K$ increases the equal risk price is getting closer to the binomial tree price. This is happening for all types of options. 
We see that the equal risk price without commitment is larger than the price with commitment by a factor as large as  $4\%$ for ATM and ITM options, and  $10\%$ for OTM options. This non-negligible difference between the prices of the two cases highlights the importance of commitment as a factor to be considered in the negotiation between the two parties regarding the transaction price. This also indicates that the value of the buyer's commitment to an exercise policy is particularly high for an OTM option. 
\begin{table}
\begin{center}
\small
\caption{The option prices resulting from the equal risk pricing model (ERP) under $\U_2$ compared to the binomial tree model (BTM) for American put options. The models with and without commitment are identified respectively as WC and NC.}
{\begin{tabular}{@{}lcccc|cccc|cccc}\toprule
& \multicolumn{4}{c}{\textbf{ITM}} & \multicolumn{4}{c}{\textbf{ATM}} & \multicolumn{4}{c}{\textbf{OTM}}\\
\hline
Periods & 16       & 49       & 100      & 225  & 16       & 49       & 100      & 225 & 16       & 49       & 100      & 225\\
\hline
$\Gamma$ & 0.083 &	0.022   &	0.006 &	0.004  & 0.083 &	0.022   &	0.006 &	0.004 & 0.083 &	0.022   &	0.006 &	0.004      \\
\hline
FPI-Upper-NC	&	134.18	&	117.54	&	114.16	&	106.79	&	103.52	&	87.56	&	84.06	&	77.20	&	75.65	&	62.72	&	59.25	&	53.22	\\
ERP-NC	&	92.09	&	84.73	&	86.85	&	84.10	&	51.76	&	52.95	&	56.49	&	54.25	&	37.83	&	32.95	&	34.11	&	32.16	\\
FPI-Lower-NC	&	50.00	&	51.93	&	59.53	&	61.42	&	0.00	&	18.34	&	28.92	&	31.31	&	0.00	&	3.18	&	8.97	&	11.11	\\
\hline										
FPI-Upper-WC	&	134.18	&	116.14	&	108.54	&	100.43	&	103.52	&	85.12	&	81.36	&	72.95	&	75.65	&	60.84	&	55.76	&	47.14	\\
ERP-WC	&	92.09	&	84.03	&	84.04	&	80.92	&	51.76	&	51.73	&	55.14	&	52.13	&	37.83	&	32.01	&	32.37	&	29.12	\\
FPI-Lower-WC	&	50.00	&	51.93	&	59.53	&	61.42	&	0.00	&	18.34	&	28.92	&	31.31	&	0.00	&	3.18	&	8.97	&	11.11	\\

\hline
BTM price	&	81.52	&	81.19	&	81.13	&	81.21	&	50.36	&	51.41	&	51.02	&	51.21	&	29.00	&	28.73	&	28.68	&	28.85	\\
\hline
\end{tabular}}
\label{table:ERPriceAndBSPriceCalibratedAme}
\end{center}
\end{table}



\begin{figure}
\begin{center}
\begin{minipage}{160mm}
\subfigure[Average (ITM)]{
\resizebox*{50mm}{!}{\includegraphics{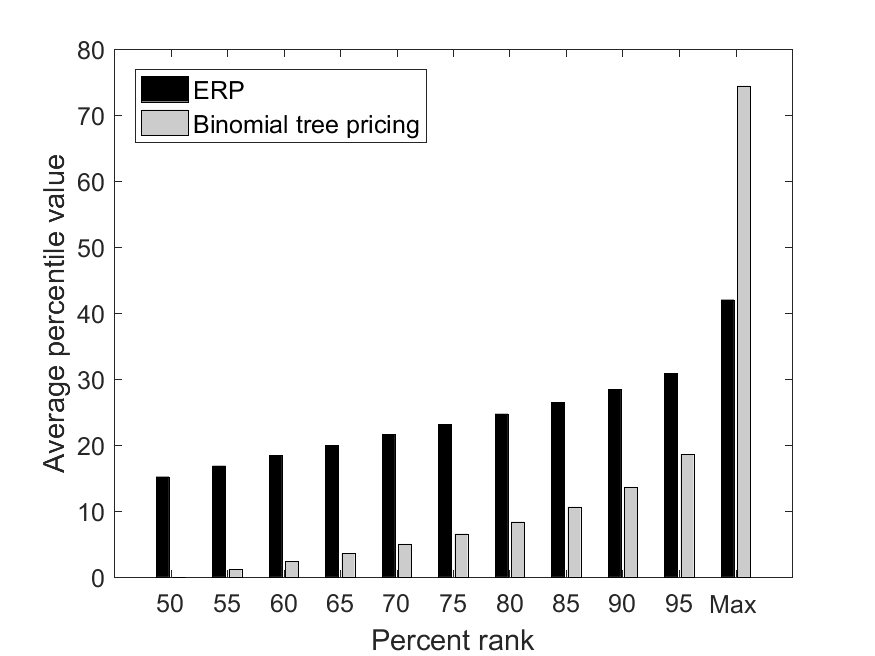}}}
\subfigure[Average (ITM)]{
\resizebox*{50mm}{!}{\includegraphics{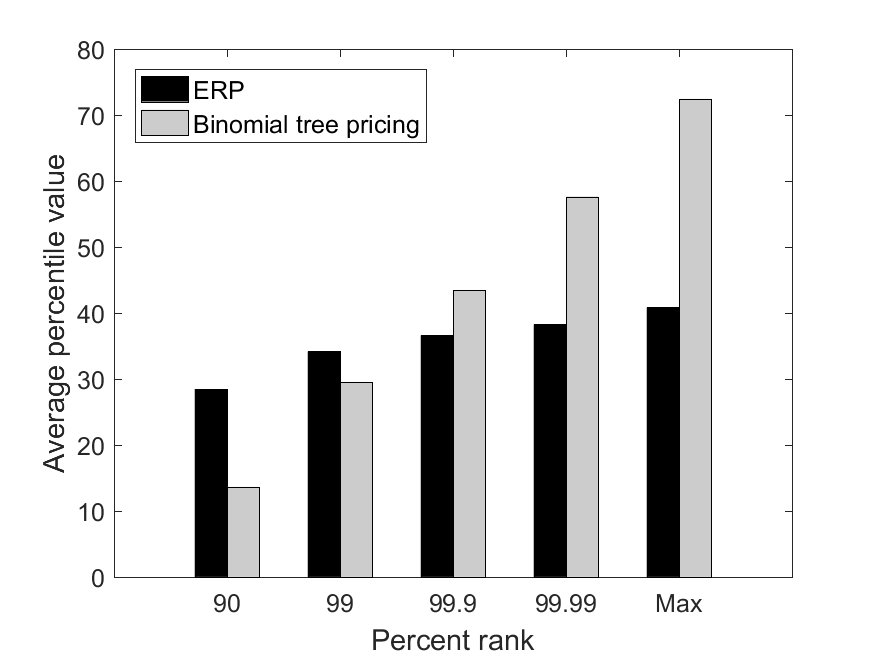}}}
\subfigure[Difference (ITM)]{
\resizebox*{50mm}{!}{\includegraphics{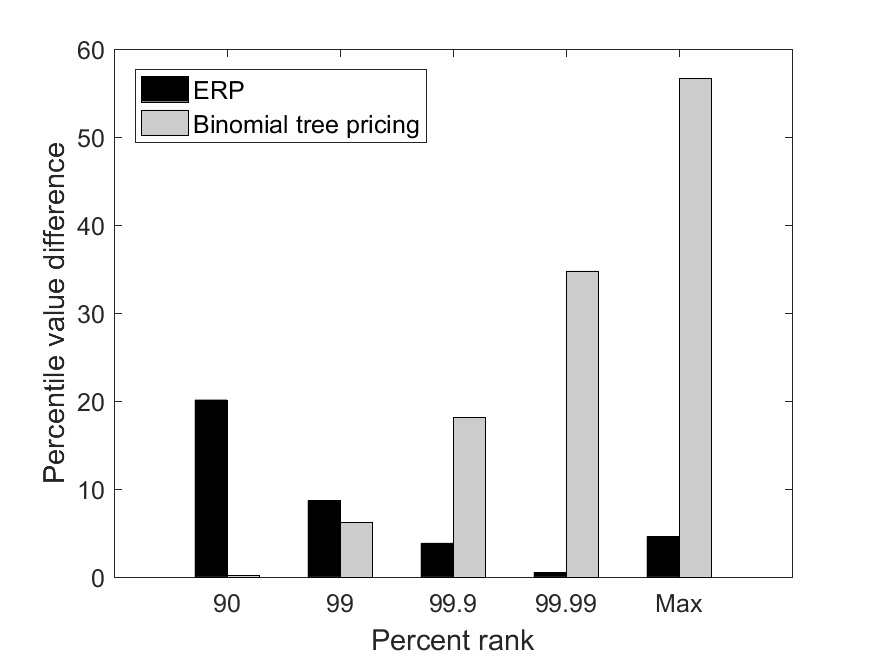}}}

\subfigure[Average (ATM)]{
\resizebox*{50mm}{!}{\includegraphics{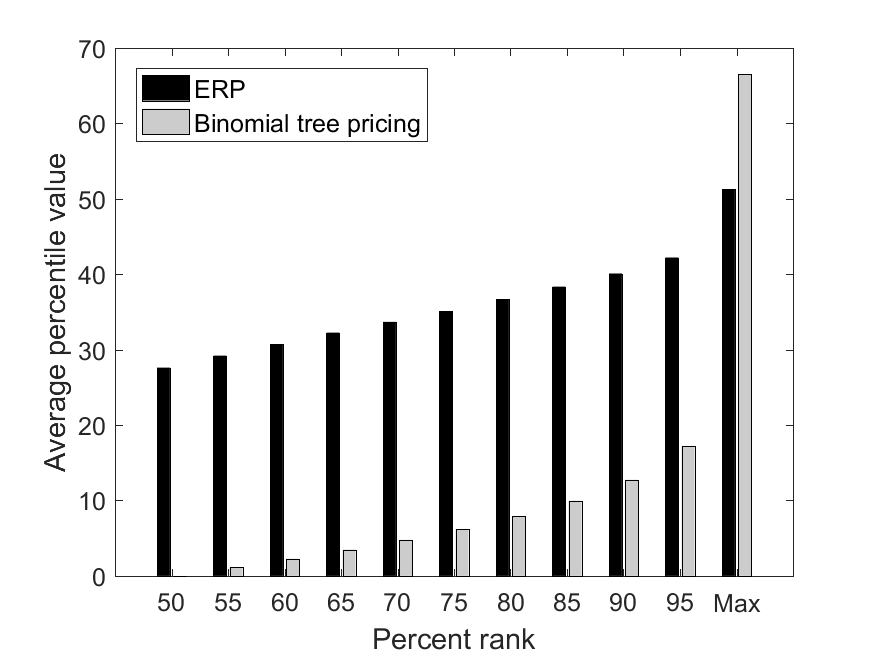}}}
\subfigure[Average (ATM)]{
\resizebox*{50mm}{!}{\includegraphics{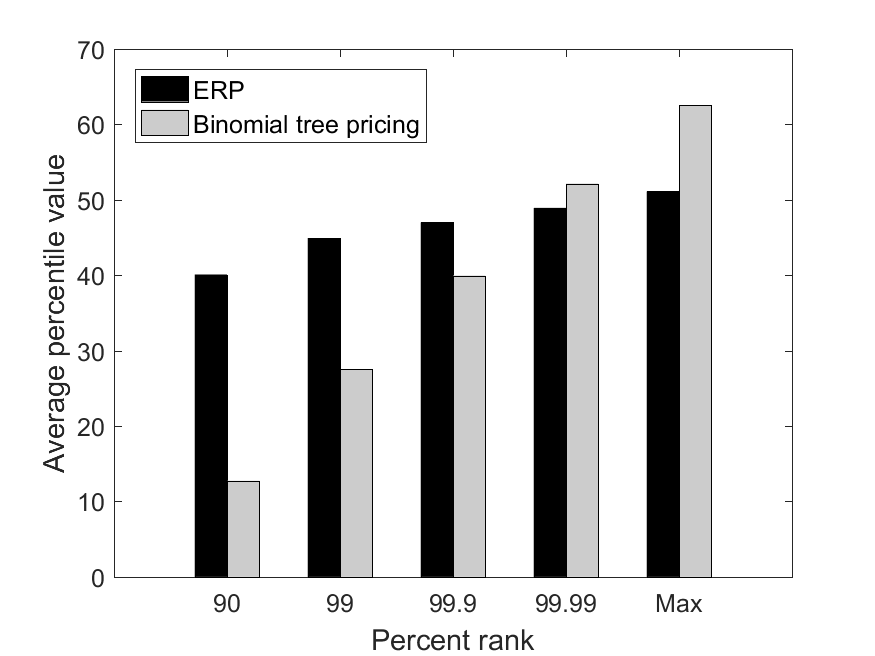}}}
\subfigure[Difference (ATM)]{
\resizebox*{50mm}{!}{\includegraphics{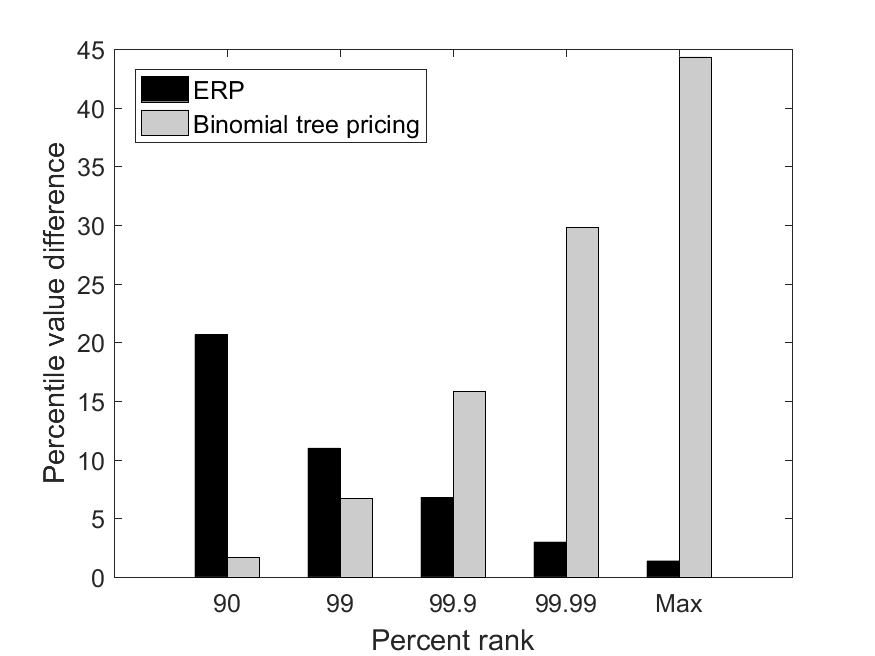}}}

\subfigure[Average (OTM)]{
\resizebox*{50mm}{!}{\includegraphics{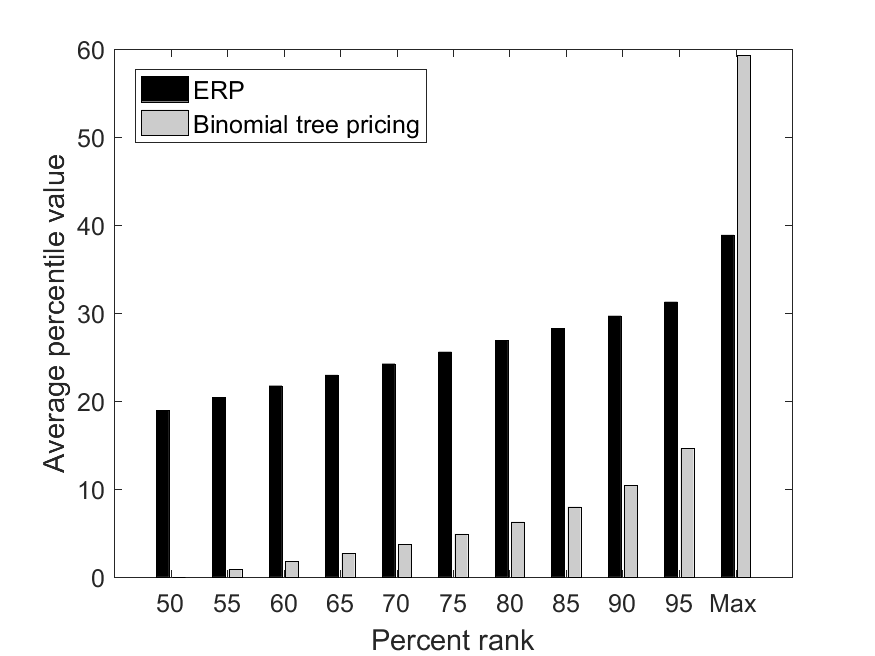}}}
\subfigure[Average (OTM)]{
\resizebox*{50mm}{!}{\includegraphics{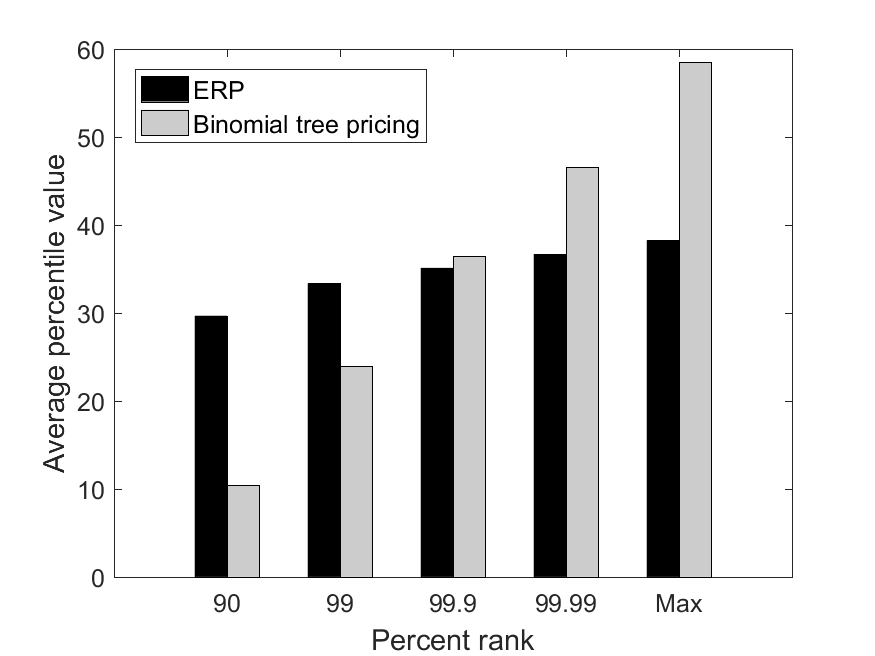}}}
\subfigure[Difference (OTM)]{
\resizebox*{50mm}{!}{\includegraphics{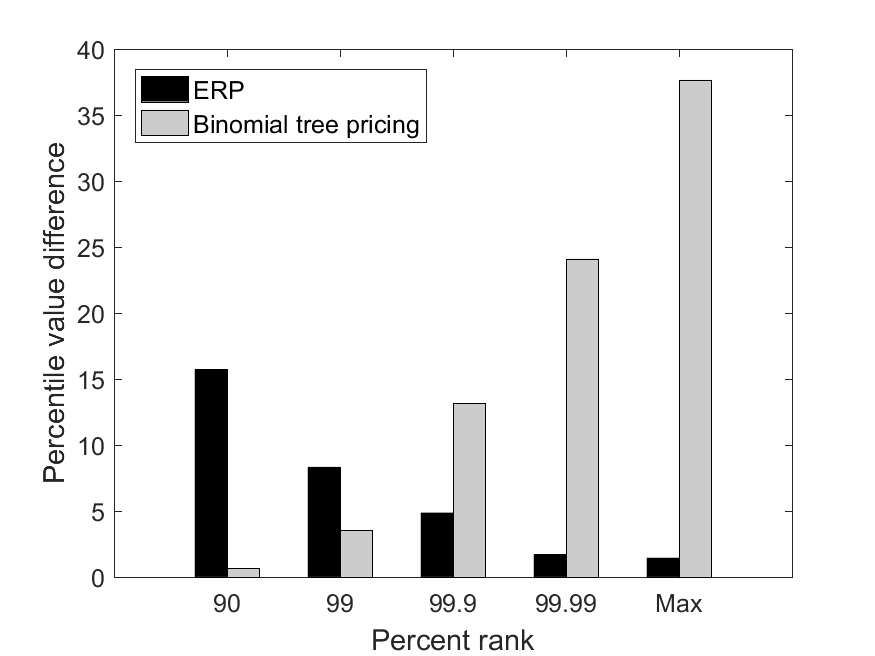}}}
\caption{Comparison of hedging performance achieved under equal risk pricing, with $\U_2$, and a binomial tree model, of an American put option with $K=16$  rebalancing periods. (a),(d), and (g) present for different percentile ranks $q$, the average among the $q$-percentile of the loss incurred for the writer and buyer of the ITM, ATM, and OTM options respectively. (b),(e), and (h) presents similar information but focusing on higher percentiles. (c), (f), and (i) present the difference between the same $q$-percentile losses. \label{fig:HedgingResultsHighPartitionAme}}
\end{minipage}
\end{center}
\end{figure}

In terms of hedging, the equal risk model shows similar results to the case of European options. Figure \ref{fig:HedgingResultsHighPartitionAme} shows that for high quantiles of loss, the equal risk model outperforms the binomial tree model. This is understood from comparing the graphs that focus on the quantiles at the tails of the loss distributions. The higher performance of the equal risk model is specifically more outstanding in terms of the equality of hedging loss for the two sides.  However, having a GBM price process prepares the ground for the binomial tree model to perform well in terms of lower quantiles as shown in figure  \ref{fig:HedgingResultsHighPartitionAme} (a),(d), and (g).

\section{Conclusion} \label{sec:conlusion}

In this paper, we explore the famous problem of pricing and hedging options in an incomplete market under a recently proposed framework called equal risk pricing. Under this framework, the pricing of an option requires that the risk of both sides of the contract be considered in order to make them equal. We consider for the first time the special case of equal risk pricing under convex risk measures for which we show that ERP conveniently reduces to the center of the fair price interval. This price can thus be obtained by solving two dynamic derivative hedging problems, i.e. for the writer and the buyer. By further imposing that the risk measures be one-step decomposable, Markovian, and satisfy a bounded conditional market risk condition, we derive finite dimensional risk averse dynamic programming equations that can be used to solve the discrete time hedging problems for both European and American options. With the latter type of option, the resulting Bellman equations further depend on whether the buyer is willing to commit or not to an exercise strategy \modified{upfront}. All of our results are general enough to accommodate situations where the risk is measured using a worst-case risk measure that considers only a subset of realizations from the outcome space, as typically done in robust optimization.

In our numerical study, we compare the performance of using equal risk pricing with a worst-case risk measure to the performance of $\epsilon$-arbitrage pricing and pricing using the Black-Scholes model in a market that is based on a discretized geometric Brownian motion. In particular, the numerical results confirm that, when using the equal risk price, both the writer and the buyer end up having risks that are more similar and on average smaller than the risks that they would experience by the two other approaches. In addition, by proposing a new uncertainty set inspired from the work of \cite{bernhard2003robust}, we show that the prices generated from equal risk pricing have the potential to converge to Black-Scholes prices as the hedging frequency increases. In the case of pricing American put options, we show how to calculate the value of commitment to an exercise policy, which ranges between 0\% and 10\% for the instances we considered. The evidence seems to indicate that this relative value decreases as the ERP without commitment increases.

\modified{Finally, it is worth mentioning that the results presented in this paper have natural extensions to more general settings than the one that is considered, i.e. with a single underlying asset, zero risk free rate, frictionless market. In some cases, the Bellman equations might need to be extended to account for a larger state space, which is likely to increase the computational efforts needed to identify the equal risk price and the hedging strategies. To circumvent this issue, one might resort to approximate dynamic programming methods. One interesting recent attempt in this direction can be found in \cite{carbonneau2020equal} that proposes a deep reinforcement learning approach to approximate the equal risk price in a variety of market dynamics (including GARCH and Merton jump-diffusion processes) and exotic options with multiple underlying assets.}


\bibliographystyle{plain}
\bibliography{References}

\appendix

\section{Analytical Solutions of One-period Example}\label{app:analyticalSolution}

We recall from \cite{bandi2014robust} that the $\epsilon$-arbitrage model under a worst-case risk measure can be defined as follows for an European option:
\begin{equation}
\begin{aligned}
&\min_{\xi,\price_0}\max_{S_1\in \U} \left \lvert (S_1-K)^+-\price_0-\xi(S_1-S_0)\right \lvert\,,
\end{aligned}\label{eq:epsilonArbitrageOnePeriod}
\end{equation}
where $S_0$ is the initial stock price, $ K $ is the strike price of the option, $S_1$ is the price at the next time period, and $ \U \subseteq \mathbb{R} =[l,u] $. Without loss of generality, we set $ l \le K \le u $ and consider the risk free rate to be zero. 

In the framework of equal risk pricing (ERP), we consider modeling separately the hedging problem of the writer and the buyer. When considering a one period problem, the equal risk model is as follows:
\begin{equation}\nonumber
\begin{aligned}
&\varrho^w(\price_0) :=\min_{\xi_w} \max_{S_1\in\U} (S_1-K)^+-\price_0-\xi_w (S_1-S_0) \\
&\varrho^b(\price_0) :=\min_{\xi_b} \max_{S_1\in\U} -(S_1-K)^++\price_0-\xi_b (S_1-S_0).
\end{aligned}
\end{equation}
The equal risk price is set to be the initial wealth $\price_0$ that leads to $\varrho^w(\price_0) = \varrho^b(\price_0)$. 

\subsection{Analytical Solution for the One Period Equal Risk Model}\label{proveAnalyticalER}
The analytical solution of the one period equal risk model is as follows:
\begin{equation}\nonumber
\begin{aligned}
&\xi_w^*=\frac{u-K}{u-l}, \quad \xi_b^*=\begin{cases}
0, \quad \text{if } S_0<K\\
-1, \quad \text{if } S_0 \ge K
\end{cases}\\
&\price_0^*=(1/2)(S_0-l)\frac{u-K}{u-l}+(1/2)(S_0-K)^+\,.
\end{aligned}
\end{equation}
Considering the writer's side of the equal risk model, since $ (S_1-K)^+ - \xi_w(S_1-S_0) $ is a convex function of $ S_1 $, the maximum in the interval of $ \U=[l,u] $ is at the boundaries, resulting in
\begin{equation}\nonumber
\begin{aligned}
\varrho^w(\price_0)&=\min_{ \xi_w }\max_{S_1 \in \U}(S_1-K)^+ -\price_0-\xi_w(S_1-S_0)\\
&=-\price_0+\min_{ \xi_w }\max \{u-K-\xi_w(u-S_0), -\xi_w(l-S_0)\}\,.
\end{aligned}
\end{equation}
Since the first argument is decreasing in $ \xi_w $ and the second one is increasing, the minimum is at the intersection of the two functions, which results in
\begin{equation}\nonumber
	\begin{aligned}
	&\xi^*_w=\frac{u-K}{u-l}, \quad \varrho^w(\price_0)=-\price_0+(S_0-l)\frac{u-K}{u-l}\,.
	\end{aligned}	
\end{equation}
On the other hand, for the buyer of the option we can show that $\varrho^b(\price_0)=-(S_0-K)^+ +\price_0$ and is achieved using the described hedging strategy $\xi_b^*$. In particular, we can first establish that for all hedging strategies $\xi_b\in\Re$
\begin{equation}\nonumber
\max_{S_1 \in \U} -(S_1-K)^+ +\price_0-\xi_b(S_1-S_0)\ge -(S_0-K)^+ +\price_0\,,
\end{equation}
where we simply use the fact that $S_0\in\U$.

Now, since $ g(y)=-(y-K)^+$ is a concave function of $y$, if $\nabla g(S_0)$ is a supergradient of $g(y)$ at $S_0$ then we have that:
\begin{equation}\nonumber
g(S_1)\leq g(S_0)+\nabla g(S_0)^T(S1-S_0)\,,
\end{equation}
which means that since it can be verified that $\xi_b^*$ is a valide candidate for $\nabla g(S_0)$, we have that:
\begin{equation}\nonumber
\max_{S_1\in\U}-(S_1-K)^+ +\price_0-\xi_b^*(S_1-S_0) \le -(S_0-K)^+ +\price_0\,.
\end{equation}
This proves that $\xi_b^*$ achieves the minimum value of $-(S_0-K)^+ +\price_0$. 

We conclude this discussion with verifying that for $\price_0^*$ indeed leads to the same risk for both the writer and the buyer:
\begin{align*}
\varrho^w(\price_0^*)&\;=\;-\price_0^*+(S_0-l)\frac{u-K}{u-l}\\
&=\;-\left((1/2)(S_0-l)\frac{u-K}{u-l}+(1/2)(S_0-K)^+\right)+\frac{u-K}{u-l}(S_0-l)\\
&=\;\frac{1}{2}(S_0-l)\frac{u-K}{u-l}-(1/2)(S_0-K)^+ \;=\;\price_0^*-(S_0-K)^+=\varrho^b(\price_0^*)\,.
\end{align*}

\subsection{Analytical Solution for the One Period $ \epsilon$-arbitrage Model}
In this section, we will demonstrate that an optimal solution for the $\epsilon$-arbitrage model takes the following form:
\begin{equation}\nonumber
\begin{aligned}
&\xi^*=\frac{u-K}{u-l}, \quad \price_0^*=\frac{u-K}{u-l}\left(S_0-\frac{1}{2}(K+l)\right).
\end{aligned}
\end{equation}
To do so, we will exploit the following lemma, which appears as proposition 3.1.4 in \cite{bertsekas2015convex}
\begin{lemma}\label{thm:BertsekasProp}
	A vector $x^*$ minimizes a convex function $f:R^n \rightarrow R$ over a convex set $\mathcal{X} \subset R^n$ if and only if there exists a subgradient $\nabla f(x^*)$ of $f$ at $x^*$ such that $\nabla f(x^*)^T(x-x^*) \ge 0, \forall x \in \mathcal{X}$.
\end{lemma}
In other words, we will be able to conclude that the $(\xi^*,\price_0^*)$ pair is a minimizer of problem \eqref{eq:epsilonArbitrageOnePeriod}, if we can show that $0$ is a subgradient of the objective function at $(\xi^*,\price_0^*)$. Based on Section 3.1.1 of \cite{bertsekas2015convex}, one can actually show that the set of all subgradients at $(\xi^*,\price_0^*)$ include
\begin{eqnarray}
\left\{ \nabla g\in \Re^2\,\middle|\,\exists \lambda\in\Re^4,\;\lambda\geq 0,\; \sum_{i=1}^4 \lambda_i=1,\;\nabla g = \begin{bmatrix}
-u+S_0 &\;& -l+S_0 &\;& K-S_0 &\;& K-S_0 \\
-1 &\;& -1 &\;& 1 &\;& 1
\end{bmatrix} \lambda \right\}\,.\label{eq:subgradientset}
\end{eqnarray}
One can then readily verify that $0$ is a member of this set using $\lambda_1:=(1/2)(K-l)/(u-l)$, $\lambda_2:=(1/2)(u-K)/(u-l)$, $\lambda_3:=0$, and $\lambda_4:=1/2$ as a certificate.

To provide more details on obtaining the set in \eqref{eq:subgradientset}, we start by recalling that when $f$ is the maximum of $m$ subdifferentiable convex functions $\phi_1,...,\phi_m$:
\begin{equation}\label{eq:maxFunc}
f(x)=\max\{\phi_1(x),...,\phi_m(x)\}, x \in \Re^n\,,
\end{equation}
then a subset of the subdifferential of $f$ can be described as:
\begin{equation}\label{eq:convFunc}
\partial f(x)=\mbox{conv}\{\nabla \phi_j(x)|j\in \mathcal{J}(x)\}\,,
\end{equation}
where $\mathcal{J}(x):=\{ j\in\{1,\dots,m\}\,|\, \phi_j(x)=f(x)\}$, and each $\nabla \phi_j(x)$ is a subgradient of $\phi_j(\cdot)$ at $x$. To obtain the set in \eqref{eq:subgradientset}, we first formulate the objective function in the form of equation \eqref{eq:maxFunc}, we then identify a subgradient $\nabla \phi_j(x^*)$ of each $j\in\mathcal{J}(x^*)$ at our proposed solution $x^*$ to compose the set described in \eqref{eq:convFunc}. 

 \noindent\textbf{Step 1.} We can rewrite the objective function of problem \eqref{eq:epsilonArbitrageOnePeriod} by exploiting a partition of $\U$ as follows:
\begin{equation}\nonumber
\max_{S_1\in\U} | (S_1-K)^+-\price_0-\xi(S_1-S_0)|  = \max \left\{\phi_1(\xi,\price_0), \phi_2(\xi,\price_0), \phi_3(\xi,\price_0), \phi_4(\xi,\price_0)\right\}\,,
\end{equation}
where
\begin{align*}
\phi_1(\xi,\price_0)&:= \max_{S_1 \in [K,u]} | (S_1-K)^+-\price_0-\xi(S_1-S_0)|  = \max_{S_1 \in [K,u]} (S_1-K)-\xi(S_1-S_0)-\price_0\\
\phi_2(\xi,\price_0)&:=\max_{S_1 \in [l,K]} | (S_1-K)^+-\price_0-\xi(S_1-S_0)|  = \max_{S_1 \in [l,K]} -\xi(S_1-S_0)-\price_0\\
\phi_3(\xi,\price_0)&:=\max_{S_1 \in [K,u]} | (S_1-K)^+-\price_0-\xi(S_1-S_0)|  = \max_{S_1 \in [K,u]} -(S_1-K)+\xi(S_1-S_0)+\price_0\\
\phi_4(\xi,\price_0)&:=\max_{S_1 \in [l,K]} | (S_1-K)^+-\price_0-\xi(S_1-S_0)|  = \max_{S_1 \in [l,K]} \xi(S_1-S_0)+\price_0\,.
\end{align*}

 \noindent\textbf{Step 2.} In order to find $\mathcal{J}(x^*)$, we study the maximum of all four functions when $\xi=\xi^*$ and $\price_0=\price_0^*$. Specifically, we have: 
\begin{align*}
&\phi_1(\xi^*,\price_0^*)= \max_{S_1 \in [K,u]} S_1-K-\frac{u-K}{u-l}(S_1-S_0)-\frac{u-K}{u-l}(S_0-\frac{1}{2}(K+l))=\frac{u-K}{u-l}\frac{K-l}{2}\\
&\phi_2(\xi^*,\price_0^*)= \max_{S_1 \in [l,K]} -\frac{u-K}{u-l}(S_1-S_0)-\frac{u-K}{u-l}(S_0-\frac{1}{2}(K+l))=\frac{u-K}{u-l}\frac{K-l}{2}\\
&\phi_3(\xi^*,\price_0^*)=\max_{S_1 \in [K,u]}  -S_1+K+\frac{u-K}{u-l}(S_1-S_0)+\frac{u-K}{u-l}(S_0-\frac{1}{2}(K+l))=\frac{u-K}{u-l}\frac{K-l}{2}\\
&\phi_4(\xi^*,\price_0^*)=\max_{S_1 \in [l,K]}\frac{u-K}{u-l}(S_1-S_0)+\frac{u-K}{u-l}(S_0-\frac{1}{2}(K+l))=\frac{u-K}{u-l}\frac{K-l}{2}\,,
\end{align*}
where we exploited the fact that the functions that are maximized are either non-decreasing for the case of $\phi_1$ and $\phi_4$ or non-increasing for $\phi_2$ and $\phi_3$. In each case, the maximum is achieved at  $S_1^*=u$ for $\phi_1$, $S_1^*=l$ for $\phi_2$, and $S_1^*=K$ for $\phi_3$ and $\phi_4$. Based on this conclusion, we get the following four subgradients:
\begin{align*}
\nabla \phi_1(\xi^*,\price_0^*) &:= \begin{bmatrix}S_0-u\\-1\end{bmatrix}& \nabla \phi_2(\xi^*,\price_0^*) &:=\begin{bmatrix}S_0-l\\-1\end{bmatrix}\\ \nabla \phi_3(\xi^*,\price_0^*) &:=\begin{bmatrix}K-S_0\\1\end{bmatrix}& \nabla \phi_4(\xi^*,\price_0^*) &:=\begin{bmatrix}K-S_0\\1\end{bmatrix}\,.
\end{align*}
 This completes our proof.
 
\section{Proofs for Section \ref{sec:Motivation}}

\subsection{Proof of Proposition \ref{pr:TrInv}}\label{app:midPrice}

	This proof mainly relies on the translation invariance property together with the following property of $ \mathcal{X}(\price_0) $:
	\begin{equation*}
	\begin{aligned}
	\mathcal{X}(\price_0)&=\left\{X\,\middle|\,\exists \xi_s,\,\exists c\in\Re, \quad X_t=\price_0+ \int_{0}^{t}\xi_sdS_s \geq c\,; \forall\,t\in[0,\,T]\right\}\\\
	&=\price_0+\left\{X'\,\middle|\,\exists \xi_s,\,\exists c\in\Re, \quad  X'_t=\int_{0}^{t}\xi_sdS_s  \geq c\,; \forall\,t\in[0,\,T] \right\}\\
	&=\price_0+\mathcal{X}(0)\,,
	\end{aligned}
	\end{equation*}
	where $\price_0+\mathcal{X}(0)$ refers to a set addition. These two properties can be used to show that both
	\begin{align}
	\varrho^w(\price_0)&=\inf_{X\in\mathcal{X}(0)+\price_0} \rho^{w}(F(S_T,Y_T)-X_T)\notag\\
	&=\inf_{X\in\mathcal{X}(0)} \rho^{w}(F(S_T,Y_T)-X_T-\price_0)\notag\\
	&=\inf_{X\in\mathcal{X}(0)} \rho^{w}(F(S_T,Y_T)-X_T)-\price_0\notag\\
	&=\inf \{s | \inf_{X\in\mathcal{X}(0)} \rho^{w}(F(S_T,Y_T)-X_T) \le s\}  - \price_0\notag\\
	&=\inf \{s | \varrho^w(s) \le 0 \}-\price_0 = 	\price_0^w - \price_0\,,\label{eq:proofMidprice1}
	\end{align}
	and similarly,
	\begin{align}
	\varrho^b(\price_0)&=\inf_{X\in\mathcal{X}(0)-\price_0} \rho^{b}(-X_T-F(S_T,Y_T))\notag\\
	&=\inf_{X\in\mathcal{X}(0)} \rho^{b}(-X_T-F(S_T,Y_T))+\price_0\notag\\
	&= \inf \{s | \varrho^b(0) \le s \} + \price_0\notag\\
	&= \inf \{s | \varrho^b(-s) \le 0 \} + \price_0\notag\\
	&= - \sup \{-s | \varrho^b(-s) \le 0 \} + \price_0\notag\\
	&= 	-\price_0^b + \price_0\,. \label{eq:proofMidprice2}
	\end{align}
	
	Hence, we can obtain our result by verifying both directions of the biconditional logical connective. First, given that an equal risk price exists, say $\price_0^*\in\Re$, it must be that both $\varrho^b(\price_0^*)$ and $\varrho^b(\price_0^*)$ are members of $\Re$. This necessarily implies that $\price_0^w\in\Re$ and $\price_0^b\in\Re$ thus that the fair price interval is bounded. Conversely, if the fair price interval is bounded, then one can verify that the midpoint $\price_0^*:=(\price_0^w+\price_0^b)/2$ does satisfy the equal risk price condition:
	\[\varrho^w(\price_0^*)=\price_0^w - \price_0^* =\price_0^w/2 - \price_0^b/2 = -\price_0^b + \price_0^* =  \varrho^b(\price_0^*)\,.\]
	Furthermore, this midpoint can be calculated as:
	\begin{align*}
	\price_0^*&=(1/2)(\inf \{\price_0 | \varrho^w(\price_0) \le 0 \}-\sup \{-s | \varrho^b(-s) \le 0 \})\\
	&=(1/2)(\price_0^w+\price_0^b)=(1/2)(\varrho^w(0) - \varrho^b(0)))\,,
	\end{align*}
	following exactly the same arguments as in \eqref{eq:proofMidprice1} and \eqref{eq:proofMidprice2}.
	This completes our proof. 

\subsection{Proof of Lemma \ref{thm:ERPisNoArbitrage}}

The proof follows directly from Property 2 in \cite{xu2006risk}. In particular, we have that $\price_0^w$ is bounded above by the super-hedging price and $\price_0^b$ is bounded below by the sub-hedging price. Hence, since $\price_0^b\leq \price_0^w$, we must have that the fair price interval is a subset of the no-arbitrage interval. This lets us conclude that the equal risk price is also a member of the no-arbitrage interval.



\section{Proofs for Section \ref{sec:EuropeanAmerican}}

\subsection{Proof of Proposition \ref{prop:bellman}}\label{app:EurOptDPproof}

We focus on providing the arguments supporting the claims for the writer model as these are analogous for the buyer model. In doing so, we will closely follow the theory presented in \cite{pichler2018risk}. We start by constructing the so-called additive preference system $\{\mathcal{R}_{k,l}\}_{(k,l)\in \{0,\dots,K\}^2:k<l}$ (a.k.a. a dynamic risk measure) based on $\rho^w$, where each $\mathcal{R}_{k,l}:\mathcal{L}_p(\Omega,\mathcal{F}_{k},\mathbb{P})\times \mathcal{L}_p(\Omega,\mathcal{F}_{k+1},\mathbb{P})\times\cdots\times \mathcal{L}_p(\Omega,\mathcal{F}_{l},\mathbb{P})$ takes the form:
\[\mathcal{R}_{k,\ell}(Z_k,Z_{k+1},\dots,Z_{\ell}):=\rho_k^w(\rho_{k+1}^w(\cdots\rho_K^w(\sum_{k'=k}^\ell Z_{k'})\cdots)\,,\,\forall\,0\leq k < \ell \leq K\,.\]
Based on this definition of $\mathcal{R}_{k,\ell}$ it is easy to see that:
\[	\varrho^w(0)=\inf_{X\in\mathcal{X}(0)} \mathcal{R}_{0,K}(0,0,\dots,0,F(S_T,Y_T)-X_K) \,.
\]
Given that $\rho^w$ is one-step decomposable, it is easy to show that $\mathcal{R}$ is both \quoteIt{Monotone} and \quoteIt{Recursive} (see definitions 2.3 and 4.1 in \cite{pichler2018risk}). In particular, for monotonicity we have that:
\begin{align*}
\forall (Z_k,\dots,Z_\ell), (Z_k',\dots,Z_\ell'), \;Z_{k'}&\geq Z_{k'}' \mbox{ a.s. } \forall\,k'=k,\dots,\ell \Rightarrow\\
\mathcal{R}_{k,\ell}(Z_k,\dots,Z_\ell)&=\rho_k^w(\rho_{k+1}^w(\cdots\rho_K^w(\sum_{k'=k}^\ell Z_{k'})\cdots)\\
&\geq \rho_k^w(\rho_{k+1}^w(\cdots\rho_K^w(\sum_{k'=k}^\ell Z_{k'}')\cdots) \\
&= \mathcal{R}_{k,\ell}(Z_k',\dots,Z_\ell')\,,
\end{align*}
given that each $\rho_k^w$ is monotone. On the other hand, for recursivity, we have
\begin{align*}
\forall Z_k,\dots,Z_\ell &\Rightarrow\\
\mathcal{R}_{k,\ell}(Z_k,\dots,Z_\ell)
&= \rho_k^w(\rho_{k+1}^w(\cdots\rho_{v-1}^w(\rho_v^w(\cdots\rho_K^w(\sum_{k'=k}^\ell Z_{k'})\cdots)\\
&= \rho_k^w(\rho_{k+1}^w(\cdots\rho_{v-1}^w(\sum_{k'=k}^{v-1} Z_{k'}+\rho_v^w(\cdots\rho_K^w(\sum_{k'=v}^\ell Z_{k'})\cdots)\\
&= \rho_k^w(\rho_{k+1}^w(\cdots\rho_{v-1}^w(\sum_{k'=k}^{v-1} Z_{k'}+\mathcal{R}_{v,\ell}(Z_v,\dots,Z_\ell))\cdots)\\
&= \rho_k^w(\rho_{k+1}^w(\cdots\rho_{v-1}^w(\rho_v^w(\cdots\rho_K^w(\sum_{k'=k}^{v-1} Z_{k'}+\mathcal{R}_{v,\ell}(Z_v,\dots,Z_\ell))\cdots)\\
&= \mathcal{R}_{k,v}(Z_k,\dots,Z_{v-1},\mathcal{R}_{v,\ell}(Z_v,\dots,Z_\ell)))\,,
\end{align*}
where we exploited monotonicity, the definition of $\mathcal{R}_{v,\ell}$, and the conditional transition invariance of $\rho_k^w$ for all $k=v,\dots,K$.

Having verified these conditions, Proposition 3.1 and the discussion that follows in Section 4 of \cite{pichler2018risk} allows us to conclude that:  
\begin{align*}
\varrho^w(0)&=\inf_{X\in\mathcal{X}(0)} \mathcal{R}_{0,K}(0,0,\dots,0,F(S_T,Y_T)-X_K)\\ 
&=\inf_{\xi_0}\mathcal{R}_{0,1}(0,\inf_{\xi_1}\mathcal{R}_{1,2}(0,\cdots\inf_{\xi_{K-1}}\mathcal{R}_{K-1,K}(0,F(S_T,Y_T)-X_K)\cdots)\,.
\end{align*}
Hence, $\varrho^w(0)=\bar{V}_0^w(0,\omega)$ where
\begin{align*}
&\bar{V}^w_k(X_k,\omega):=\inf_{\xi_k}\rho^w_k(\bar{V}^w_{k+1}(X_k+\xi_k \Delta S_{k+1}),\omega)\\
&\bar{V}^w_K(X_K,\omega):=F(S_T(\omega),Y_T(\omega))-X_K\modified{(\omega)}\,.
\end{align*}
Furthermore, the set of optimal policies for problem \eqref{equation:riskMeasures1} must contain the following hedging policies: 
\begin{align*}
\bar{\xi}_k^{w*}(X_k,\omega)\in\arg\min_{\xi_k}\rho^w_k(\bar{V}^w_{k+1}(X_k+\xi_k \Delta S_{k+1}),\omega),\,\forall\,k=0,\dots,K-1\,.
\end{align*}
Yet, by conditional translation invariance, we know that: 
\[\bar{V}^w_K(X_K,\omega)=V^w_K(\omega)-X_K\modified{(\omega)}\,,\]
and recursively that:
\begin{align*}
\bar{V}^w_k(X_k,\omega)&=\inf_{\xi_k}\rho^w_k(\bar{V}^w_{k+1}(X_k+\xi_k \Delta S_{k+1}),\omega)\\
&=\inf_{\xi_k}\rho^w_k(V^w_{k+1}-X_k-\xi_k \Delta S_{k+1},\omega)\\
&=\inf_{\xi_k}\rho^w_k(V^w_{k+1}-\xi_k \Delta S_{k+1},\omega)-X_k\modified{(\omega)}=V^w_{k}(\omega)-X_k\modified{(\omega)}\,.
\end{align*}
Hence, $\varrho^w(0)=\bar{V}_0^w(0,\omega)=V_0^w(\omega)$.
A similar reasoning confirm that the set of optimal policies for problem \eqref{equation:riskMeasures1} contains for all $k=0,\dots,K-1$:
\begin{align*}
\bar{\xi}_k^{w*}(X_k,\omega)&\in\arg\min_{\xi_k}\rho^w_k(\bar{V}^w_{k+1}(X_k+\xi_k \Delta S_{k+1}),\omega)\\
&=\arg\min_{\xi_k}\rho^w_k(V^w_{k+1}-\xi_k \Delta S_{k+1},\omega)-X_k\modified{(\omega)} \\
&=\arg\min_{\xi_k}\rho^w_k(V^w_{k+1}-\xi_k \Delta S_{k+1},\omega)\,,
\end{align*}
thus equivalent to $\xi_k^{w*}(\omega)$.

\subsection{\modified{
Dynamic Programming Equations for the Case Non-translation Invariant Risk Measures}}\label{app:DP:nonTransInv}

In the case where $\rho^w$ and $\rho^b$ do not satisfy the translation invariance property, the bisection method described in Remark \ref{rem:nonTI} relies on computing the value of $\varrho^w(\price_0)$ and $\varrho^b(\price_0)$ for any value of $\price_0$. Similar arguments as used in sections \ref{subsec:EUoption} and \ref{sec:amOption} to identify dynamic programming equations, which now depend on accumulated wealth, when the risk measures are one-step decomposable and Markovian. 

In particular, focusing on the case of the writer's problem associated to a European option, by following the steps in Section \ref{app:EurOptDPproof}, one can simply work with the following unreduced value functions:
\begin{subequations}
\begin{align}
&\bar{V}^w_k(X_k,\omega):=\inf_{\xi_k}\rho^w_k(\bar{V}^w_{k+1}(X_k +\xi_k \Delta S_{k+1}),\omega)\,,\,k=0,\dots,K-1\\
&V^w_K(X_K,\omega):=F(S_{{T}}(\omega),Y_{{T}}(\omega)) - X_K(\omega)\,,
\end{align}
\end{subequations}
serving the purpose of computing $\varrho^w(\price_0) = \bar{V}^w_0(\price_0)$.

Furthermore, by exploiting the Markovian risk measure assumption, one can easily reduce the representation to the following Bellman equations:
\[\tilde{V}^w_K(X_K,S_K,Y_K,\theta_K^w):= F(S_K,Y_K) - X_K\,,\]
and recursively
\[\tilde{V}^w_k(X_k,S_k,Y_k,\theta_k^w):= \inf_{\xi_k} \bar{\rho}_k(\tilde{V}(X_k + \xi_k\Delta S_{k+1},S_k+\Delta S_{k+1},Y_k+\Delta Y_{k+1},f_k(\theta_k^w)),\theta_k)\,.\]
Then, considering that
\[\bar{V}^w_K(X_K,\omega)=\tilde{V}^w_K(X_K(\omega),S_K(\omega),Y_K(\omega),\theta_K^w(\omega))\,,\]
and recursively that if $\bar{V}^w_{k+1}(X_{k+1},\omega)=\tilde{V}^w_{k+1}(X_{k+1}(\omega),S_{k+1}(\omega),Y_{k+1}(\omega),\theta_{k+1}^w(\omega))$, then we have that:
\begin{align*}
&\bar{V}^w_k(X_k,\omega)=\inf_{\xi_k}\rho^w_k(\bar{V}^w_{k+1}(X_k +\xi_k \Delta S_{k+1}),\omega)\\
&=\inf_{\xi_k}\rho^w_k(\tilde{V}^w_{k+1}(X_k +\xi_k \Delta S_{k+1},S_{k+1},Y_{k+1},\theta_{k+1}^w),\omega)\\
&=\inf_{\xi_k}\bar{\rho}^w_k(\bar{\Pi}_k(\tilde{V}^w_{k+1}(X_k +\xi_k \Delta S_{k+1},S_{k+1},Y_{k+1},\theta_{k+1}^w)),\theta_k^w(\omega))\\
&=\inf_{\xi_k}\bar{\rho}^w_k(\tilde{V}^w_{k+1}(X_k(\omega) +\xi_k \Delta S_{k+1},S_k(\omega)+\Delta S_{k+1},Y_k(\omega)+\Delta Y_{k+1},f(\theta_{k}^w)),\theta_k^w(\omega))\\
&=\tilde{V}_k(X_k(\omega),S_k(\omega),Y_k(\omega),\theta_k(\omega))\,.
\end{align*}

Now, we see that $\varrho^w(\price_0)=\bar{V}^w_0(\price_0)=\tilde{V}^w_0(\price_0,S_0,Y_0,\theta_0^w)$. In the case of the buyer, similar derivations lead to the Bellman equations:
\[\tilde{V}^b_K(X_K,S_K,Y_K,\theta_K^b):= -F(S_K,Y_K) - X_K\]
and 
\[\tilde{V}^b_k(X_k,S_k,Y_k,\theta_k^b):= \inf_{\xi_k} \bar{\rho}_k^b( \tilde{V}^b_{k+1}(X_k-\xi_k\Delta S_{k+1} ,S_k+\Delta S_{k+1},Y_k+\Delta Y_{k+1},f_k(\theta_k^b)),\theta_k^b)\,,\]
which can be used to compute $\varrho^b(\price_0)=\tilde{V}^b_0(\price_0,S_0,Y_0,\theta_0^b)$.



\erickdelete{
If the risk measures are not translation invariant, the value functions of the dynamic programming formulations are dependent on the accumulated wealth, $X_k$. So, the writer and the buyer's models in Proposition \ref{prop:bellman} would change as follows:\\
\textbf{Writer's model:}
\begin{subequations}
\begin{align}
&V^w_k(X_k,\omega):=\inf_{\xi_k}\rho^w_k(X_k-\xi_k \Delta S_{k+1}+ V^w_{k+1}(X_k +\xi_k \Delta S_{k+1}, \omega),\omega)\,,\,k=0,\dots,K-1\\
&V^w_K(X_K,\omega):=F(S_{\saeedmodified{K}}(\omega),Y_{\saeedmodified{K}}(\omega)) - X_K\,.
\end{align}
\end{subequations}
\textbf{Buyer's model:}
\begin{subequations}
\begin{align}
&V^b_k(X_k,\omega):=\inf_{\xi_k}\rho^b_k(X_k- \xi_k \Delta S_{k+1} + V^b_{k+1}(X_k +\xi_k \Delta S_{k+1},\omega),\omega)\,,\,k=0,\dots,K-1\\
&V^b_K(X_K,\omega):=-F(S_{\saeedmodified{K}}(\omega),Y_{\saeedmodified{K}}(\omega)) - X_K\, ,
\end{align}
\end{subequations}
Then, the equal risk price is $\price_0^* = \{ \price_0 \in \mathbb{R} | \Delta(\price_0) = \varrho^w(\price_0) - \varrho^b(\price_0) = 0 \}$ where $\varrho^w(\price_0) = V^w_0(\price_0)$ and $\varrho^b(\price_0) = V^b_0(\price_0)$. Furthermore, given that the accumulated wealth up to time $k \in \{1, \cdots , K-1 \}$ is $X_k$, the minimal risk hedging policy for both the writer and the buyer can be described respectively as:
\begin{align*}
\xi_k^{w*}(X_k,\omega)\in\arg\min_{\xi_k}\rho^w_k(X_k-\xi_k \Delta S_{k+1}+ V^w_{k+1},\omega),\,\forall\,k=1,\dots,K-1\\
\xi_k^{b*}(X_k,\omega)\in\arg\min_{\xi_k}\rho^b_k(X_k-\xi_k \Delta S_{k+1}+ V^b_{k+1},\omega),\,\forall\,k=1,\dots,K-1\,.
\end{align*}

In order to find the finite dimensional Bellman equations for option's equal risk price under non-translation invariant risk measures, we can again follow the equations in section \ref{subsec:EUoption}, without removing the accumulated portfolio wealth state variable, $X_k$, from the value functions. With assumptions similar to the translation invariant case, and further assuming that $X_k$ is measureable on $\Sigma_k$, the finite dimensional Bellman equations are then defined as follows:

\[\tilde{V}^w_K(X_K,S_K,Y_K,\theta_K^w):= F(S_K,Y_K) - X_K\,,\]
and recursively
\[\tilde{V}^w_k(X_k,S_k,Y_k,\theta_k^w):= \inf_{\xi_k} \bar{\rho}_k(X_k-\xi_k\Delta S_{k+1} + \tilde{V}(X_k + \xi_k\Delta S_{k+1},S_k+\Delta S_{k+1},Y_k+\Delta Y_{k+1},f_k(\theta_k^w)),\theta_k)\,,\]
with $\tilde{f}_k(\theta_k^w)$ be a random variable in $\mathcal{L}_p(\Omega_{k+1},\Sigma_{k+1},\mathbb{P}_{k+1})$. Then, considering that
\[V^w_K(X_K,\omega)=\tilde{V}^w_K(X_K,S_K(\omega),Y_K(\omega),\theta_K^w(\omega))\,,\]
and recursively that if $V^w_{k+1}(X_{k+1},\omega)=\tilde{V}^w_{k+1}(X_{k+1},S_{k+1}(\omega),Y_{k+1}(\omega),\theta_{k+1}^w(\omega))$ we have that:
\begin{align*}
&V^w_k(X_k,\omega)=\inf_{\xi_k}\rho^w_k(X_k,X_k-\xi_k \Delta S_{k+1}+ V^w_{k+1},\omega)\\
&=\inf_{\xi_k}\rho^w_k(X_k-\xi_k \Delta S_{k+1}+ \tilde{V}^w_{k+1}(X_{k+1},S_{k+1},Y_{k+1},\theta_{k+1}^w),\omega)\\
&=\inf_{\xi_k}\bar{\rho}^w_k(\bar{\Pi}_k(X_{k}-\xi_k \Delta S_{k+1}+ \tilde{V}^w_{k+1}(X_{k+1},S_{k+1},Y_{k+1},\theta_{k+1}^w)),\theta_k^w(\omega)),\\
&=\inf_{\xi_k}\bar{\rho}^w_k(X_k(\omega)-\xi_k \Delta S_{k+1}+ \tilde{V}^w_{k+1}(X_k(\omega),S_k(\omega)+\Delta S_{k+1},Y_k(\omega)+\Delta Y_{k+1},f(\theta_k^w)),\theta_k^w(\omega))\\
&=\tilde{V}_k(X_k(\omega),S_k(\omega),Y_k(\omega),\theta_k(\omega))\,.
\end{align*}
Now, we see that $\varrho^w(\price_0)=V^w_0(\price_0)=\tilde{V}^w_0(\price_0,S_0,Y_0,\theta_0^w)$. In the case of the buyer, similar derivations lead to the Bellman equations:
\[\tilde{V}^b_K(X_K,S_K,Y_K,\theta_K^b):= -F(S_K,Y_K) - X_K\]
and 
\[\tilde{V}^b_k(X_k,S_k,Y_k,\theta_k^b):= \inf_{\xi_k} \bar{\rho}_k^b(X_k-\xi_k\Delta S_{k+1} + \tilde{V}^b_{k+1}(X_{k+1},S_k+\Delta S_{k+1},Y_k+\Delta Y_{k+1},f_k(\theta_k^b,\omega_{k+1})),\theta_k)\,,\]
which can be used to compute $\varrho^b(\price_0)=V^b_0(\price_0)=\tilde{V}^b_0(\price_0,S_0,Y_0,\theta_0^b)$. Now, according to Remark \ref{rem:nonTI}, we can use a bisection method to find the equal risk price that guarantees $\Delta(\price_0) = 0$.
}

\subsection{Proof of Lemma \ref{thm:noComitIsMoreExpansive}}\label{app:proof:noComitIsMoreExpansive}

This can be shown by contradiction. Let us assume that $\wcstar>\wncstar$ and denote with $\tcstar$ a risk minimizing stopping time strategy for the buyer when the price of the option is set to $\wcstar$. One can straightforwardly establish that either:
\begin{align*}
\varrho^w(\wcstar,\tcstar)&\leq \sup_\tau \varrho^w(\wcstar,\tau) \leq \varrho_\tau^w(\wcstar)\leq\varrho_\tau^w(\wncstar)=\varrho_\tau^b(\wncstar)<\varrho_\tau^b(\wcstar)=\varrho^b(\wcstar,\tcstar)\,,
\end{align*}
or
\begin{align*}
\varrho^w(\wcstar,\tcstar)&\leq \sup_\tau \varrho^w(\wcstar,\tau) \leq \varrho_\tau^w(\wcstar)<\varrho_\tau^w(\wncstar)=\varrho_\tau^b(\wncstar)\leq\varrho_\tau^b(\wcstar)=\varrho^b(\wcstar,\tcstar)\,,
\end{align*}
where in the second inequality, we used the fact that the risk can only increase when the supremum over $\tau$ is evaluated after the hedging policy has been fixed. We also used in the following two strict inequalities the fact that both $\rho^w$ and $\rho^b$ are monotone and that $\wncstar$ is the \modified{unique} equal risk price without commitment, which implies that: 
\[\varrho_\tau^w(\wcstar)=\varrho_\tau^w(\wncstar)=\varrho_\tau^b(\wncstar)=\varrho_\tau^b(\wcstar)\]
is not possible. Our analysis leads to a contradiction since it implies that $\varrho^w(\wcstar,\tcstar) < \varrho^b(\wcstar,\tcstar)$ while by definition $\varrho^w(\wcstar,\tcstar) = \varrho^b(\wcstar,\tcstar)$.

\subsection{Proof of Lemma \ref{thm:midPriceWithCommitment}}\label{app:proof:midPriceWithCommitment}

This proof follows similar arguments as the proof of Proposition \ref{pr:TrInv}. In particular, one can again demonstrate that for any $\price_0\in\Re$ and any $\tau$, the minimal risk achievable are $\varrho^w(\price_0,\tau)=\varrho^w(0,\tau)+\price_0$ and $\varrho^b(\price_0,\tau)=\varrho^w(0,\tau)-\price_0$ because of the translation invariance property of $\rho^w$ and $\rho^b$. We can then prove the two conditional statements. 

First, in the case that an equal risk price $\price_0^*$ exists, based on the definition of $\price_0^*$, there must also exist a $\tau^*\in\arg\min_{\tau}\varrho^b(\price_0^*,\tau)$. This further implies that $\tau^*\in\arg\min_{\tau}\varrho^b(0,\tau)-\price_0^*$ hence that $\tau^*\in\arg\min_{\tau}\varrho^b(0,\tau)$ and can therefore play the role of $\tau_0$. Next, the definition of $\price_0^*$ also ensures that $\varrho^b(\price_0^*,\tau^*)=\varrho^w(\price_0^*,\tau^*)\in\Re$ which implies that both $\varrho^w(0,\tau^*)$ and $\varrho^b(0,\tau^*)$ are finite.

Reversely, in the case that the fair price interval is bounded and $\tau_0\in\arg\min_{\tau}\varrho^b(0,\tau)$ exists, then we can construct $\price_0^*:= (\varrho^w(0,\tau_0)-\varrho^b(0,\tau_0))/2\in\Re$. Necessarily, $\tau_0\in\arg\min_{\tau}\varrho^b(0,\tau)-\price_0^*=\arg\min_{\tau}\varrho^b(\price_0^*,\tau)$. Finally, we have that:
\[\varrho^w(\price_0^*,\tau_0)=\varrho^w(0,\tau_0)-\price_0^*=\frac{\varrho^w(0,\tau_0)-\varrho^b(0,\tau_0)}{2}=\varrho^b(0,\tau_0)+\price_0^*=\varrho^b(\price_0^*,\tau_0)\,.\]

\subsection{Proof of Proposition \ref{thm:AmWithCommitDP}}\label{app:proof:AmWithCommitDP}

In the case of the writer, the argument are exactly analogous as for the proof of Proposition \ref{prop:bellman}. In particular, one can use Proposition 3.1 and the discussion in Section 4 of \cite{pichler2018risk} to conclude that: 
\[\varrho^w(0,\tau)=\bar{V}_0^w(0,\tau)\,,\]
with 
\begin{align*}
\bar{V}^w_k(X_k,\tau,\omega)&:=\inf_{\xi_k,\{\hat{\xi}^i\}_{i=0}^k}\rho^w_k(\bar{V}^w_{k+1}(X_k+(\1\{\tau>k\}\xi_k + \sum_{\ell=0}^k\1\{\tau=\ell\}\hat{\xi}_k^\ell) \Delta S_{k+1},\tau),\omega)\\
&=\inf_{\bar{\xi}_k}\rho^w_k(\bar{V}^w_{k+1}(X_k+\bar{\xi}_k \Delta S_{k+1},\tau),\omega)\\
\bar{V}^w_K(X_K,\tau,\omega)&:=F(S_{\tau(\omega)}(\omega),Y_{\tau(\omega)}(\omega))-X_K\modified{(\omega)}\,,
\end{align*}
\modified{where the first argument of each $\bar{V}^w_k$ is a random variable in $\mathcal{L}_p(\Omega,\mathcal{F}_{t_k},\mathbb{P})$.}
By exploiting conditional translation invariance, one then easily obtains:
\[\bar{V}^w_K(X_K,\tau,\omega)=V_K^w(\tau,\omega)+\sum_{k=0}^{K-1}\1\{\tau(\omega)=k\}F(S_k(\omega),Y_k(\omega))-X_K\modified{(\omega)}\,,\]
and recursively that:
\begin{align*}
\bar{V}^w_k(X_k,\tau,\omega)&=\inf_{\bar{\xi}_k}\rho^w_k(\bar{V}^w_{k+1}(X_k+\bar{\xi}_k \Delta S_{k+1},\tau),\omega)\\
&=\inf_{\bar{\xi}_k}\rho^w_k(V_{k+1}^w(\tau)+\sum_{\ell=0}^{k}\1\{\tau=\ell\}F(S_\ell,Y_\ell)-X_k-\bar{\xi}_k \Delta S_{k+1},\omega)\\
&=\inf_{\bar{\xi}_k}\rho^w_k(V_{k+1}^w(\tau)-\bar{\xi}_k \Delta S_{k+1},\omega)+\sum_{\ell=0}^{k}\1\{\tau(\omega)=\ell\}F(S_\ell(\omega),Y_\ell(\omega))-X_k\modified{(\omega)}\\
&=V_k^w(\tau,\omega)+\sum_{\ell=0}^{k-1}\1\{\tau(\omega)=\ell\}F(S_\ell(\omega),Y_\ell(\omega))-X_k\modified{(\omega)}\,.
\end{align*}
Hence, we have that:
\[\varrho^w(0,\tau)=\bar{V}_0^w(0,\tau)=V_0^w(\tau)\,.\]

In the case of the buyer's equations, the proof is more challenging yet follows similar arguments\footnote{Note that here we diverge from the arguments used in Section 6.1.2 of \cite{pichler2018risk} to simplify exposition}. In particular, we first define an operator that computes the minimal risk under an initial price of $\price_0$, and presents an equivalent reformulation:
\begin{align*}
\varrho^b(\price_0)&:=\min_\tau \inf_{X\in\bar{\mathcal{X}}_\tau(\price_0)} \rho^b(F(S_\tau,Y_\tau)-X_K)\\
&=\inf_{Z\in\mathcal{Z},\,X\in\bar{\mathcal{X}}(\price_0,Z)} \rho^b(\sum_{k=0}^{K-1}Z_k F(S_K,Y_K)-X_K)\,,
\end{align*}
where $\mathcal{Z}:=\{Z:\Omega\rightarrow \{0,1\}^{K} \,|\, \sum_{k=0}^{K-1}Z_k \leq 1\}$ and each $Z_k$ is $\mathbb{F}_k$-adapted and captures $Z_k:=\1\{\tau=k\}$, and where
\[	\bar{\mathcal{X}}(\price_0,Z):= \left\{X:\Omega\rightarrow \Re^K \middle| \begin{array}{l}\exists X_0=\price_0,  \forall k=1,\dots,K-1,\exists \xi_k, \{\hat{\xi}_k^i\}_{i=0}^{k}\quad\\ X_{k+1}=X_{k}+(\xi_{k}+\sum_{i=0}^{k}(\hat{\xi}_{k}^i-\xi_k) Z_i) \Delta S_{k+1}\end{array} \right\}\,.\]
Once again, we use the arguments in \cite{pichler2018risk} to conclude that:  
\[\varrho^b(0)=\bar{V}_0^b(0)\,,\]
with 
\begin{align*}
\bar{V}^b_K(X_K,Z_{0:K-1},\omega)&:=-\sum_{k=0}^{K-1}Z_k\modified{(\omega)} F(S_k(\omega),Y_k(\omega))-(1-\sum_{k=0}^{K-1}Z_k\modified{(\omega)})F(S_K(\omega),Y_K(\omega))-X_K(\omega)\\
\bar{V}^b_k(X_k,Z_{0:k-1},\omega)&:=\inf_{Z_k,\xi_k,\{\hat{\xi}^i\}_{i=0}^k:Z_k\leq 1-\sum_{\ell=0}^{k-1}Z_\ell}\rho^b_k(\bar{V}^b_{k+1}(X_k+((1-\sum_{\ell=0}^kZ_\ell)\xi_k + \sum_{\ell=0}^kZ_k\hat{\xi}_k^\ell) \Delta S_{k+1},Z_{0:k}),\omega)\\
&=\inf_{Z_k,\bar{\xi}_k:Z_k\leq 1-\sum_{\ell=0}^{k-1}Z_\ell}\rho^b_k(\bar{V}^b_{k+1}(X_k+\bar{\xi}_k \Delta S_{k+1},Z_{0:k}),\omega)\,.
\end{align*}
By exploiting conditional translation invariance, one then easily obtains:
\[\bar{V}^b_K(X_K,Z_{0:K-1},\omega)=V_K^b(\sum_{k=0}^{K-1}Z_k,\omega)-\sum_{k=0}^{K-1}Z_k \modified{(\omega)}F(S_k(\omega),Y_k(\omega))-X_K(\omega)\,,\]
and recursively that:
\begin{align*}
\bar{V}^b_k(X_k,Z_{0:k-1},\omega)&=\inf_{Z_k,\bar{\xi}_k:Z_k\leq 1-\sum_{\ell=0}^{k-1}Z_\ell }\rho^b_k(\bar{V}^w_{k+1}(X_k+\bar{\xi}_k \Delta S_{k+1},Z_{0:k}),\omega)\\
&=\inf_{Z_k,\bar{\xi}_k:Z_k\leq 1-\sum_{\ell=0}^{k-1}Z_\ell }\rho^b_k(V_{k+1}^b(\sum_{\ell=0}^{k}Z_k)-\sum_{\ell=0}^{k}Z_\ell F(S_\ell,Y_\ell)-X_k-\bar{\xi}_k \Delta S_{k+1},\omega)\\
&=\inf_{Z_k,\bar{\xi}_k:Z_k\leq 1-\sum_{\ell=0}^{k-1}Z_\ell }\rho^b_k(V_{k+1}^b(\sum_{\ell=0}^{k}Z_k)-\bar{\xi}_k \Delta S_{k+1},\omega)-\sum_{\ell=0}^{k}Z_\ell\modified{(\omega)} F(S_\ell(\omega),Y_\ell(\omega))-X_k(\omega)\\
&=V_k^b(\sum_{\ell=0}^{k-1}Z_k,\omega)-\sum_{\ell=0}^{k-1}Z_\ell\modified{(\omega)} F(S_\ell(\omega),Y_\ell(\omega))-X_k\modified{(\omega)}\,.
\end{align*}
Hence, we have that $\varrho^b(0)=\bar{V}_0^b(0)$.
Optimal policies for each problem can be identified similarly as was done in the proof of Proposition \ref{prop:bellman}.

\subsection{Proof of Lemma \ref{thm:midPriceWithoutCommitment}}\label{app:proof:midPriceWithoutCommitment}

This proof mainly relies on the translation invariance property together with the property that once again $\bar{\mathcal{X}}_\tau(\price_0)=\price_0+\bar{\mathcal{X}}_\tau(0)$. These two properties can be used to show that:
	\begin{align}
	\varrho_\tau^w(\price_0)&=\inf_{X\in\bar{\mathcal{X}}_\tau(0)+\price_0}\sup_\tau \rho^{w}(F(S_\tau,Y_\tau)-X_K\modified{(\tau)})\notag\\
	&=\inf_{X\in\bar{\mathcal{X}}_\tau(0)}\sup_\tau \rho^{w}(F(S_\tau,Y_\tau)-X_K\modified{(\tau)}-\price_0)\notag\\
	&=\inf_{X\in\bar{\mathcal{X}}_\tau(0)}\sup_\tau \rho^{w}(F(S_\tau,Y_\tau)-X_K\modified{(\tau)})-\price_0\notag\\
	&=\inf \{s | \inf_{X\in\bar{\mathcal{X}}_\tau(0)} \sup_\tau \rho^{w}(F(S_\tau,Y_\tau)-X_K\modified{(\tau)}) \le s\}  - \price_0\notag\\
	&=\inf \{s | \varrho_\tau^w(s) \le 0 \}-\price_0 = 	\price_0^w - \price_0\,.
	\end{align}
	Similarly for the buyer, we have that $\varrho_\tau^b(\price_0)=-\price_0^b + \price_0$. The rest follows as in the proof of Proposition \ref{pr:TrInv}.
	
\removed{	Hence, we can obtain our result by verifying both directions of the biconditional logical connective. First, given that an equal risk price exists, say $\price_0^*\in\Re$, it must be that both $\varrho_\tau^b(\price_0^*)$ and $\varrho_\tau^b(\price_0^*)$ are members of $\Re$. This necessarily implies that $\price_0^w\in\Re$ and $\price_0^b\in\Re$ thus that the fair price interval is bounded. Conversely, if the fair price interval is bounded, then one can verify that the midpoint $\price_0^*:=(\price_0^w+\price_0^b)/2$ does satisfy the equal risk price condition:
	\[\varrho^w(\price_0^*)=\price_0^w - \price_0^* =\price_0^w/2 - \price_0^b/2 = -w^b_0 + \price_0^* =  \varrho^b(\price_0^*)\,.\]
	Furthermore, this midpoint can be calculated as
	\begin{align*}
	\price_0^*&=(1/2)(\inf \{\price_0 | \varrho^w(\price_0) \le 0 \}-\sup \{-s | \varrho^b(-s) \le 0 \})\\
	&=(1/2)(\price_0^w+\price_0^b)=(1/2)(\varrho^w(0) - \varrho^b(0)))\,.
	\end{align*}
		This completes our proof. }

\subsection{Proof of Proposition \ref{thm:DPwithoutCommitment}}\label{app:proof:DPwithoutCommitment}

This proof focuses on the case of the writer given that the buyer's problem was already studied in the proof of Proposition \ref{thm:AmWithCommitDP}. In particular, the proof follows similar lines as in Section 6.1.1 of \cite{pichler2018risk} and in fact extends that result to a case where hedging is allowed to pursue past the exercise time up to the end of the horizon. We start by considering that the self financing hedging policy described by $\xi$ and $\hat{\xi}$ is fixed and reformulate the worst-case exercise time problem. We then look into reformulating the optimization of the hedging policy.\\

\noindent\textbf{Step 1 (Worst-case exercise time problem):} For any fixed hedging strategy, one can define the writer's worst-case exercise time problem as: 
\[\nu_0 := \sup_\tau \rho^{w}\left(F(S_{\tau},Y_{\tau})-\sum_{\ell=0}^{\tau-1}\xi_\ell\Delta S_{\ell+1} - \sum_{\ell=\tau}^{K-1}\hat{\xi}_\ell^\tau\Delta S_{\ell+1}\right)\,.\]
In order to find a dynamic programming formulation for this problem, we start by reformulating it in the form of an \quoteIt{optimal stopping problem} as defined in \cite{pichler2018risk}. In particular, we can consider that:
\[\nu_0=\sup_\tau \rho^w(E_\tau^1)\,,\]
where $E_k^1(\omega):=\rho^w_{k,K}(F_k(S_k,Y_k)-\sum_{\ell=0}^{k-1}\xi_\ell\Delta S_{\ell+1} - \sum_{\ell=k}^{K-1}\hat{\xi}_\ell^k\Delta S_{\ell+1},\omega)$
with $\rho^w_{k,K}(X):=\rho^w_k(\rho^w_{k+1}(\cdots\rho_{K-1}(X)\cdots)$. Hence, based on Theorem 6.4 in \cite{pichler2018risk}, we can conclude that if we define:
\begin{align*}
E_K^2(\omega)&:=E_K^1(\omega)\\
E_k^2(\omega)&:=\max(E_k^1(\omega),\,\rho_k(E_{k+1}^2,\omega))\,,
\end{align*}
and
\[\tau_m^*(\omega):=\min\{ k\,|\,E_k^2(\omega)=E_k^1(\omega),\,m\leq t \leq T \}\,,\]
then $\tau_m^*$ is an optimal solution to: 
\[\sup_{\tau:\tau\geq m} \rho^{w}(F(S_{\tau},Y_{\tau})-\sum_{\ell=0}^{\tau-1}\xi_\ell\Delta S_{\ell+1} - \sum_{\ell=\tau}^{K-1}\hat{\xi}_\ell^\tau\Delta S_{\ell+1})\,,\]
and $\nu_0=E_0^2$.\\

\noindent\textbf{Step 2 (Optimal hedging optimization problem):} Based on our analysis of the worst-case exercise time problem, we have found that the optimal hedging problem has the following form:
\begin{align*}
\varrho_\tau^w(0)&=\inf_{\xi,\hat{\xi}} E_0^2(\xi,\hat{\xi})\\
&=\inf_{\xi,\hat{\xi}}\mathcal{R}_{0,1}(E_0^1(\hat{\xi}^0),\mathcal{R}_{1,2}(E_1^1(\xi_0,\hat{\xi}^1),\cdots,\mathcal{R}_{K-1,K}(E_{K-1}^1(\xi_{0:K-1},\hat{\xi}^{K-1}),E_{K}^1(\xi))\,,
\end{align*}
where $\mathcal{R}_{k,k+1}(X,Y,\omega):=\max(X(\omega),\rho_k^w(Y,\omega))$ and where we made explicit the influence of $\xi$ and $\hat{\xi}$ on each $E_k^1$. As argued in \cite{pichler2018risk}, given that each $\mathcal{R}_{t,t+1}(\cdot,\cdot)$ is monotone, one can apply the interchangeability principle to generate the following reformulation:
\begin{align*}
\varrho_\tau^w(0)&=\mathcal{R}_{0,1}(\inf_{\hat{\xi}^0}E_0^1(\hat{\xi}^0),\inf_{\xi_0}\mathcal{R}_{1,2}(\inf_{\hat{\xi}^1}E_1^1(\xi_0,\hat{\xi}^1),\cdots,\inf_{\xi_{K-1}}\mathcal{R}_{K-1,K}(\inf_{\hat{\xi}^{K-1}}E_{K-1}^1(\xi_{0:K-1},\hat{\xi}^{K-1}),E_{K}^1(\xi))\,.
\end{align*}
Based on this argument, we create the following operators:
\begin{align*}
\bar{V}^w_k(1,X_k,\omega)&:=\inf_{\hat{\xi}_{k:K-1}^k}\rho^w_{k,K}(F_k(S_k,Y_k)-X_k -\sum_{\ell=k}^{K-1} \hat{\xi}_\ell^k \Delta S_{\ell+1},\omega)\\
\bar{V}^w_k(0,X_k,\omega)&:=\max(\bar{V}_k^w(1,X_k,\omega),\;\inf_{\xi_k}\rho^w_k(\bar{V}^w_{k+1}(0,X_k+\xi_k \Delta S_{k+1}),\omega))\\
\bar{V}^w_K(0,X_K,\omega)&:=F_K(S_K(\omega),Y_K(\omega))-X_K\modified(\omega)\,,
\end{align*}
in order to have that $\varrho_\tau^w(0)=\bar{V}^w_0(0,0)$\modified{, and where we again let the second argument of $\bar{V}^w_k$ be a random variable in $\mathcal{L}_p(\Omega,\mathcal{F}_{t_k},\mathbb{P})$}.

In the case of $\bar{V}^w_k(1,X_k,\omega)$, we can further apply the interchangeability principle to get that for all $k=0,\dots,K-1$,
\begin{align*}
\bar{V}^w_k(1,X_k,\omega)&=\inf_{\hat{\xi}_k^k}\rho^w_k(\inf_{\hat{\xi}_{k+1}^k}\rho^w_{k+1}(\cdots \inf_{\hat{\xi}_{K-1}^k}\rho^w_{K-1}(F_k(S_k,Y_k)-X_k -\sum_{\ell=k}^{K-1} \hat{\xi}_\ell^k \Delta S_{\ell+1})\cdots),\omega)\\
&=F_k(S_k(\omega),Y_k(\omega))-X_k\modified{(\omega)}+\\
&\quad\quad\quad\inf_{\hat{\xi}_k^k} \rho^w_k(- \hat{\xi}_k^k \Delta S_{k+1} + \inf_{\hat{\xi}_{k+1}^k}\rho^w_{k+1}(- \hat{\xi}_{k+1}^k \Delta S_{k+2}+\cdots \inf_{\hat{\xi}_{K-1}^k}\rho^w_{K-1}(-\hat{\xi}_{K-1}^k \Delta S_{K})),\omega)\\
&=F_k(S_k(\omega),Y_k(\omega))-X_k\modified{(\omega)}+V_k^w(1,\omega)\,,
\end{align*}
where we applied conditional translation invariance. While one can verify that we also have:
\begin{align*}
\bar{V}^w_K(0,X_K,\omega)&= V^w_K(0,\omega)-X_K\modified{(\omega)}\\
\bar{V}^w_k(0,X_k,\omega)&= \max(V_k^w(1,\omega)+F_k(S_k(\omega),Y_k(\omega))-X_k\modified{(\omega)},\;\inf_{\xi_k}\rho^w_k(V^w_{k+1}(0)-X_k-\xi_k \Delta S_{k+1},\omega))\\
&=\max(V_k^w(1,\omega)+F_k(S_k(\omega),Y_k(\omega)),\;\inf_{\xi_k}\rho^w_k(V^w_{k+1}(0)-\xi_k \Delta S_{k+1},\omega))-X_k\modified{(\omega)}\\
&=V^w_k(0,\omega)-X_k\modified{(\omega)}\,.
\end{align*}
Hence, $\varrho_\tau^w(0)=\bar{V}^w_0(0,0)=V^w_0(0)$. An optimal hedging policy can be identified with an optimal solution to the infimum operations in equation \eqref{eq:erpnoComV1} or \eqref{eq:erpnoComV2} depending on whether the option was exercised at a period smaller or equal to $k$.

\subsection{Proof of Corollary \ref{thm:AmStopHedge}}\label{app:proof:AmStopHedge}

We start by looking at the case of an American option with commitment. Based on Proposition \ref{thm:AmWithCommitDP}, we can first prove that, given any exercise policy $\tau$, the writer should stop hedging after exercise by studying, for each $k$, whether $\xi_k^{i*}(\omega)=0$ is optimal when $\tau(\omega)\leq k$. Specifically, if $\tau(\omega)\leq k$, then we have that: 
\begin{align*}
\arg\min_{\xi_k}\rho^w_k( V^w_{k+1}(\tau)-&\xi_k \Delta S_{k+1},\omega)\\
&=\arg\min_{\xi_k}\inf_{\xi_{k+1:K-1}}\rho^w_{k,K}(\sum_{\ell=k+1}^K\1\{\tau=\ell\}F(S_{\ell},Y_\ell) -\sum_{\ell=k}^{K-1}\xi_\ell \Delta S_{\ell+1},\omega)\\
&=\arg\min_{\xi_k}\inf_{\xi_{k+1:K-1}}\rho^w_{k,K}(-\sum_{\ell=k}^{K-1}\xi_\ell \Delta S_{\ell+1},\omega)\;\supseteq\; \{ 0\}\,,
\end{align*}
since the $\rho_k^w$ satisfies the bounded conditional market risk assumption and is conditionally coherent, thus $\inf_{\xi_{k:K-1}}\rho^w_{k,K}(-\sum_{\ell=k}^{K-1}\xi_\ell \Delta S_{\ell+1},\omega)=0$. This confirms that if $\tau(\omega)\leq k$, then $\hat{\xi}_k^{i*}(\tau,\omega):=0$ is optimal for all $i\leq k$ so that the number of shares of the risky asset becomes:
\[\xi_{k}\1\{\tau>k\}+\sum_{i=0}^{k}\hat{\xi}_{k}^i\1\{\tau=i\}=\xi_{k}\cdot 0+ \sum_{i=0}^{k}0\cdot\1\{\tau=i\}=0\,.\]
In other words, it is optimal to stop hedging at $\tau(\omega)$.

A similar argument can be used for the buyer. Namely, we can study, for each $k$, the structure of $\xi_k^{i*}(\omega)$ when $\tau_0(\omega)\leq k$. This is done as follows:
\begin{align*}
\arg\min_{\xi_k}\rho^b_k(-&\xi_k \Delta S_{k+1}+ V^b_{k+1}(\1\{\tau_0\leq k\}),\omega)=\arg\min_{\xi_k}\inf_{\xi_{k+1:K-1}}\rho^b_{k:K}(-\sum_{\ell=k}^{K-1}\xi_\ell \Delta S_{\ell+1},\omega)\;\supseteq\; \{ 0\}\,,
\end{align*}
which again implies that it is optimal to stop hedging at $\tau_0$.

For the case of an American option without commitment, the same argument has for the case with commitment applies for the buyer. On the other hand, for the writer we can retrieve the optimal hedging policy from Proposition \ref{thm:DPwithoutCommitment}. Looking carefully, for each $k$, at the structure of $\xi_k^{i*}(\omega)$ when $\tau(\omega)\leq k$, we realize that the same arguments apply
\begin{align*}
\arg\min_{\xi_k}\rho^w_k(-\xi_k \Delta S_{k+1}+ V^w_{k+1}(\1\{\tau\leq k\}),\omega)
&=\arg\min_{\xi_k}\inf_{\xi_{k+1:K-1}}\rho^w_{k,K}\left(-\sum_{\ell=k}^{K-1}\xi_\ell \Delta S_{\ell+1},\omega\right)\;\supseteq\; \{ 0\}\,.
\end{align*}
Hence, once again it is optimal to stop hedging starting at $\tau$.

\section{Appendix for Section \ref{sec:numInvestigate}}

\subsection{Verifying the Bounded (Conditional) Market Risk Property for Worst-case Risk Measures}
In this section we identify sufficient conditions under which the one-step decomposition of a worst-case risk measure satisfies the bounded conditional market risk property. However, before studying such conditions we need to first define a useful projection operator.

\begin{definition}\label{def:USTimeProjection}
	Given an uncertainty set $ \U \subseteq \mathbb{R}^K $, and a history of observations $ \hat{r}_{1:k-1}\in \mathbb{R}^{k-1} $, we define the operation of projecting $ \U $ over the time interval $\{k,\dots,k'\}$ with $ 1\le k\le k'\le K $ as follows:
	\begin{equation}\nonumber
	\begin{aligned}
	&\U_{k:k'}({\hat{r}_{1:k-1}}):=\left\{r \in \mathbb{R}^{k'-k+1} \left| 
	\begin{aligned}
	&\text{If } k'<K, \quad\exists \bar{r}\in\Re^{K-k'}, [\hat{r}_{1:k-1}^T\;\;r^T\;\;\bar{r}^T]^T\in \U\\
	&\text{If } k'=K, \quad [\hat{r}_{1:k-1}^T \;\;r^T]^T\in \U\\
	\end{aligned}
	\right. \right\}\,.
	\end{aligned}
	\end{equation}
\end{definition}

This definition is helpful in describing, for a given worst-case risk measure that exploits some uncertainty set $\U$, the set of all realisations of the return vector for which the bounded conditional market risk property is satisfied.
\begin{definition}\label{def:AUset}
	Given a worst-case risk measure using an uncertainty set $ \U \subseteq \mathbb{R}^K $, we define the set of returns with bounded conditional market risk as follows:
	\begin{equation}\label{eq:necessaryCondition1}
	\mathcal{A}(\U) := \left\{r \in \mathbb{R}^K \middle| \forall k \in \{0,...,K-1\},\, \inf_{\myzeta_k,\dots,\myzeta_{K-1}}\rho_{k,K}(-\sum_{\ell=k}^K\myzeta_\ell r_{\ell+1},r) \in\;]-\infty,\;0]\right\}\,.
	\end{equation}
\end{definition}
In particular, one can also reformulate the definition of $\mathcal{A}(\U)$ as follows:
\[	\mathcal{A}(\U) = \left\{r \in \mathbb{R}^K \middle| \begin{array}{c}\forall k \in \{0,...,K-1\},\,\\ \U_{k+1:K}(r_{1:k})=\emptyset \;\vee\; \inf_{\myzeta_k,\dots,\myzeta_{K-1}} \sup_{\bar{r}_{k+1:K} \in \U_{k+1:K}(r_{1:k})} -\sum_{\ell=k}^{K-1}\myzeta_\ell \bar{r}_{\ell+1} \in\;]-\infty,\;0]\end{array}\right\}\,,\]
where $\bar{r}_{k+1:K}$ refers to a vector in $\Re^{K-k-1}$ with indexes in the range $\{k+1,\dots,K\}$, and where we use $\myzeta_\ell$ as shorthand notation for $\myzeta_\ell(\bar{r}_{k:\ell})$. We will repeat this abuse of notation throughout the section to simplify the presentation of equations.

Based on \modified{Definition} \ref{def:AUset}, it is clear that a worst-case risk measure will satisfy the bounded conditional market risk condition if $\U\subseteq \mathcal{A}(\U)$. This is formally stated by the following lemma.
\begin{lemma}\label{thm:BCMRwithAu}
	Let $\rho$ be a worst-case risk measure that uses an uncertainty set $ \U \subseteq \mathbb{R}^K$ such that $\U\subseteq \mathcal{A}(\U)$, then $\rho$ necessarily satisfies the bounded conditional market risk property.
\end{lemma}
\begin{proof}
This result simply follows from the fact that for any $r\in\Re^K$ and any $k\in\{1,\dots,K\}$, two situation can occur. First, the set $\U_{k+1:K}(r_{1:k})$ might be empty, which leads to $\inf_{\myzeta_k,\dots,\myzeta_{K-1}}\rho_{k,K}(-\sum_{\ell=k}^K\myzeta_\ell r_{\ell+1},r)=0$ by the definition of $\rho_k(X,r)$ thus the market risk is bounded for this realization. Secondly, one should investigate the case where $\U_{k+1:K}(r)$ is non-empty. In this case, there exists a $\hat{r}\in\U\subseteq \mathcal{A}(\U)$ such that $r_{1:k}=\hat{r}_{1:k}$. Hence, one can verify that:
\begin{align*}
\inf_{\myzeta_k,\dots,\myzeta_{K-1}}\rho_{k,K}(-\sum_{\ell=k}^K\myzeta_\ell r_{\ell+1},r)&=\inf_{\myzeta_k,...,\myzeta_{K-1}} \sup_{\bar{r}_{k+1:K} \in \U_{k+1:K}(r_{1:k})} -\sum_{\ell=k}^{K-1}\myzeta_\ell \bar{r}_{\ell+1}\\
&=\inf_{\myzeta_k,...,\myzeta_{K-1}}  \sup_{\bar{r}_{k+1:K} \in \U_{k+1:K}(\hat{r}_{1:k})} -\sum_{\ell=k}^{K-1}\myzeta_k \bar{r}_{k+1}\in\;]-\infty,\,0]\,,
\end{align*}
based on the fact that $\hat{r}\in\mathcal{A}(\U)$. This implies that the conditional market risk is bounded on all of $\Re^K$ which is a stronger condition than in Assumption \ref{ass:BCMR} where the condition is only imposed with probability one.
\end{proof}

Based on the above discussion, given an arbitrary uncertainty set which might not satisfy the condition $\U\subseteq \mathcal{A}(\U)$, it therefore appears that we are in need of a procedure that would select a subset $\U'$ of $\U$ for which this property is satisfied. One attractive candidate takes the form of the following set which we will call the no-arbitrage subset of $\U$, when it exists.

\begin{definition}\label{def:ArbitrageFreeSupportSet}
	Given an uncertainty set $\U$, we define the no-arbitrage subset $\U^{na}$ of $\U$ as the largest set $\U'\subseteq \U$ that satisfies $\U'\subseteq \mathcal{A}(\U')$. Mathematically, $\U^{na}$ satisfies the following two properties:
	\begin{enumerate}
		\item $\U^{na} \subseteq  \mathcal{A}(\U^{na})$
		\item $ \forall \U' \subseteq \U, \U'\subseteq \mathcal{A}(\U')$ we have that $ \U' \subseteq \U^{na} $.
	\end{enumerate}
\end{definition}

Considering the previous definitions, one might wonder if such a no-arbitrage subset always exists. The following theorem confirms that it does always exist when $\U$ is both closed and  convex.


\begin{theorem}\label{thm:ArbitrageFreeSupportSet}
	Given that $\U$ is convex and closed, the no-arbitrage subset of $ \U $ is equal to $ \mathcal{V}(\U)\cap \U $ where
	\[	\mathcal{V}(\U)=\left\{r \in \mathbb{R} ^K\middle| 0\in \U, \;\left(\sum_{j=1}^{k}e_j r_j \in \U\right) \vee \left(r_{1:k}\notin\U_{1:k}\right), \forall k=1,...,K-1 \right\}\,,
\]
	with $ e_j\in \mathbb{R}^K $ as the $ j $-th column of identity matrix, and where $\mathcal{V}(\U)=\emptyset$ if $0\notin\U$.
\end{theorem}

\begin{proof}
	The proof of this theorem is divided in four parts. First, we show that $ \mathcal{V}(\U)=\mathcal{A}(\U) $.  This step is itself divided in two parts, namely first that $ \mathcal{V}(\U)\subseteq \mathcal{A}(\U) $ and then that $ \mathcal{V}(\U) \supseteq \mathcal{A}(\U) $. The second step consists in proving that $\mathcal{V}(\U)\cap \U $ satisfies the two conditions of the no-arbitrage subset $\U^{na}$. 
	
\noindent\textbf{Step 1.a ($ \mathcal{V}(\U) \subseteq \mathcal{A}(\U) $).} 
Given any member $r$ of $\mathcal{V}(\U)$, we know that for all $k=1,\dots,K$, either $r_{1:k}\notin \U_{1:k}$ which leads to: 
\[\inf_{\myzeta_k,\dots,\myzeta_{K-1}}\rho_{k,K}(-\sum_{\ell=k}^K\myzeta_\ell r_{\ell+1},r)=0\,,\]
by definition. Otherwise, the vector $[r_{1:k-1}^T\;\; 0_{1:K-k+1}^T]^T\in\U$ thus we can conclude that:
	\begin{equation}\nonumber
	\begin{aligned}
	&\inf_{{\myzeta}_k,...,\myzeta_{K-1}}\sup_{\bar{r}_{k+1:K} \in \U_{k+1:K}(r_{1:k})} -\sum_{\ell=k}^{K-1}\myzeta_\ell r_{\ell+1} \ge \inf_{{\myzeta}_k,...,\myzeta_{K-1}} -\sum_{\ell=k}^{K-1}\myzeta_\ell \cdot 0 =0 > - \infty\,.
	\end{aligned}
	\end{equation}
	From this, we conclude that $ \mathcal{V}(\U) \subseteq \mathcal{A}(\U) $.\\
	
\noindent\textbf{Step 1.b ($ \mathcal{A}(\U) \subseteq \mathcal{V}(\U)$).} Given any member $r$ of $\mathcal{A}(\U)$, for any $k=1,\dots,K-1$, we have that: 
\[\inf_{\myzeta_k,\dots,\myzeta_{K-1}}\rho_{k,K}(-\sum_{\ell=k}^K\myzeta_\ell r_{\ell+1},r)> -\infty\,.\]
This means that either $r_{1:k}\notin\U_{1:k}$ or
\[\inf_{\myzeta_k,\dots,\myzeta_{K-1}} \sup_{\bar{r}_{k+1:K} \in \U_{k+1:K}(r_{1:k})} -\sum_{\ell=k}^{K-1}\myzeta_\ell \bar{r}_{\ell+1} > -\infty\,.
\]
We can further process this second condition by rewriting it as:
\[\inf_{\myzeta_k}\sup_{\bar{r}_{k+1} \in \U_{k+1}(r_{1:k})} \myzeta_k \bar{r}_{k+1} + \pi_{k+1}([r_{1:k}^T\;\; \bar{r}_{k+1}]^T) > -\infty\,,\]
where $\U_{k+1}(r_{1:k})$ is short for $\U_{k+1:k+1}(r_{1:k})$, and with
\[\pi_{k}(r_{1:k}):=\inf_{\myzeta_{k}}\sup_{\bar{r}_{k+1} \in \U_{k+1}(r_{1:k})} \myzeta_k \bar{r}_{k+1} + \pi_{k+1}([r_{1:k}^T \bar{r}_{k+1}]^T)\,,\]
for all $k=0,\dots, K-1$ while $\pi_{K}(r):=0$. Yet, one quickly realizes that:
\begin{align*}
\pi_{K-1}(r_{1:K-1})&=\inf_{\myzeta_{K-1}}\sup_{\bar{r}_{K} \in \U_{K}(r_{1:K-1})}\myzeta_{K-1}\bar{r}_{K}\\
&=\inf_{\myzeta_{K-1}}\max\left(\myzeta_{K-1}\inf_{\bar{r}_{K} \in \U_{K}(r_{1:K-1})}\bar{r}_{K},\;\myzeta_{K-1}\sup_{\bar{r}_{K} \in \U_{K}(r_{1:K-1})}\bar{r}_{K}\right)\\
&=\left\{\begin{array}{cl}0 & \mbox{if $0\in \U_{K}(r_{1:K-1})$}\\-\infty & \mbox{otherwise}\end{array}\right.\,,
\end{align*}
where the second equality follows from the fact that $\U_{K}(r_{1:K-1})$ is a closed interval given that $\U$ is convex and closed. The third equality comes from the fact if $0\in\U_{K}(r_{1:K-1})$ then the infimum over $\myzeta_k$ is reached by $\myzeta_k=0$, while when it $\U_{K}(r_{1:K-1})$ does not include zero, then the infimum can reached an arbitrarily low value since the sign of $\bar{r}_K$ is determined. Consequently, by induction, for any $k=0,\dots,K-1$, it must actually be that:
\begin{align*}
\pi_{k}(r_{1:k})&=\inf_{\myzeta_{k}}\sup_{\bar{r}_{k+1} \in \U_{k+1}(r_{1:k})} \myzeta_{k}\bar{r}_{k+1} + \pi_{k+1}([r_{1:k}^T\;\; \bar{r}_{k+1}]^T)\\
&= \inf_{\myzeta_{k}}\max\left(\myzeta_{k}\inf_{\bar{r}_{k+1}:[\bar{r}_{k+1} \;\;0_{k+2:K}^T]^T \in \U_{k+1:K}(r_{1:k})}\bar{r}_{k+1},\;\myzeta_{k}\sup_{\bar{r}_{k+1}:[\bar{r}_{k+1} \;\;0_{k+2:K}^T]^T \in \U_{k+1:K}(r_{1:k})}\bar{r}_{k+1}\right)\\
&=\left\{\begin{array}{cl}0 & \mbox{if $0\in \U_{k+1:K}(r_{1:k})$}\\-\infty & \mbox{otherwise}\end{array}\right.\,.
\end{align*}

Based of this argument, we must therefore conclude that if $r \in \mathcal{A}(\U)$, then for all $k$, either $r_{1:k}\notin \U_{1:k}$ or $\pi_{k}(r_{1:k})>-\infty$ hence that $0\in \U_{k+1:K}(r_{1:k})$. Overall, this confirms that $r\in\mathcal{V}(\U)$.

\noindent\textbf{Step 2.a ($\U^{na}\subseteq \mathcal{A}(\U^{na})$).}	To prove this property, we need to show that:
\[\U\cap \mathcal{A}(\U) \subseteq \mathcal{A}(\U\cap \mathcal{A}(\U))\,,\]
which is equivalent to showing that:
\[\U\cap \mathcal{V}(\U) \subseteq \mathcal{V}(\U\cap \mathcal{V}(\U))\,,\]
since $\U$ is convex and closed so that $\mathcal{A}(\U)=\mathcal{V}(\U)$ and therefore $\U\cap \mathcal{A}(\U)=\U\cap \mathcal{V}(\U)$ is also convex and closed so that $\mathcal{A}(\U\cap \mathcal{A}(\U))=\mathcal{V}(\U\cap \mathcal{A}(\U))=\mathcal{V}(\U\cap \mathcal{V}(\U))$. We will tackle the second equivalent condition, where we will make use the following representation:
\[	\mathcal{V}(\mathcal{V}(\U)\cap \U) =\left\{r \in \Re^K \middle| \sum_{j=1}^{k}e_j r_j \in \mathcal{V}(\U)\cap \U\;\vee \;\left(r_{1:k}\notin(\mathcal{V}(\U)\cap \U)_{1:k}\right), \forall k=1,...,K-1 \right\}\,.
\]
Specifically, given any $r\in\U\cap \mathcal{V}(\U)\subseteq \Re^K$, and for all $k=1,\dots,K$, we will confirm that $\sum_{j=1}^{k}e_j r_j \in \mathcal{V}(\U)\cap \U$. We can first check that:
\[\sum_{j=1}^k e_j r_j \in \U\,,\]
since $r\in\mathcal{V}(\U)$. Furthermore, letting $w:=\sum_{j=1}^k e_j r_j $, we can further check that for all $\ell=1,\dots,K$,
\[\sum_{i=1}^\ell e_i w_i = \sum_{i=1}^\ell e_i e_i^T \left(\sum_{j=1}^k e_j r_j\right)=\sum_{j=1}^k \sum_{i=1}^\ell e_i e_i^T e_j r_j= \sum_{j=1}^{\min(k,\ell)} e_j r_j \in \U\,,\]
since again $r\in\mathcal{V}(\U)$, which implies that $w\in\U\cap \mathcal{V}(\U)$. Based on these arguments, we can conclude that $r\in\mathcal{V}(\U\cap \mathcal{V}(\U))$.

\noindent\textbf{Step 2.b ($\U^{na}$ is the largest).}	The second property is proved as follows:
	\begin{equation}\nonumber
	\begin{aligned}
	\U' \subseteq \U &\Rightarrow \U'\cap\mathcal{V}(\U') \subseteq \U\cap\mathcal{V}(\U)\Rightarrow \U'=\U'\cap\mathcal{A}(\U')=\U'\cap\mathcal{V}(\U') \subseteq \U\cap \mathcal{V}(\U)\,,
	\end{aligned}
	\end{equation}
	where the first implication comes from the definition of $\mathcal{V}(\mathcal{U'})$ and the fact that $ \mathcal{U'} \subseteq \U$. The second implication first  exploits the fact that $\U'\subseteq \mathcal{A}(\U')$ and then that $\mathcal{A}(\U')=\mathcal{V}(\U')$. This concludes the proof.
	
\end{proof}

\subsubsection{Bounded Conditional Market Risk Property for $\U_1$}\label{proveConditionalBoundedMktRiskBandi}

Exploiting the result of Theorem \ref{thm:ArbitrageFreeSupportSet}, we can now provide a proof of Lemma \ref{lemu1}. Specifically, since $\U_1$ is a closed convex set, the theorem provides a recipe to construct the no-arbitrage subset of $\U_1$, i.e. $\U^{na}:=\U\cap \mathcal{V}(\U)$. In the context that is studied $\mathcal{V}(\U)$ reduces to
\[	\mathcal{V}(\U)=\left \{r \in \mathbb{R}^K \middle| \begin{array}{c}\max_{k\in\{1,\dots,K\}} \mu \sqrt{kT/K}/\sigma  -\Gamma\le 0\\\left(\max_{k'\in\{k,\dots,K\}} \left| \frac{\sum_{\ell=1}^{k}\log (1+r_\ell)-\mu k'T/K}{\sigma\sqrt{k'T/K}} \right| -\Gamma\le 0\right)\vee(r_{1:k}\notin \U_{1:k}), \forall k \in \{1,...,K\}\end{array}\,. \right \}
\]
One can easily verify that $\mathcal{W}\cap\U=\mathcal{V}(\U)\cap\U$ under the conditions of Lemma \ref{lemu1}.

\subsubsection{Bounded Conditional Market Risk Property for $\U_2$}\label{proveConditionalBoundedMktRisk}

In this section, we show that the worst-case risk measure defined based on the set ${\U}_2$ satisfies the bounded conditional market risk property by showing that $\mathcal{A}(\U_2)=\Re^K$ and exploiting Lemma \ref{thm:BCMRwithAu}. Specifically, we can show that for all $r\in\Re^K$ and all $k\in\{0,\dots,K-1\}$, either $r_{1:k}\notin\U_{1:k}$ or
\[
\pi_{k}(r_{1:k}):=\inf_{\myzeta_{k},...,\myzeta_{K-1}}\sup_{\bar{r}_{k+1:K}\in{\U}_{k+1:K}(r_{1:k})}-\sum_{\ell=k}^{K-1}\myzeta_{\ell}\bar{r}_{\ell+1}=0\,.
\]
Indeed, when $r_{1:k}\notin\U_{1:k}$, one can rewrite
\[\pi_{k}(r_{1:k})=\inf_{\myzeta_{k}}\sup_{\bar{r}_{k+1}\in{\U}_{k+1}(r_{1:k})}-\myzeta_{k}\bar{r}_{k+1}+\pi_{k+1}([r_{1:k}^T\;\;\bar{r}_{k+1}]^T)\,,\]
where 
\[\pi_{K-1}(r_{1:K-1})=\inf_{\myzeta_{K-1}}\sup_{\bar{r}_{K}\in{\U}_{K:K}(r_{1:K-1})}-\myzeta_{K-1}\bar{r}_{K}\,.\]
Given that for all $k$, the function $\pi_k$ is evaluated with some non-empty and symmetric ${\U}_{k+1}(r_{1:k})$, we thus have that:
\[\pi_{K-1}(r_{1:K-1})=\inf_{\myzeta_{K-1}}|\myzeta_{K-1}|\sup_{\bar{r}_{K}\in{\U}_{K:K}(r_{1:K-1})} \bar{r}_{K}=0\,,\]
and recursively for $k=K-1,\dots,0$,
\begin{align*}
\pi_{k}(r_{1:k})&=\inf_{\myzeta_{k}}\sup_{\bar{r}_{k+1}\in{\U}_{k+1}(r_{1:k})}-\myzeta_{k}\bar{r}_{k+1}+\pi_{k+1}([r_{1:k}^T\;\;\bar{r}_{k+1}]^T)=\inf_{\myzeta_{k}}\sup_{\bar{r}_{k+1}\in{\U}_{k+1}(r_{1:k})}-\myzeta_{k}\bar{r}_{k+1}+0\\
&=\inf_{\myzeta_{k}} |\myzeta_{k}|\sup_{\bar{r}_{k+1}\in{\U}_{k+1}(r_{1:k})}\bar{r}_{k+1}=0\,.
\end{align*}
This confirms that the worst-case risk measure defined based on the set ${\U}_{2}$
satisfies bounded conditional market risk property.

\subsection{Worst-case Risk Measures with $\U_1$ or $\U_2$ Satisfying the Markov Property}\label{proveMarkov}

In this section we identify two state processes $\theta_k:\Re^K\rightarrow \Re$ under which the worst-case risk measures with $\U_1$ and $\U_2$ are respectively Markovian.

Starting with the set inspired from \cite{bandi2014robust}, we let $\theta_k:=\sum_{\ell=1}^k \log(1+r_{\ell})$. With this definition in hand, we can demonstrate the properties that are described in Definition \ref{def:MarkovianRiskMeasure}. First, we have that $\theta_{k+1} = \sum_{\ell=1}^{k+1} \log(1+r_{\ell})= \log(1+r_{k+1})+\theta_k$ and hence can be measured directly from $(\theta_k,r_k+1)$. Second, we can confirm that for all $X\in\mathcal{L}_p(\Omega,\mathcal{F}_{k+1},\mathbb{P})$, if $r_{1:k}\in \U_{1:k}$, then: 
\begin{align*}
\rho_k(X,r)&=\sup_{r'\in \U:r'_{1:k}=r_{1:k}} X(r')\\
&=\sup_{\bar{r}_{k+1}\in \U_k^\theta(\sum_{\ell=1}^k\log(1+r_\ell))} X([r_{1:k}^T\;\; \bar{r}_{k+1}\;\; r_{k+2:K}]^T)\\
&=\sup_{\bar{r}_{k+1}\in \U_k^\theta(\theta_k(r))} \Pi_k(X,\bar{r}_{k+1}) = \bar{\rho}_k(\Pi_k(X,r),\theta_k(r))\,,
\end{align*}
where
\[\U_k^\theta(\theta_k):=\left\{ r\in\Re\,\middle|\, \left|\frac{\theta_k+\log(1+r)-\mu k'T/K}{\sigma\sqrt{k'T/K}}\right|\leq \Gamma, \quad \forall k'\geq k\right\}\]
and
\[\bar{\rho}_k(X,\theta_k):=\left\{\begin{array}{cl}\sup_{r_{k+1}\in\U_k^\theta(\theta_k)} X(r_{k+1}) & \mbox{if $\U_k^\theta(\theta_k)\neq\emptyset$}\\ X(0) & \mbox{otherwise}\end{array}\right.\,.\]

In the case of the set $\U_2$ inspired from \cite{bernhard2003robust}, we let instead $\theta_k:=\sum_{\ell=1}^k r_\ell^2$, with $\theta_{k+1} := \theta_{k}+r_{k+1}^2$ and $\bar{\rho}_k(X,\theta_k):=\sup_{r_{k+1}\in\U_k^\theta(\theta_k)} X(r_{k+1})$ where
\[\U_k^\theta(\theta_k):=\left\{ r\in\Re\,\middle|\, \theta_k+r^2 \in [\sigma iNT/K - \Gamma \sqrt{iN}, \sigma iNT/K + \Gamma \sqrt{iN}]\,,\,\forall\,i\geq k/N\right\}\,.\]
The rest of the details are very similar as previously.

\subsection{\modified{Implementation Details Regarding How the Dynamic Program Was Solved}}\label{app:implementationDetails}

{In order to solve the dynamic programs of the models, in the first step, we divide our simulated stock paths into a training, $D_{train}$, and a test set, $D_{test}$. Then we calibrate the uncertainty set parameter $\Gamma$ in a way that 95 percent of the train paths fall into the uncertainty set. Depending on the type of the uncertainty set, $\U_1$ or $\U_2$, we consider either cumulative log returns $\sum_{l=1}^k \log (1+r_l)$, or cumulative square returns $\sum_{l=1}^k r_l^2$ along with the stock price $S_k$ as state variables of the DP. Next, we generate a two dimensional grid for the state variables at each time step. The upper and lower bounds of the grid for the stock prices are obtained by simply considering the price bounds in $D_{train}$. The bounds could be computed in a more conservative way by considering some deviations from the minimum and maximum values of $D_{train}$ so that they contain with a higher probability the paths of the $D_{test}$, however, we did not observe improvements in terms of the hedging performance by considering such bounds in our setting. Another issue worth mentioning is regarding the time dependency of the grid bounds. In our implementation we consider the same bounds for all time periods. This would allow the model to consider the cases where the whole budget of the uncertainty sets are used up in the very first steps. However, in general, the bounds of the grid could be time dependent and determined at each time by using the price paths in that time. For the other state variable, we first use the $D_{train}$ to compute the associated paths of $\sum_{l=1}^k \log (1+r_l)$ or $\sum_{l=1}^k r_l^2$. The upper and lower bounds of the grid could be computed from these values. 

In order to solve the dynamic model, starting from the last period, for each combination of state variables, and for each side of the contract, the \quoteIt{best hedging risk to go} is computed and assigned to the point. For periods other than the last period, we need to solve an optimization model to obtain the optimal allocation of wealth. Since we discretized the state space, at each point we find the reachable values in the next period for the two state variables $(\theta_k,S_k)$. To do so, we start with the reachable values for $r_{k+1}$ using the definition of the projected uncertainty set. Specifically, for $\U_1$ we will have $r_{k+1} \in [\underline{r}_{k+1},\overline{r}_{k+1}]$, while for $\U_2$ we will have $r_{k+1} \in [\underline{r}_{k+1}^1,\overline{r}_{k+1}^1] \cup [\underline{r}_{k+1}^2,\overline{r}_{k+1}^2]$. When choosing the reachable grid points, we always include up to the first grid point that falls outside the intervals in order to induce a conservative bias to our approximation. Having these points, the next step is to find the optimal wealth allocation by solving a piece-wise linear convex optimization problem. We pursue recursively until $k=0$.
}


\end{document}